\documentclass[11pt, reqno]{article}
\pdfoutput=1
\usepackage{amssymb,amsmath,amsthm,amsfonts,latexsym}
\usepackage{epsfig,esint}
\usepackage{mathrsfs}
\usepackage{mathtools,textcomp}
\usepackage{graphicx}
\usepackage{epic,eepic,wrapfig,color, ifthen} 
\usepackage{url,cite}
\usepackage[colorlinks=true, citecolor=blue, filecolor=cyan, linkcolor=red, urlcolor=magenta]{hyperref}
\usepackage{pmboxdraw}
\usepackage{verbatim}
\usepackage[normalem]{ulem}
\usepackage{enumerate}
\usepackage[bottom]{footmisc}

\numberwithin{equation}{section}

\mathtoolsset{showonlyrefs}

\def\pf{\noindent{\bf Proof. }}
\def\qed{{\hfill $\Box$ \bigskip}}

\def\1{{\bf 1}}
\def\nn{\nonumber}
\def\sA {{\cal A}}  \def\sC {{\cal C}}

 \def\sK {{\cal K}}

\def\R {{\mathbb R}} 
\def\N{{\mathbb N}} 
\def\E {{\mathbb E}}  \def \P{{\mathbb P}}

\def\eps{\varepsilon}

\def\wt{\widetilde}

\def\Log{\text{\rm Log}}
\def\Lg{\text{\rm Log}\,}

\theoremstyle{plain}
\newtheorem{thm}{Theorem}[section]
\newtheorem{lem}[thm]{Lemma}
\newtheorem{cor}[thm]{Corollary}
\newtheorem{remark}[thm]{Remark}
\newtheorem{prop}[thm]{Proposition}
\newtheorem{defn}[thm]{Definition}

\theoremstyle{definition}
\newtheorem*{eg*}{Example}
\newtheorem*{egs*}{Examples}
\newtheorem*{def*}{Definition}
\theoremstyle{remark}

\begin{document}
	
	\title{Fractional Laplacian with supercritical killings}	
	\date{}

\author{
	{\bf Soobin Cho}
	\quad {\bf Renming Song\thanks{Research supported in part by a grant from the Simons Foundation (\#960480, Renming Song).}
	}
}

	\maketitle
	
	\begin{abstract}
		In this paper, we study Feynman-Kac semigroups of symmetric $\alpha$-stable processes with supercritical killing potentials belonging to a large class of functions containing functions of the form $b|x|^{-\beta}$,
		where $b>0$ and $\beta>\alpha$. 
		We obtain two-sided estimates on the densities $p(t, x, y)$ of these semigroups for all $t>0$, along with estimates for the corresponding Green functions.
		 
		\medskip
		
		\medskip
		
		
	\end{abstract}
	\allowdisplaybreaks

	\section{Introduction}\label{s:intro}	
	Suppose 
	$Y=\{Y_t, t\ge 0; \, \P_x, x \in \R^d\}$ is a Markov process on 	$\R^d$ and $\kappa$ is a Borel function  on $\R^d$.
Under some conditions on $\kappa$, the operators
$$
P^{\kappa}_tf(x)=\E_x\left[\exp\bigg(-\int^t_0\kappa(Y_s)ds\bigg)f(Y_t)\right], \quad t\ge 0,
$$
form a semigroup. The semigroup $(P^{\kappa}_t)_{t\ge 0}$ is called a Feynman-Kac semigroup or a Schr\"odinger semigroup, 
with potential $\kappa$.

If $L$ is the generator of $Y$, then the generator of $(P^{\kappa}_t)_{t\ge 0}$ is $L-\kappa$. 
Feynman-Kac semigroups are very important in probability, statistical physics and other fields. They have been intensively studied, see, for instance, \cite{CZ, Simon} and the references therein for early contributions, and  \cite{CW23a, CW23b, CKS15, CKSV} and the references therein for recent contributions. 
If $\kappa$ is a non-negative function, then $(P^{\kappa}_t)_{t\ge 0}$ is a (sub)-Markov semigroup, and so there is a Markov process associated with it. In this case, we say that $\kappa$ is a killing potential.

A function $\kappa$ is said to  belong to the Kato class of $Y$ if 
$$
\lim_{t\to0}\sup_{x\in \R^d}\E_x\int^t_0|\kappa(Y_s)|ds=0.
$$
For equivalent analytic characterizations of the Kato class of Brownian motion, 
we refer to   \cite{CZ, Simon}.
A potential belonging to the Kato class is small,  in the sense of quadratic forms,  compared to  the generator $L$ of  $Y$, 
so $L-\kappa$ can be regarded as a small (or subcritical) perturbation of $L$.  
For various Hunt processes, it is well known that Kato class perturbations preserve short time heat kernel estimates, see, for instance, 
\cite{BM,Tak07,Dev19} for the cases of Gaussian heat kernel estimates,  \cite{Song06,Wan08} for $\alpha$-stable-like heat kernel estimates,  and  \cite{KKK23} for more general heat kernel estimates.

In this paper, we will concentrate on Feynman-Kac semigroups with 
killing potentials
belonging to a certain  class, containing functions of the form $\kappa(x)=b|x|^{-\beta}$ with $b, \beta>0$ as examples, in the case when $Y$ is an isotropic $\alpha$-stable process (with generator $-(-\Delta)^{\alpha/2}$), $\alpha\in (0, 2)$. Potentials of the form $\kappa(x)=b|x|^{-\beta}$, particularly the case $\beta=\alpha$, 
are very important and have been studied intensively,
see, for instance, \cite{AMPP16a, BG, BFT, BDK,BGJP, LZ17} and the references therein.

In the remainder of this paper, we always assume that $Y$ is an isotropic $\alpha$-stable process with generator $-(-\Delta)^{\alpha/2}$.
The Dirichlet form of  $Y$ is $(\mathcal{E}, D(\mathcal{E}))$, where
$$
D(\mathcal{E})=\left\{u\in L^2(\R^d): \int_{\R^d}\int_{\R^d}\frac{(u(x)-u(y))^2}{|x-y|^{d+\alpha}}dxdy<\infty\right\}
$$
and
$$
\mathcal{E}(u, u)=\frac12\sA(d, -\alpha)\int_{\R^d}\int_{\R^d}\frac{(u(x)-u(y))^2}{|x-y|^{d+\alpha}}dxdy, \quad u\in D(\mathcal{E})
$$
with 
$$
\sA(d, -\alpha)=\frac{\alpha 2^{\alpha-1}\Gamma((d+\alpha)/2)}{\pi^{d/2}\Gamma(1-\alpha)}.
$$
The 
process $Y$ has a transition density $q(t, x, y)$ such that
\begin{align}\label{e:HKE-alpha-stable}
	C^{-1}	\wt q(t,x,y) \le 	q(t,x,y) \le C	\wt q(t,x,y),  \qquad (t, x, y)\in (0, \infty)\times\R^d\times\R^d
\end{align}
for some $C=C(d, \alpha)>1$, where
\begin{align}\label{e:def-wtq}
	\wt q(t,x,y):=t^{-d/\alpha}\bigg( 1 \wedge \frac{t^{1/\alpha}}{|x-y|}\bigg)^{d+\alpha}.
\end{align}

When $\beta\in (0, \alpha)$, the function $\kappa(x)=b|x|^{-\beta}$ 
with $b\in \R$
belongs to the Kato  class of $Y$.  In this case,
the corresponding Feynman-Kac semigroup $(P^{\kappa}_t)_{t\ge 0}$ has a density $p^{\kappa}(t, x, y)$ and, for small $t$, $p^\kappa(t,x,y)$ is comparable with $q(t, x, y)$. 
More precisely, for any $T > 0$, there exists a constant $C = C(d, \alpha, b, \beta, T) > 1$ such that 
$$
C^{-1}	\wt q(t,x,y) \le 	q(t,x,y) \le C	\wt q(t,x,y), \qquad (t, x, y)\in (0, T]\times\R^d\times\R^d.
$$
See \cite[Theorem 3.4]{Song06}. 

The function $\kappa(x)=b|x|^{-\alpha}$ does not belong to the Kato class of $Y$. According to the fractional Hardy inequality, the quadratic form below
$$
 \mathcal{E}(u, u)+b\int_{\R^d}\frac{u^2(x)}{|x|^\alpha}dx, \quad u\in C^\infty_c(\R^d)
$$
is non-negative definite
 if and only if 
$$
b\ge -b^*_{d,\alpha} \quad \text{where} \;\; b^*_{d,\alpha}:=\frac{2^\alpha \Gamma((d+\alpha)/4))^2}{\Gamma((d-\alpha)/4))^2}.
$$ 
In this sense, $\kappa(x)=b|x|^{-\alpha}$ is a critical potential.
 In this case, the effective state space of the
Markov process associated with $(P^\kappa_t)_{t\ge 0}$ is $\R^d_0:=\R^d\setminus\{0\}$. 
When $d>\alpha$ and $b \ge -b^*_{d,\alpha}$, 
according to \cite[Theorem 1.1]{BGJP}, \cite[Theorem 3.9]{CKSV} and \cite[Theorem 1.1]{JW},
the corresponding 
Feynman-Kac semigroup $(P^\kappa_t)_{t\ge 0}$ has a density $p^\kappa(t, x, y)$ satisfying the following estimates:
there exists $C=C(d, \alpha, b, \alpha)>1$ such that for all $(t, x, y)\in (0, \infty)\times\R^d\times \R^d$,
\begin{equation*}
C^{-1}\left(1\wedge \frac{|x|}{t^{1/\alpha}}\right)^\delta  \left(1\wedge \frac{|y|}{t^{1/\alpha}}\right)^\delta \widetilde{q}(t, x, y)\le p^\kappa(t, x, y)\le C\left(1\wedge \frac{|x|}{t^{1/\alpha}}\right)^\delta  \left(1\wedge \frac{|y|}{t^{1/\alpha}}\right)^\delta \widetilde{q}(t, x, y) ,
\end{equation*}
where
$\delta \in [-(d-\alpha)/2,\alpha)$
 is uniquely determined by the one-to-one map
$$ 
 \delta \mapsto -
\frac{2^\alpha\Gamma((\alpha-\delta)/2)\,\Gamma( (d+\delta)/2)}{\Gamma(-\delta/2)\, \Gamma((d+\delta-\alpha)/2)} = b
$$
from $ [-(d-\alpha)/2,\alpha)$ to $[-b^*_{d,\alpha},\infty)$.

When $\beta>\alpha$, the function 
$\kappa(x)=b|x|^{-\beta}$ 
is said to be supercritical. 
Note that if $b<0$, the quadratic form
$$
\mathcal{E}(u, u)+b\int_{\R^d}\frac{u^2(x)}{|x|^\beta}dx, \quad u\in C^\infty_c(\R^d)
$$
is not non-negative definite. Thus, we will concentrate on the case $b>0$. The purpose of this paper is to obtain two-sided estimates on the density of the corresponding Feynman-Kac semigroup
 In fact, we will deal with a general class of supercritical potentials which includes
$\kappa(x)=b|x|^{-\beta}$, 
with $\beta>\alpha$ and $b>0$,  as examples. 
We state here our main result in the particular case of $\kappa(x)=b|x|^{-\beta}$ 
with $\beta>\alpha$ and $b>0$. We will use
$p^{\beta, b}(t, x, y)$ to denote the density of the corresponding Feynman-Kac 
semigroup and use $G^{\beta,b}(x,y)$ to denote the corresponding Green function
defined as $G^{\beta,b}(x,y) = \int_0^\infty 	p^{\beta,b}(t,x,y)\, dt$.
For $a, b\in \R$, we use the usual notation $a\wedge b:=\min\{a,b\}$ and $a\vee b:=\max\{a,b\}$. 
For two non-negative functions $f$ and $g$, the notation $f\asymp g$ means that there exists a constant $c>1$  such that $c^{-1}g(x) \le f(x) \le cg(x)$ in the common domain of definition of $f$ and $g$. In this paper, we will use the following notation
\begin{equation}\label{e:Log}
\Lg r:= \log(e-1+r), \quad r\ge 0.
\end{equation}

\begin{thm}\label{t:special-case}
(i) {\rm \bf (Small time estimates)}
 Let $T>0$. There exist constants $\lambda_1=\lambda_1(\beta,b)>0$, $\lambda_2=\lambda_2(\beta,b)>0$ and $C= C(d,\alpha,\beta,b,T)>1$ such that for all $t \in (0,T]$ and $x,y \in \R^d_0$,
	\begin{align}\label{e:specil-case-upper}
	p^{\beta,b}(t,x,y) &\le C\bigg(1 \wedge \frac{|x|^\beta}{t}\bigg) \bigg(1 \wedge \frac{|y|^\beta}{t}\bigg)\nn\\
	&\quad \times  \bigg[\, e^{-\lambda_1t/(|x|\vee |y|)^\beta}t^{-d/\alpha}\bigg( 1 \wedge \frac{t^{1/\alpha}}{|x-y|}\bigg)^{d+\alpha} +  t^{-(d+2\alpha-2\beta)/\beta}\bigg(1 \wedge \frac{t^{1/\beta}}{|x-y|}\bigg)^{d+2\alpha}\, \bigg] \qquad 
\end{align}
and
\begin{align}\label{e:specil-case-lower}
	p^{\beta,b}(t,x,y) &\ge C^{-1}\bigg(1 \wedge \frac{|x|^\beta}{t}\bigg) \bigg(1 \wedge \frac{|y|^\beta}{t}\bigg)\nn\\
	&\quad \times  \bigg[\, e^{-\lambda_2t/(|x|\vee |y|)^\beta}t^{-d/\alpha}\bigg( 1 \wedge \frac{t^{1/\alpha}}{|x-y|}\bigg)^{d+\alpha} +   t^{-(d+2\alpha-2\beta)/\beta}\bigg(1 \wedge \frac{t^{1/\beta}}{|x-y|}\bigg)^{d+2\alpha}\, \bigg]. \qquad 
\end{align}

\noindent (ii) {\rm \bf(Large time estimates)}   There exist comparison constants depending only on $d,\alpha,\beta$ and $b$  such that the following estimates hold for all $t \in [2,\infty)$ and $x,y \in \R^d_0$: (1) If $d>\alpha$, then
\begin{align*}
	p^{\beta,b}(t,x,y) \asymp (1\wedge |x|^\beta) (1\wedge |y|^\beta)\,t^{-d/\alpha}\bigg( 1 \wedge \frac{t^{1/\alpha}}{|x-y|}\bigg)^{d+\alpha};
\end{align*}
(2) if $d=1<\alpha$, then
	\begin{align*}
		p^{\beta,b}(t,x,y) \asymp \bigg( 1 \wedge \frac{|x|^{\beta} \wedge |x|^{\alpha-1} }{t^{(\alpha-1)/\alpha}} \bigg) \bigg( 1 \wedge \frac{|y|^\beta \wedge |y|^{\alpha-1} }{t^{(\alpha-1)/\alpha}} \bigg) \,  t^{-1/\alpha}\bigg( 1 \wedge \frac{t^{1/\alpha}}{|x-y|}\bigg)^{1+\alpha};
\end{align*}
and (3) if $d=1=\alpha$, then 
	\begin{align*}
p^{\beta,b}(t,x,y) \asymp \bigg( 1 \wedge \frac{|x|^\beta \wedge \Lg |x| }{\Lg t} \bigg) \bigg( 1 \wedge \frac{|y|^\beta \wedge \Lg |y| }{\Lg t} \bigg) \, t^{-1}\bigg( 1 \wedge  \frac{t}{|x-y|}\bigg)^{2}.
\end{align*}

\noindent (iii) {\rm \bf(Green function estimate)} There exist comparison constants depending only on $d,\alpha,\beta$ and $b$ such that for all  $x,y \in \R^d_0$, (1) if $d>\alpha$, then
	\begin{align*}
	G^{\beta,b}(x,y) \asymp \bigg( 1\wedge \frac{(|x|\wedge 1)^\beta}{(|x-y|\wedge 1)^{\alpha}}\bigg) \bigg( 1\wedge \frac{(|y| \wedge 1)^\beta}{(|x-y|\wedge 1)^{\alpha}}\bigg) \frac{1}{|x-y|^{d-\alpha}};
\end{align*}
(2) if $d=1<\alpha$, then
\begin{align*}
	G^{\beta,b}(x,y) &\asymp 	\left( 1\wedge \frac{(|x|\wedge 1)^\beta}{(|x-y|\wedge 1)^\alpha}\right) 	\left( 1\wedge \frac{(|y|\wedge 1)^\beta}{(|x-y|\wedge 1)^\alpha}\right) 	\left( 1\wedge \frac{|x| \vee 1}{|x-y| \vee 1}\right)^{\alpha-1} 	\left( 1\wedge \frac{|y| \vee 1}{|x-y| \vee 1}\right)^{\alpha-1} \nn\\
	&\;\;\quad \times \big( (|x|^\beta \wedge  |x|^\alpha) \vee |x-y|^\alpha \big)^{(\alpha-1)/(2\alpha)}  \big( (|y|^\beta \wedge  |y|^\alpha) \vee |x-y|^\alpha \big)^{(\alpha-1)/(2\alpha)};
\end{align*}
and (3) if $d=1=\alpha$, then
\begin{align*}
	G^{\beta,b}(x,y) &\asymp 	\left( 1\wedge \frac{(|x|\wedge 1)^\beta}{|x-y|\wedge 1}\right) 	\left( 1\wedge \frac{(|y|\wedge 1)^\beta}{|x-y|\wedge 1}\right) 	\left( 1\wedge \frac{\Lg |x|}{\Lg |x-y|}\right)^{1/2} \left( 1\wedge \frac{\Lg |y|}{\Lg |x-y|}\right)^{1/2}  \nn\\
	&\;\;\quad \times \left[   \Log \left( \frac{|x|^\beta \wedge |x| }{|x-y|\wedge 1}\right)  \Log \left( \frac{|y|^\beta \wedge |y| }{|x-y| \wedge 1}\right) \right]^{1/2}.
\end{align*}
\end{thm}

\begin{remark}
\rm Both terms in the square brackets on the right hand sides of  \eqref{e:specil-case-upper} and \eqref{e:specil-case-lower}
are needed. In some region of $(0, T]\times\R^d\times \R^d$, one term dominates and otherwise the other term dominates.
\end{remark}

It is natural to study heat kernel estimates of Feynman-Kac semigroups of Brownian motion
with supercritical killing potentials of the $\kappa(x)=b|x|^{-\beta}$ with $\beta>2$. We believe
that the corresponding heat kernel estimates are  drastically different.
One reason for this is that,  according
to \cite{MS04, MS05, MSS17},
 the density $p^{2, b}(t, x, y)$ of the Feynman-Kac semigroup corresponding to the generator $\Delta-b|x|^{-
	2}$ admits the following estimates: There exist positive constants $c_i, i=1, 2, 3, 4,$ such that for all $(t, x, y)\in (0, \infty)\times \R^d\times \R^d$,
\begin{align*}
&c_1\left(1\wedge\frac{|x|}{\sqrt{t}}\right)^\delta \left(1\wedge\frac{|y|}{\sqrt{t}}\right)^\delta t^{-d/2}\exp\left(-c_2\frac{|x-y|^2}{t}\right)\\
&\le p^{2, b}(t, x, y)
\le c_3\left(1\wedge\frac{|x|}{\sqrt{t}}\right)^\delta \left(1\wedge\frac{|y|}{\sqrt{t}}\right)^\delta t^{-d/2}\exp\left(-c_4\frac{|x-y|^2}{t}\right),
\end{align*}
where $\delta = \sqrt{(d-2)^2/4+ b } - (d-2)/2$. 
Another reason is that, according to \cite[Lemma 2.1]{LZ17}, if $\beta>2$, then
$$
u(x)=|x|^{-\frac{d}2+1}K_{(d-2)/(\beta-2)}\left(\frac{2}{\beta-2}\sqrt{b}|x|^{-\frac{\beta-2}2}\right)
$$
is a solution of
$$
\Delta u-b|x|^{-\beta} u=0,
$$
where $K_{(d-2)/(\beta-2)}$ is the modified   Bessel function of the second kind. It is known that 
$$
\lim_{r\to \infty} K_{(d-2)/(\beta-2)} (r) / (r^{-1/2}e^{-r}) = (\pi/2)^{1/2}.
$$ 
See \cite[9.7.2]{AS64}.  
Using this asymptotic property of $K_{(d-2)/(\beta-2)}$, one can easily check that the function $u$ decays exponentially at the origin.
Due to the two reasons above, we believe that the density of the Feynman-Kac semigroup corresponding to the generator 
$\Delta-b|x|^{-\beta}$, $\beta>2$, decays to 0 exponentially at the origin, instead of algebraically as in
Theorem \ref{t:special-case}. Because of this, heat kernel estimates of Feynman-Kac semigroups of Brownian motion
with supercritical killing potentials of the $\kappa(x)=b|x|^{-\beta}$ with $\beta>2$ are more delicate. We intend to tackle this
in a separate project.

The approach of this paper is probabilistic. The probabilistic representation of the Feynman-Kac semigroup, the L\'evy
system formula, and the sharp two-sided Dirichlet heat estimates for the fractional Laplacian in balls and exterior balls
obtained in \cite{BGR10, CKS10} play essential roles. 
To help us  get the behavior of the heat kernel near the origin, we need to construct appropriate barrier functions,
see Section \ref{s:3}.

Now we introduce the class of potentials that we will deal with in this paper.

	For $\beta_2\ge \beta_1>\alpha$ and  $\Lambda\ge 1$, let 	 $\sC_\alpha(\beta_1,\beta_2,\Lambda)$
	be the family of all strictly increasing continuous functions $\psi:(0,\infty)\to (0,\infty)$ with $\psi(1)=1$ satisfying the following  property:
	\begin{align}\label{e:psi-scale}
	\Lambda^{-1}\bigg(\frac{R}{r}\bigg)^{\beta_1}\le 	\frac{\psi(R)}{\psi(r)} \le 	\Lambda\bigg(\frac{R}{r}\bigg)^{\beta_2} 
	\quad \text{for all} \;\, 0<r\le R.
	\end{align}

\begin{defn}
\rm	 Let  $\psi \in \sC_\alpha(\beta_1,\beta_2,\Lambda)$
and   $\kappa$ be a non-negative Borel function on $\R^d_0$.

\smallskip

\noindent(i) We say that $\kappa$ belongs to the class  $\sK^0_\alpha(\psi,\Lambda)$ if
\begin{align}\label{e:kappa-cond-1}
	\frac{\Lambda^{-1}}{ \psi(|x|)}\le 	\kappa(x) \le 	\frac{\Lambda}{ \psi(|x|)}  \quad \text{for all} \;\, x \in \R^d_0 \text{ with } |x|\le 1
\end{align}
and
\begin{align}\label{e:kappa-cond-2}
	\sup_{x\in \R^d \setminus B(0,1)}	\kappa(x) \le \Lambda.
\end{align}

\noindent (ii) 	 We say that $\kappa$ belongs to the class  $\sK_{\alpha}(\psi,\Lambda)$ if  \eqref{e:kappa-cond-1} holds and 
\begin{align}\label{e:kappa-cond-3}
	\kappa(x) \le 	\frac{\Lambda}{\psi(|x|)}  \quad \text{for all} \;\, x \in \R^d_0 \text{ with } |x|\ge 1.
\end{align}
\end{defn}

The inclusion $\sK_{\alpha}(\psi,\Lambda) \subset \sK^0_\alpha(\psi,\Lambda)$ is obvious.

	 In this paper, we consider heat kernel estimates for non-local operators of the form
	\begin{align}\label{e:def-generator}
	L_{\alpha}^{\kappa} := -(-\Delta)^{\alpha/2} - \kappa(x)
	\end{align}
when $\kappa\in \sK^0_\alpha(\psi,\Lambda)$ or  $\kappa\in\sK_{\alpha}(\psi,\Lambda)$.
It is obvious that the function $\kappa(x)=b|x|^{-\beta}$, with $\beta>\alpha$ and $b>0$, belongs to $\sK_{\alpha}(r^\beta,b \vee b^{-1})$.

The rest of this paper is organized as follows. In Section \ref{s:2}, we collect some preliminary results. In Section \ref{s:3}, we construct appropriate barrier functions and prove some survival probability estimates. Small time heat kernel estimates are proved in Section \ref{s:4}.  In Section \ref{s:5}, we prove a key proposition needed to get large time heat kernel estimates.
Large time heat kernel estimates are proved in Section \ref{s:6}. In Section \ref{s:7}, we prove the Green function estimates. In the Appendix, we prove a lemma which is used in getting the large time heat kernel upper bound in the case $d=1=\alpha$. We believe this result is of independent interest.

\section{Preliminaries}\label{s:2} 
In the remainder of this paper, 
we assume that $\beta_2\ge \beta_1>\alpha$  and $\Lambda\ge 1$ are given
constants,  $\psi \in  \sC_\alpha(\beta_1,\beta_2,\Lambda)$ and $\kappa \in \sK^0_\alpha(\psi,\Lambda)$.

Recall that $Y=\{Y_t, t\ge 0; \, \P_x, x \in \R^d\}$ is an isotropic $\alpha$-stable process on $\R^d$ with generator $ -(-\Delta)^{\alpha/2}$ and the density $q(t, x, y)$ of $Y$ satisfies \eqref{e:HKE-alpha-stable} with $\wt{q}$ defined by
\eqref{e:def-wtq}. 
From \eqref{e:HKE-alpha-stable},  one sees that there exists $C=C(d,\alpha)> 1$ such that for all $t,s>0$ and $x,y \in \R^d$,
\begin{equation}\label{e:wtq-integral}
	\int_{\R^d} \wt q(t,x,z) dz \le C \int_{\R^d} q(t,x,z) dz  =C
\end{equation}
and
\begin{equation}\label{e:wtq-semigroup}
	\int_{\R^d} \wt q(t,x,z) \wt q(s,z,y) dz \le  C^2 	\int_{\R^d}  q(t,x,z)  q(s,z,y) dz =C^2 q(t+s,x,y) \le C^3\wt q(t+s,x,y).
\end{equation}
Note that
\begin{align}\label{e:wtq-compare}
	\wt q(t,x,y) \asymp \wt q(t,0,y) \quad \text{for all $t>0$ and $x,y \in \R^d$ with $|x|\le 2t^{1/\alpha}$}.
\end{align}
Indeed, if $|y|<4t^{1/\alpha}$, then $\wt q(t,x,y) \asymp t^{-d/\alpha} \asymp \wt q(t,0,y)$ and if $|y|\ge 4t^{1/\alpha}$, then $\wt q(t,x,y) = t|y-x|^{-d-\alpha} \asymp t|y|^{-d-\alpha} \asymp \wt q(t,0,y)$. We also note that
\begin{align}\label{e:wtq-compare-2}
	2^{-1}\wt q(t,x,y)\le 	\wt q(s,x,y) \le 2^{d/\alpha}\wt q(t,x,y) \quad \text{for all $t>0$, $t/2\le s\le t$ and $x,y \in \R^d$}.
	\end{align}
For an open subset $U$ of $\R^d$, let $\tau^Y_U:=\inf\{t>0:\, Y_t\notin U\}$ be  the first exit time from $U$ for $Y$. The killed process $Y^U$ on $U$ is defined by $Y^U_t=Y_t$ if $t<\tau^Y_U$ and $Y^U_t=\partial$ if $t\ge \tau^Y_U$, where $\partial$ is the cemetery point. It is known that the process $Y^U$ has a jointly continuous transition density $q^U:(0,\infty)\times U \times U \to (0,\infty)$ such that $q^U(t,x,y)\le q(t,x,y)$ for all $t>0$ and $x,y \in U$.

Define for  $t\ge 0$, $x \in \R^d$ and a non-negative Borel function $f$ on $\R^d$,
\begin{align}\label{e:Feynman-Kac-0}
	P^{\kappa}_t f(x) =  \E_x \left[ \,\exp \bigg(-  \int_0^t   \kappa(Y_s)ds \bigg)   f(Y_t) \,\right].
\end{align}
Since $\kappa$ is non-negative, we can extend $(P^\kappa_t)_{t\ge 0}$ to a strongly continuous Markov semigroup on  $L^2(\R^d)$.
Let $X^\kappa$ be the Markov process on $\R^d$ associated with the semigroup $(P^\kappa_t)_{t\ge 0}$.  
When $d\ge \alpha$, singletons are polar, and so 
$\P_x(\tau^Y_{\R^d_0}=\infty)=1$.
 So in this case, we can regard $X^\kappa$ as a Markov process
on $\R^d_0$. 
Now let us consider the case $d=1<\alpha$. In this case, since $\beta_1>\alpha>1$,
 by \eqref{e:kappa-cond-1} and \eqref{e:psi-scale} (with $\psi(1)=1$),  we have for all $\delta\in (0,1)$,
$$
\int_{-\delta}^\delta \kappa(y)dy  \ge \frac{1}{\Lambda} \int_{-\delta}^\delta \frac{dy}{\psi(|y|)} \ge \frac{1}{\Lambda^2} \int_{-\delta}^\delta \frac{dy}{|y|^{\beta_1}}=\infty.
$$ 
Hence,  it follows from \cite[Lemma 1.6]{Zanzotto97} that
$$
\P_0\left(\int^t_0\kappa(Y_s)ds=\infty \mbox{ for all } t>0\right)=1.
$$
Combining the above with the strong Markov property, we get that for any $x\in \R_0$, 
$$
\P_x\left(\int^t_0\kappa(Y_s)ds=\infty \mbox{ for all } t>
\tau^Y_{\R^d_0}\right)=1.
$$
Thus, in both cases, it holds that for any $t>0$ and $x\in\R^d_0$,  
$$
 P^\kappa_tf(x)=\E_x\left[\exp\left(-\int^{t}_0\kappa(Y_s)ds\right)f(Y_t) : t<\tau^Y_{\R^d_0}\right]=\E_x\left[\exp\left(-\int^{t}_0\kappa(Y^{\R^d_0}_s)ds\right)f(Y^{\R^d_0}_t) \right].
$$
Therefore, $X^\kappa$ can be  considered as a process on $\R^d_0$ and it can be obtained as follows:
first kill the process $Y$ upon exiting $\R^d_0$ and then kill the resulting process using the killing potential $\kappa$.
From now on, we always regard $X^\kappa$ as a process on $\R^d_0$.

For an open subset $U$ of $\R^d_0$, we denote by $\tau^{\kappa}_U=\tau^{X^\kappa}_U$ the first exit time from $U$ for $X^\kappa$ and by $X^{\kappa,U}$ the killed process of $X^\kappa$ on $U$. 
Let $U$ be an open subset of $\R^d_0$ such that $\overline U \subset \R^d_0$.  By \eqref{e:kappa-cond-1} and \eqref{e:kappa-cond-2}, we have  $\sup_{x\in U} \kappa(x)<\infty$, and therefore,
$$
\lim_{t\to 0}\sup_{x\in U}\left| \int^t_0\int_U q^U(s, x, y) \kappa(y)  dyds\right| \le \lim_{t\to 0} \, t\sup_{y\in U} \kappa(y)=0.
$$
Hence, by general theory, the semigroup $(P^{\kappa,U}_t)_{t\ge 0}$ of $X^{\kappa,U}$ can be represented by the following Feynman-Kac formula: For any $t\ge 0$, $x \in U$ and any non-negative Borel function $f$ on $U$,
\begin{align}\label{e:Feynman-Kac}
	P^{\kappa,U}_t f(x) &=  \E_x \left[ \,\exp \bigg(-  \int_0^t   \kappa(Y^U_s)ds \bigg)  f(Y^U_t) \,\right].
\end{align}
We refer to \cite[Section 1.2]{CKS15} for more details. Define $p^{\kappa,U}_0(t,x,y) := q^U(t, x, y)
$ and, for $k \ge 1$, 
\begin{align}\label{e:p-construction}
	p^{\kappa,U}_{k}(t,x,y) =&-  \int_0^t \int_U q^U(s, x, z)
	p^{\kappa,U}_{k-1}(t-s,z,y) \kappa(z)  dzds.
\end{align}
Set $p^{\kappa,U}(t,x,y) := \sum_{k=0}^\infty p^{\kappa,U}_{k}(t,x,y).$ According to  \cite[Theorem 3.4]{CKS15},   $p^{\kappa,U}(t,x,y) $ is jointly continuous on $(0,\infty) \times U \times U$ and is the transition density for $P^{\kappa,U}_t$.  Define  for $t>0$ and $x,y \in \R^d_0$,
$$
p^\kappa(t,x,y) :=\lim_{n \to \infty} p^{\kappa, B(0,1/n)^c} (t,x,y) .
$$
Using the monotone convergence theorem,  one sees that $p^\kappa(t,x,y)$ is the transition density of $X^\kappa$. 
Consequently, we have $p^\kappa(t,x,y)=p^\kappa(t,y,x) \le q(t,x,y)$ for all $t>0$ and $x,y \in \R^d_0$. Moreover, as an increasing limit of continuous functions, for every fixed $t>0$ and $z \in \R^d_0$, the maps $x\mapsto p^\kappa(t,x,z)$ and $y \mapsto p^\kappa(t,z,y)$ are lower semi-continuous.

\subsection{Some properties of  $Y$}

The following result is well known, see \cite[Theorem 5.1]{BSW13} for the corresponding result for L\'evy-type processes.

\begin{prop}\label{p:EP-alpha-stable}
There exists  $C=C(d,\alpha)>0$ such that for all $x\in \R^d$, $t>0$ and $R>0$,
\begin{align*}
	\P_x(\tau^Y_{B(x,R)}\le t) \le Ct R^{-\alpha}.
\end{align*}
\end{prop}

For an open set $E\subset \R^d$ and $x \in E$, we let 
$$\delta_E(x):=\inf\{|x-y|:y \in E^c\}.$$ 
\begin{prop}\label{p:DHKE-1}  
For any $k \ge 1$, there   exist comparison constants depending only on $d,\alpha$ and $k$ such that for any $z \in \R^d$, $R>0$, $0<t\le (kR)^\alpha$ and $x,y \in B(z,R)$,
	\begin{align*}
		\frac{q^{B(z,R)}(t,x,y)}{ \wt q(t,x,y)} \asymp \bigg( 1 \wedge \frac{\delta_{B(z,R)}(x)^{\alpha/2}}{t^{1/2}}\bigg)\bigg( 1 \wedge \frac{\delta_{B(z,R)}(y)^{\alpha/2}}{t^{1/2}}\bigg).
	\end{align*}
\end{prop}
\pf Without loss of generality, by the translation-invariance of $Y$, we assume that $z=0$. By the scaling property of $Y$, we have
\begin{align}\label{e:DHKE-scaling}
		q^{B(0,R)}(t,x,y)= R^{-d} q^{B(0,1)}(t/R^\alpha,x/R,y/R).
\end{align}
Applying \cite[Thoerem 1.1(i)]{CKS10} with $T=k^\alpha$, since  $\delta_{B(0,1)}(w/R)=\delta_{B(0,R)}(w)/R$ for all $w \in B(0,R)$, we obtain
\begin{align*}
	&q^{B(0,1)}(t/R^\alpha,x/R,y/R) \\
	& \asymp \bigg( 1 \wedge \frac{\delta_{B(0,1)}(x/R)^{\alpha/2}}{(t/R^\alpha)^{1/2}}\bigg)\bigg( 1 \wedge \frac{\delta_{B(0,1)}(y/R)^{\alpha/2}}{(t/R^\alpha)^{1/2}}\bigg) (t/R^\alpha)^{-d/\alpha}  \bigg( 1 \wedge \frac{(t/R^\alpha)^{1/\alpha}}{|x-y|/R}\bigg)^{d+\alpha}\\
		&=  R^d \bigg( 1 \wedge \frac{\delta_{B(z,R)}(x)^{\alpha/2}}{t^{1/2}}\bigg)\bigg( 1 \wedge \frac{\delta_{B(z,R)}(y)^{\alpha/2}}{t^{1/2}}\bigg) \,\wt q(t,x,y),
\end{align*}
where the comparison constants above depend only on $d,\alpha$ and $k$. Combining this with \eqref{e:DHKE-scaling}, we get the desired result. 
\qed

As a consequence of Proposition \ref{p:DHKE-1}, we obtain

\begin{cor}\label{c:NDL-alpha-stable}
	For any $k \ge 1$,	there exists  $C=C(d,\alpha,k)>0$  such that for all  $x \in \R^d$, $R>0$,  $0<t\le (kR)^\alpha$ and $y \in B(x, t^{1/\alpha}/(2k))$,
	\begin{align*}
		q^{B(x,R)} (t,x,y) \ge Ct^{-d/\alpha}.
	\end{align*}
\end{cor}

 Similar to Proposition \ref{p:DHKE-1}, using  the translation-invariance and the scaling property of $Y$,  we get the following results from \cite[Theorem 3 and Corollary 2]{BGR10}.
\begin{prop}\label{p:DHKE-2}
	There  are comparability constants depending only on $d$ and $\alpha$ such that the following estimates hold for all $z \in \R^d$, $R>0$, $t>0$ and $x,y \in B(z,R)^c$.

	\smallskip
	
	\noindent 	  (i) If $d>\alpha$, then
	\begin{align*}
		\frac{q^{B(z,R)^c}(t,x,y)}{	\wt q(t,x,y)} \asymp \bigg( 1 \wedge \frac{(\delta_{B(z,R)^c}(x)\wedge R)^{\alpha/2}}{(t\wedge R^\alpha)^{1/2}}\bigg)  \bigg( 1 \wedge \frac{(\delta_{B(z,R)^c}(y)\wedge R)^{\alpha/2}}{(t\wedge R^\alpha)^{1/2}}\bigg).
	\end{align*}
	\noindent (ii) If $d=1<\alpha$, then
	\begin{align*}
		\frac{q^{B(z,R)^c}(t,x,y)}{	\wt q(t,x,y)} &\asymp \bigg( 1 \wedge \frac{\delta_{B(z,R)^c}(x)^{\alpha-1} ( \delta_{B(z,R)^c}(x) \wedge R)^{(2-\alpha)/2}}{  t^{(\alpha-1)/\alpha} (t\wedge R^\alpha)^{(2-\alpha)/(2\alpha)} }\bigg) \\
		&\qquad \times \bigg( 1 \wedge \frac{\delta_{B(z,R)^c}(y)^{\alpha-1} ( \delta_{B(z,R)^c}(y) \wedge R)^{(2-\alpha)/2}}{  t^{(\alpha-1)/\alpha} (t\wedge R^\alpha)^{(2-\alpha)/(2\alpha)} }\bigg).
	\end{align*}
	\noindent (iii) If $d=1=\alpha$, then
	\begin{align*}
		\frac{q^{B(z,R)^c}(t,x,y)}{	\wt q(t,x,y)} &\asymp  \bigg( 1 \wedge \frac{(\delta_{B(z,R)^c}(x) \wedge R)^{1/2} \,\Lg ( \delta_{B(z,R)^c}(x)/R)}{(t \wedge R)^{1/2}\, \Lg ( t/R)}\bigg) \\
		&\qquad \times  \bigg( 1 \wedge \frac{(\delta_{B(z,R)^c}(y) \wedge R)^{1/2} \,\Lg (\delta_{B(z,R)^c}(y)/R)}{(t \wedge R)^{1/2}\, \Lg ( t/R)}\bigg).
	\end{align*}
\end{prop}
\pf Since the proofs are similar, we only present the proof for (iii).
 Suppose that $d=1=\alpha$. Without loss of generality, by the translation-invariance of $Y$, we assume that $z=0$. Using the scaling property of $Y$ in the first line below and \cite[Corollary 2]{BGR10} in the second, since  $\delta_{B(0,1)^c}(w/R)=\delta_{B(0,R)^c}(w)/R$ for all $w \in B(0,R)^c$,  we get that  for all  $R>0$, $t>0$ and $x,y \in B(0,R)^c$,
\begin{align}\label{e:DHKE2}
		&q^{B(0,R)^c}(t,x,y) =  	R^{-1}q^{B(0,1)^c}(t/R,x/R,y/R)  \nn\\
		&\asymp \bigg( 1 \wedge \frac{ \log (1+ (\delta_{B(0,R)^c}(x)/R)^{1/2})}{ \log (1+ (t/R)^{1/2})}\bigg) \bigg( 1 \wedge \frac{ \log (1+ (\delta_{B(0,R)^c}(y)/R)^{1/2})}{ \log (1+ (t/R)^{1/2})}\bigg) \wt q(t,x,y).
\end{align}
Note that for all $a>0$ and $r>0$,
\begin{align*}
	\log (1+ (a/r)^{1/2}) &\asymp  \begin{cases}
		(a/r)^{1/2} &\mbox{ if } a/r \le 1,\\
			\Lg (a/r) &\mbox{ if } a/r > 1
	\end{cases}\\
	&\asymp r^{-1/2}(a \wedge r)^{1/2} \Lg (a/r).
\end{align*}
Hence, for all $R>0$, $t>0$ and $w\in B(0,R)^c$, we have 
\begin{align*}
	 1 \wedge \frac{ \log (1+ (\delta_{B(0,R)^c}(w)/R)^{1/2})}{ \log (1+ (t/R)^{1/2})} \asymp 1 \wedge \frac{(\delta_{B(0,R)^c}(w) \wedge R)^{1/2} \, \Lg (\delta_{B(0,R)^c}(w)/R)}{(t \wedge R)^{1/2} \,\Lg ( t/R)}.
\end{align*}
Combining this with \eqref{e:DHKE2}, we arrive at the result. \qed

\subsection{Preliminary estimates for $X^\kappa$}

 For $x\in \R^d_0$ and a Borel subset $A$ of $\R^d_0 \cup\{\partial\}$,  define 
 $$N(x, A)=\sA(d, -\alpha)\int_{A\cap \R^d_0}|x-y|^{-d-\alpha}dy+\kappa(x)1_A(\partial).$$ Then $(N, t)$ is a L\'evy system for $X^\kappa$ (cf. \cite[Theorem 5.3.1]{FOT} and the argument on  \cite[p. 40]{CK03}), that is, 
 for any $x \in  \R^d_0$, non-negative Borel function $f$ on  $\R^d_0  \times(\R^d_0 \cup\{\partial\})$
 vanishing on  $\{(x,x): x \in \R^d_0\} \cup \{(x,\partial):x\in \R^d_0\}$,
 and any stopping time $\tau$,
	\begin{align}\label{e:Levysystem}
		\E_x\bigg[  \sum_{s \le \tau} f( X^\kappa_{s-},  X^\kappa_s)  \bigg] = \E_x \bigg[ \int_0^\tau  \int_{\R^d_0} 
		f( X^\kappa_s, y)N(X^\kappa_s, dy)
		 \,ds  \bigg].
	\end{align}

\begin{prop}\label{p:interior-HKE-lower}
For any  $a\ge 1$, there exists  $C=C(d,\alpha, \Lambda,a)>0$ such that for all $R>0$, $0<r \le 1 \wedge (R/4)$ and   $0<t\le a\psi(r)$,
	\begin{align*}
		p^{\kappa,B(0, R) \setminus B(0,r)}(t,x,y) \ge C	\wt q(t,x,y) \quad \text{for all} \;\,  x,y \in B(0, R/2) \setminus B(0,2r).
	\end{align*}
\end{prop}
\pf By \eqref{e:kappa-cond-1} and \eqref{e:kappa-cond-2},  we have $	\sup_{z \in B(0,R)\setminus B(0,r)} \kappa(z) \le\Lambda/\psi(r)$.
Hence, from  \eqref{e:Feynman-Kac},  we see that for all   $0<t\le a\psi(r)$ 
and $x,y \in B(0, R/2) \setminus B(0,2r)$,
\begin{equation}\label{e:interior-HKE-lower-1}
	p^{\kappa,B(0, R) \setminus B(0,r)}(t,x,y)  \ge  e^{-\Lambda t/\psi(r)}q^{B(0, R) \setminus B(0,r)}(t,x,y)  \ge  
	e^{-\Lambda a}q^{B(0, R) \setminus B(0,r)}(t,x,y).
\end{equation}
Set $r_t:=t^{1/\alpha}/(2\Lambda^{1/\alpha}a^{1/\alpha})$. We consider  two separate cases.

\smallskip

Case 1: 
$|x-y|< r_t$. Note that by \eqref{e:psi-scale}, since $\psi(1)=1$, $r\le 1$ and $\beta_1>\alpha$,
$$
t^{1/\alpha} \le   a^{1/\alpha}  \psi(r)^{1/\alpha} \le 
\Lambda^{1/\alpha} a^{1/\alpha} \psi(1)^{1/\alpha} r^{\beta_1/\alpha}  \le  \Lambda^{1/\alpha} a^{1/\alpha} r .
$$ 
 Hence,  since $B(x,r)\subset B(0,R)\setminus B(0,r)$, using Corollary \ref{c:NDL-alpha-stable}  (with $k=\Lambda^{1/\alpha}a^{1/\alpha}$),  we obtain
\begin{align}\label{e:interior-HKE-lower-2}
	q^{B(0, R) \setminus B(0,r)}(t,x,y) \ge q^{B(x, r)}(t,x,y) \ge c_1 t^{-d/\alpha}.
\end{align}

Case 2: 
 $|x-y|\ge r_t$. Using the strong Markov property, \eqref{e:interior-HKE-lower-2} and \eqref{e:Levysystem}, we get that
\begin{align}\label{e:interior-HKE-lower-3}
	&q^{B(0, R) \setminus B(0,r)}(t,x,y) 	\ge \E_x \bigg[ q^{B(0, R) \setminus B(0,r)}(t-\tau^Y_{B(x,r_t/2)},Y_{\tau^Y_{B(x,r_t/2)}},y): \nn\\
	&\hspace{1.83in} \tau^Y_{B(x,r_t/2)}<t/2, \, Y_{\tau^Y_{B(x,r_t/2)}} \in B(y,(t/2)^{1/\alpha}/ 
	(4\Lambda^{1/\alpha}a^{1/\alpha}))  \bigg] \nn\\
	&\ge c_1 t^{-d/\alpha} \P_x \left( \tau^Y_{B(x,r_t/2)}<t/2, \, Y_{\tau^Y_{B(x,r_t/2)}} \in B(y,(t/2)^{1/\alpha}/ 
	(4\Lambda^{1/\alpha}a^{1/\alpha}))  \right) \nn\\
	&= c_1 t^{-d/\alpha} \E_x \left[\int_0^{\tau^Y_{B(x,r_t/2)} \wedge (t/2)} \int_{B(y,(t/2)^{1/\alpha}/ 
	(4\Lambda^{1/\alpha}a^{1/\alpha}))}	\frac{\sA(d, -\alpha)}{ |Y_s - w|^{d+\alpha}} \,dw ds  \right]. 
\end{align}
For all $z \in B(x,r_t/2)$ and $w \in B(y,(t/2)^{1/\alpha}/ (4\Lambda^{1/\alpha}a^{1/\alpha}))$, we have 
\begin{align}\label{e:interior-HKE-lower-4}
	|z-w| < |x-y|  + r_t/2 + r_t/2 \le 2|x-y|.
\end{align}
Besides, by  Proposition \ref{p:EP-alpha-stable}, there exists $\eps=\eps(d,\alpha)\in (0,1)$ such that
\begin{align*}
\P_x (\tau^Y_{B(x,r_t/2)} \ge  (\eps r_t/2)^\alpha ) = 1-\P_x (\tau^Y_{B(x,r_t/2)} <  (\eps r_t/2)^\alpha ) \ge 1 - c_2 \eps^{1/\alpha} \ge 1/2.
\end{align*} 
It follows that 
\begin{align}\label{e:interior-HKE-lower-5}
	\E_x [\tau^Y_{B(x,r_t/2)} \wedge (t/2) ]	&\ge ((\eps r_t/2)^\alpha \wedge (t/2)) \, \P_x (\tau^Y_{B(x,r_t/2)} \ge (\eps r_t/2)^\alpha ) \nn\\
	&\ge 2^{-1} \left( ( 2^{-2\alpha} \Lambda^{-1}a^{-1}\eps^\alpha) \wedge 2^{-1} \right) t.
\end{align}
Using \eqref{e:interior-HKE-lower-4} and \eqref{e:interior-HKE-lower-5}, we get from  \eqref{e:interior-HKE-lower-3} that
\begin{align}\label{e:interior-HKE-lower-6}
	q^{B(0, R) \setminus B(0,r)}(t,x,y) &\ge \frac{c_1(t/2)^{-d/\alpha} \E_x [\tau^Y_{B(x,r_t/2)} \wedge (t/2) ] }{(2|x-y|)^{d+\alpha}}  \int_{B(y,(t/2)^{1/\alpha}/ 
	(4\Lambda^{1/\alpha}a^{1/\alpha}))} dw\nn\\
	& \ge c_3 t |x-y|^{-d-\alpha}.
\end{align}

Now, combining \eqref{e:interior-HKE-lower-1} with \eqref{e:interior-HKE-lower-2} and \eqref{e:interior-HKE-lower-6},  we arrive at the result.
\qed

\begin{lem}\label{l:exit-time-lower-0}
There exists  $C=C(d,\alpha,\Lambda)>0$ such that for all $x \in \R^d_0$ with $|x|\le 1$, $ 0< t\le \psi(|x|/2)$ and $y \in B(x,t^{1/\alpha}/(2\Lambda^{1/\alpha}))$,
	\begin{align*}
		p^{\kappa,B(x,|x|/2)}(t,x,y) \ge Ct^{-d/\alpha}.
	\end{align*}
\end{lem}
\pf  By \eqref{e:kappa-cond-1} and  \eqref{e:kappa-cond-2},
\begin{align*}
	\sup_{z \in B(x,|x|/2)} t\kappa(z) \le 
	\Lambda \psi(|x|/2) /\psi(|x|/2)=\Lambda.
\end{align*}
Hence, from  \eqref{e:Feynman-Kac}, we obtain
\begin{align}\label{e:exit-time-lower-0}
		p^{\kappa,B(x,|x|/2)}(t,x,y) \ge e^{-\Lambda}q^{B(x,|x|/2)}(t,x,y).
\end{align}
By \eqref{e:psi-scale}, since $\psi(1)=1$, $\beta_1>\alpha$ and  $|x|\le 1$, it holds that
\begin{align*}
	t\le \Lambda \psi(1) (|x|/2)^{\beta_1}  \le \Lambda (|x|/2)^\alpha.
\end{align*}
Applying Corollary \ref{c:NDL-alpha-stable}  (with $k=\Lambda^{1/\alpha}$), 
we get that $q^{B(x,|x|/2)}(t,x,y) \ge c_1 t^{-d/\alpha}$ for some constant $c_1=c_1(d, \alpha, \Lambda)>0$.
Combining this with \eqref{e:exit-time-lower-0}, we get the desired result. 
\qed 

\begin{lem}\label{l:exit-time-lower}
There exists $C=C(d,\alpha,\Lambda)\in (0,1)$ such that for all $x \in \R^d_0$ with $|x|\le 1$,
	\begin{align*}
		\P_x \big(\tau^\kappa_{B(x,|x|/2)} \ge \psi(|x|/2) \big) \ge C.
	\end{align*}
\end{lem}
\pf  By Lemma \ref{l:exit-time-lower-0}, we obtain
\begin{align*}
	&	\P_x \big(\tau^\kappa_{B(x,|x|/2)} \ge \psi(|x|/2) \big)\\
	&=\int_{B(x,|x|/2)}  p^{\kappa,B(x,|x|/2)}( \psi(|x|/2), x,y )\, dy\ge c_1\psi(|x|/2)^{-d/\alpha}
	\int_{B(x,\, \psi(|x|/2)^{1/\alpha}/(2\Lambda^{1/\alpha}))}  dy =c_2.
\end{align*}
\qed

Using Markov's inequality, from Lemma \ref{l:exit-time-lower}, we deduce the following corollary.
\begin{cor}\label{c:exit-time-lower}
There exists $C=C(d,\alpha,\Lambda)\in (0,1)$ such that for all $x \in \R^d_0$ with $|x|\le 1$,
	\begin{align*}
		\E_x \big[\tau^\kappa_{B(x,|x|/2)}\big] \ge C\psi(|x|/2).
	\end{align*}
\end{cor}

	\section{Survival probability estimates}\label{s:3}

In this section, we continue to assume that
$\psi \in \sC_\alpha(\beta_1,\beta_2,\Lambda)$ and  $\kappa \in \sK^0_\alpha(\psi,\Lambda)$, where $\beta_2\ge \beta_1>\alpha$ and $\Lambda\ge 1$. 
Define a function $H:(0,\infty)\to (0,\infty)$ by
\begin{align}\label{e:def-H}
	H(r)=\frac{2}{r^2} \int_0^r \int_0^s  \psi(u) du ds.
\end{align}
Note that 	$H$ is twice differentiable on $(0,\infty)$.
\begin{lem}\label{l:lemma-H}
	There exists  $C_0=C_0(\beta_2,\Lambda)>1$ such that for all $r>0$,
	\begin{align}
		C_0^{-1}\psi(r)&\le	H(r) \le \psi(r),\label{e:H-psi-compare}\\[3pt]
		|rH'(r)| &\le 2C_0H(r),\nn\\[3pt]
		|r^2H''(r)| &\le 6C_0H(r).\nn
	\end{align}
\end{lem}
\pf Let $r>0$. Since  $\psi$ is increasing, we have
\begin{align*}
	H(r)\le \frac{2\psi(r)}{r^2} \int_0^r \int_0^s   du ds =\psi(r).
\end{align*}
On the other hand, using \eqref{e:psi-scale}, we  get that
\begin{align*}
	H(r) \ge \frac{2}{r^2} \int_{r/2}^r \int_{s/2}^s  \psi(u) du ds \ge \frac{2\psi(r/4)}{r^2} \int_{r/2}^r \int_{s/2}^s   du ds \ge \frac{3\psi(r/4)}{8} \ge \frac{3\psi(r)}{2^{2\beta_2+3} \Lambda }.
\end{align*}
Hence, \eqref{e:H-psi-compare} holds with 
$C_0=2^{2\beta_2+3}\Lambda/3$.
Observe that, since $\psi$ is increasing, 
\begin{align}\label{e:psi-H}
	\int_0^r \int_0^s  \psi(u) du ds   \le r \int_0^r   \psi(u) du \le r^2\psi(r).
\end{align}
Using this in the first and second inequalities below and \eqref{e:H-psi-compare} in the third, we obtain 
\begin{align*}
	&	|rH'(r)| =\left|\frac{2}{r} \int_0^r \psi(u)du - \frac{4}{r^2} \int_0^r \int_0^s  \psi(u) du ds    \right|\le \frac{2}{r} \int_0^r \psi(u)du \le 2\psi(r)\le 2C_0H(r).
\end{align*}
Similarly, using \eqref{e:psi-H} and \eqref{e:H-psi-compare}, we also get that 
\begin{align*}
	r^2H''(r) &= \frac{12}{r^2} \int_0^r \int_0^s  \psi(u) du ds    -\frac{8}{r}  \int_0^r  \psi(u) du  + 2\psi(r) \\
	&\le\frac{4}{r} \int_0^r  \psi(u) du      + 2\psi(r) \le 6 \psi(r) \le 6C_0 H(r)
\end{align*}
and $	r^2H''(r)  \ge - 8 \psi(r) +2\psi(r) \ge -6C_0H(r).$  The proof is complete. \qed

 Let $C_0>1$
	 be the constant in Lemma \ref{l:lemma-H}. 
	The fractional Laplacian can be written as 
	 \begin{align*}-(-\Delta)^{\alpha/2}f(x) = 	\sA(d, -\alpha)	\,p.v. \int_{\R^d} \frac{f(y)-f(x)}{|y-x|^{d+\alpha}} dy,	\end{align*} 	whenever the above principal value integral makes sense.  
	Let $A_{d-1} = 2\pi^{d/2}/\Gamma(d/2)$	denote the hypervolume of the unit sphere in $\R^d$.	  In the remainder of this paper, we let
	\begin{align}\label{e:def-eps-0}
	 \eps_0&:= 
	 \frac18 \wedge \bigg[  \frac{2-\alpha}{ 2^{\alpha-\beta_1+10}(1+2d)3^{\beta_1}
	 \sA(d, -\alpha)C_0^2\Lambda^2
	 A_{d-1}}  \wedge \frac{\alpha}{2^{\alpha+2\beta_1+4} \sA(d, -\alpha)C_0 \Lambda^2 
	 A_{d-1}}  \nn\\
	 &\qquad \qquad \wedge   \frac{\beta_1-\alpha }{2^{d-\alpha+2\beta_1+4} 
	 \sA(d, -\alpha)C_0\Lambda^2 
	 A_{d-1}} \wedge  \frac{\alpha }{2^{d-\alpha+2\beta_1+4}
	 \sA(d, -\alpha)C_0 \Lambda^2
	 A_{d-1}} \bigg]^{1/(\beta_1-\alpha)}
	\end{align}
	and
	\begin{align*}
		\delta_0:= \frac{\alpha \eps_0^\alpha}{2^{d-\alpha+4}
		\sA(d, -\alpha)C_0\Lambda A_{d-1}}.
	\end{align*}
	Note that
\begin{align}\label{e:eps0-delta0}
\frac{\Lambda(4\eps_0)^{\beta_1}}{\delta_0} =  \frac{2^{d-\alpha+2\beta_1+4}\sA(d,-\alpha)C_0 \Lambda^2A_{d-1}\eps_0^{\beta_1- \alpha}}{\alpha } \le 1.
\end{align}
Further, for all $R\in (0,1]$ and $x \in B(0,4\eps_0R)$, by \eqref{e:H-psi-compare},  \eqref{e:psi-scale} (with $\psi(1)=1$) and  \eqref{e:eps0-delta0}, we have
\begin{align}\label{e:generator-1}
	H(|x|) \le \psi(|x|) \le 
	\Lambda |x|^{\beta_1} 
	\le \delta_0 R^{\beta_1}
	\le \delta_0R^\alpha.
\end{align}
For each  $R\in (0,1]$,  let
 $\phi_{R}:\R^d \to [0,\infty)$  be an element of $C^2(\R^d_0)$ satisfying the following properties:

\smallskip

(1) $\phi_{R}(x)=H(|x|)$ for $x \in B(0,4\eps_0R) \setminus \{0\}$.

(2) $H(\eps_0R) \le  \phi_{R}(x) \le	\delta_0 R^{\alpha} $   for $x\in B(0,R)\setminus B(0,4\eps_0R)$.

(3) $\phi_{R}(x)=\delta_0 R^\alpha$  for $x \in B(0,2^{1/d}R)\setminus B(0,R)$.

(4) $0\le\phi_{R}(x)\le\delta_0 R^\alpha$  for $x \in B(0,4^{1/d}R) \setminus B(0,2^{1/d}R)$.

(5) $\phi_R(x)=0$ for $x \in B(0,4^{1/d}R)^c$.

\begin{lem}\label{l:genenrator}
 For all $R \in (0,1]$ and  $x \in \R^d_0$ with $|x|<\eps_0R$,  it holds that
	\begin{align*}
	 	-(-\Delta)^{\alpha/2} \phi_{R}(x) \le 
		\frac{1}{4C_0\Lambda}.
	\end{align*}
\end{lem}
\pf Let $R\in (0,1]$ and $x=(x_1,\cdots,x_d) \in B(0,\eps_0R)\setminus \{0\}$.  Observe that
\begin{align*}
	&	-(-\Delta)^{\alpha/2} \phi_{R}(x)\\
	&=	\sA(d, -\alpha)
	\bigg[\, p.v. \int_{B(x,|x|/2)}  \frac{H(|y|) - H(|x|)}{|y-x|^{d+\alpha}} dy + \int_{B(x,3|x|)\setminus B(x,|x|/2)}  \frac{H(|y|) - H(|x|)}{|y-x|^{d+\alpha}} dy  \\
		&\qquad\qquad \qquad + \int_{B(0,4\eps_0R)\setminus B(x,3|x|)}  \frac{\phi_R(y)- H(|x|)}{|y-x|^{d+\alpha}} dy  +  \int_{B(0,4^{1/d}R)\setminus B(0,4\eps_0R)}  \frac{\phi_R(y) - H(|x|)}{|y-x|^{d+\alpha}} dy \\
		&\qquad\qquad \qquad  +  \int_{B(0,4^{1/d}R)^c}  \frac{ - H(|x|)}{|y-x|^{d+\alpha}} dy \bigg] \\
		&=:		\sA(d, -\alpha)
		(I_1+I_2+I_3+I_4+I_5).
\end{align*}
By  symmetry, we have
\begin{align*}
	I_{1} & =\lim_{\eps \to 0} \bigg(  \int_{B(x, |x|/2)\setminus B(x,\eps)}  \frac{H(|y|) - H(|x|) - H' (|x|)|x|^{-1}  x \cdot (y-x) }{|y-x|^{d+\alpha}} dy \bigg).
\end{align*}
 For any $y=(y_1,\cdots,y_d) \in B(x,|x|/2)$,  using  Taylor's theorem in the first inequality below, Lemma \ref{l:lemma-H} in the second,  the inequality   $\sum_{1\le i<j\le d} |x_i-y_i||x_j-y_j| \le   d\sum_{1\le i\le d} |x_i-y_i|^2=d|x-y|^2 $ in the third, and  \eqref{e:psi-scale} and $\psi(1)=1$ in the fourth, we obtain
\begin{align*}
	&|H(|y|) - H(|x|) - H' (|x|)|x|^{-1}  x \cdot (y-x)| \\
	& \le \sup_{z=(z_1,\cdots, z_d) \in \{x + s(y-x):s\in [0,1]\} }\bigg|\, \frac12\sum_{1\le i\le d} \bigg( H''(|z|)\frac{z_i^2}{|z|^2} + H'(|z|) \frac{|z|^2-z_i^2}{|z|^3} \bigg)  |x_i-y_i|^2	 \\
	&\qquad \qquad \qquad \qquad \qquad \qquad \qquad  + \sum_{1\le i<j\le d} \bigg( H''(|z|)\frac{z_iz_j}{|z|^2} - H'(|z|) \frac{z_iz_j}{|z|^3} \bigg)  |x_i-y_i||x_j-y_j|\, \bigg| \\
		& \le \sup_{ z\in \R^d:|x|/2\le |z| \le 3|x|/2 }\bigg[\, \sum_{1\le i\le d}  \frac{4C_0\psi(|z|)}{|z|^2}   |x_i-y_i|^2	 + \sum_{1\le i<j\le d}  \frac{8C_0\psi(|z|)}{|z|^2}  |x_i-y_i||x_j-y_j|\, \bigg] \\
	& \le \sup_{z\in \R^d:|x|/2\le |z| \le 3|x|/2  }\bigg[\,  \frac{(4+8d)C_0\psi(|z|)}{|z|^2}  \bigg]  |x-y|^2	 \\
		& \le \sup_{z\in \R^d:|x|/2\le |z| \le 3|x|/2  }\bigg[\,  \frac{(4+8d)C_0 
		\Lambda |z|^{\beta_1}}{|z|^{2}}  \bigg]  |x-y|^2	 \\
	&\le 16(1+2d)(3/2)^{\beta_1}C_0  
	\Lambda |x|^{\beta_1-2} |x-y|^2  .
\end{align*}
Thus, since $|x|\le \eps_0$ and $\beta_1>\alpha$, it holds that 
\begin{align*}
	|I_{1}|&\le 16(1+2d)(3/2)^{\beta_1}C_0 
	\Lambda |x|^{\beta_1-2}\int_{B(x,|x|/2)} \frac{dy}{|y-x|^{d+\alpha-2}} \\
	&=\frac{16(1+2d)(3/2)^{\beta_1}C_0   
	\Lambda A_{d-1}  |x|^{\beta_1-2}(|x|/2)^{2-\alpha}}{2-\alpha}\\
	&  \le  \frac{16(1+2d) 3^{\beta_1}C_0  
	\Lambda  A_{d-1} 
	\eps_0^{\beta_1-\alpha}}{2^{2-\alpha+\beta_1}(2-\alpha)} \le \frac{1}{16\sA(d,-\alpha)C_0\Lambda}.
\end{align*}
For $I_2$, using \eqref{e:generator-1} and the facts that $|x|\le \eps_0$ and $\beta_1>\alpha$, we obtain
\begin{align*}
|I_2| \le  \int_{B(x,3|x|)\setminus B(x,|x|/2)}  
\frac{\Lambda(4|x|)^{\beta_1}}{|y-x|^{d+\alpha}}dy 
\le \frac{4^{\beta_1}\Lambda A_{d-1} |x|^{\beta_1}}{\alpha (|x|/2)^\alpha} 
\le \frac{2^{\alpha+2\beta_1}
\Lambda A_{d-1}   
\eps_0^{\beta_1-\alpha}}{\alpha }\le \frac{1}{16\sA(d,-\alpha)C_0\Lambda}.
\end{align*}
 Note that for  $y \in B(x,3|x|)^c$,  we have $|y| \ge 2|x|$. Hence,
 \begin{align}\label{e:generator-2}
 	|y-x| \ge |y|/2 \quad \text{for } \;\, y \in B(x,3|x|)^c.
 \end{align} 
Further, for any $y \in B(0,4\eps_0R)\setminus B(x,3|x|)$, since $|y| \ge 2|x|$, by \eqref{e:generator-1}, it holds that 
\begin{align}\label{e:generator-3}
	 |\phi_R(y)-H(|x|)| \le H(|y|)\vee H(|x|) \le 
	 \Lambda |y|^{\beta_1}.
\end{align}
 Using \eqref{e:generator-2} and \eqref{e:generator-3}, we obtain
\begin{align*}
	|I_{3}| &\le 2^{d+\alpha} 
	\Lambda 
	\int_{B(0,4\eps_0R)\setminus B(x,3|x|)}  \frac{dy}{|y|^{d+\alpha-\beta_1}}dy \le \frac{2^{d+\alpha} 
	\Lambda  
	A_{d-1}(4\eps_0R)^{\beta_1-\alpha}}{\beta_1-\alpha}   \le 
	\frac{1}{16\sA(d,-\alpha)C_0\Lambda}.
\end{align*} 
 For $I_4$, we note that by \eqref{e:generator-1} and \eqref{e:eps0-delta0}, since $\beta_1>\alpha$ and $R\le 1$,
 \begin{align}\label{e:generator-5}
 	H(|x|)\le 
	\Lambda |x|^{\beta_1}  \le \Lambda \eps_0^{\beta_1} R^\alpha 
	\le \delta_0 R^\alpha.
 \end{align} 
 Therefore, using \eqref{e:generator-2}, we get that  
\begin{align*}
|I_4|  \le 2^{d+\alpha} \delta_0 R^\alpha \int_{B(0,4^{1/d}R) \setminus B(0,4\eps_0R)} \frac{dy}{|y|^{d+\alpha}} \le \frac{ 2^{d+\alpha} \delta_0 A_{d-1}}{\alpha(4\eps_0 )^\alpha}  =
 \frac{1}{16\sA(d,-\alpha)C_0\Lambda}.
\end{align*}
For $I_5$, by \eqref{e:generator-2} and \eqref{e:generator-5}, we have
\begin{align*}
	0 \ge I_5 \ge - 2^{d+\alpha}\delta_0 R^\alpha \int_{B(0,4^{1/d}R)^c} \frac{dy}{|y|^{d+\alpha}} = - \frac{2^{d+\alpha-2\alpha/d} \delta_0 A_{d-1}}{\alpha}.
\end{align*}
Combining the estimates for $I_1,I_2,I_3,I_4$ and $I_5$ above, we deduce that $-(-\Delta)^{\alpha/2}\phi_R(x)$ is well-defined and that
\begin{align*}
	-(-\Delta)^{\alpha/2}\phi_R(x) \le c_{d,\alpha}(I_1+I_2+I_3+I_4) \le 
	\frac{1}{4C_0\Lambda}.
\end{align*} 
The proof is complete.\qed 

\begin{cor}\label{c:generator-sub}
	For all $R \in (0,1]$ and  $x \in \R^d_0$ with $|x|<\eps_0R$,  we have 
	\begin{align*}
		L_{\alpha}^{\kappa/2} \phi_{R}(x) \le 
		-\frac{1}{4C_0\Lambda}.
	\end{align*}
\end{cor}
\pf Using Lemma \ref{l:genenrator}, \eqref{e:kappa-cond-1} and \eqref{e:H-psi-compare}, we get that
\begin{align*}
	L_{\alpha}^{\kappa/2} \phi_{R}(x) = -(-\Delta)^{\alpha/2} \phi_{R}(x) - \frac12\kappa(x) H(|x|)  \le  \frac{1}{4C_0\Lambda} - \frac{\psi(|x|)}{2C_0\Lambda \psi(|x|)} = - \frac{1}{4C_0\Lambda}.
\end{align*}
\qed

Under the assumptions of this section,  conditions  \cite[{\bf (H1)}--{\bf (H4)} and (1.2)--(1.3)]{KSV23} hold. Thus, by \cite[Theorem 4.8]{KSV23}, we have the following Dynkin-type theorem for 
$L^{\kappa/2}_\alpha$. Note that,
although \cite{KSV23} imposes a stronger condition  on $\kappa$ (see (1.8) therein),  \cite[Theorem 4.8]{KSV23} only relies on the local boundedness of $\kappa$, which  is satisfied in our context since $\kappa \in \sK^0_\alpha(\psi,\Lambda)$.

\begin{thm}\label{t:Dynkin}
	Let $U$ be a relatively compact subset of $\R^d_0$. For any non-negative function $u$ defined on $\R^d_0$ satisfying $u \in C^2(\overline U)$ and any $x \in U$,
	\begin{align*}
		\E_x\Big[ u(X^{\kappa/2}_{\tau^{\kappa/2}_{U} })\Big]  = u(x) + \E_x \int_0^{\tau^{\kappa/2}_U}
		L^{\kappa/2}_\alpha u (X^{\kappa/2}_s)ds. 
	\end{align*}
\end{thm}

\begin{prop}\label{p:exit-time-upper}
There exists $C=C(d,\alpha,\beta_1,\beta_2,\Lambda)>0$ such that for all $R\in (0,1]$ and $x\in \R^d_0$ with $|x|<\eps_0R$,
	\begin{align}\label{e:exit-time-upper}
	\E_x \big[\tau^{\kappa}_{B(0,\eps_0R)}  \big]\le 	\E_x \big[\tau^{\kappa/2}_{B(0,\eps_0R)}  \big] \le C\psi(|x|).
	\end{align}
\end{prop}
\pf 
Let $R\in (0,1]$. The first inequality in \eqref{e:exit-time-upper} is obvious. We now prove the  second inequality.  Take $\eps \in (0,\eps_0R)$. For any  $x \in B(0,\eps_0 R) \setminus \{0\}$, using  Theorem \ref{t:Dynkin},  Corollary \ref{c:generator-sub} and \eqref{e:H-psi-compare},  since $\phi_R(x)=H(|x|)$ by definition, we obtain
\begin{equation}\label{e:Dynkin-sub}
	\E_x\Big[ \phi_R\big(X^{\kappa/2}_{\tau^{\kappa/2}_{B(0,\eps_0R)\setminus \overline{B(0,\eps)}} }\big)\Big]  = \phi_R(x) + \E_x \int_0^{\tau^{\kappa/2}_{B(0,\eps_0R)\setminus \overline{B(0,\eps)}}} 
	L^{\kappa/2}_\alpha u (X^{\kappa/2}_s)ds   \le H(|x|) \le \psi(|x|).
\end{equation}
On the other hand, by \eqref{e:Levysystem}, we have
\begin{align}\label{e:Dynkin-sub-2}
	\E_x\Big[\phi_R\big(X^{\kappa/2}_{\tau^{\kappa/2}_{B(0,\eps_0R)\setminus \overline{B(0,\eps)}} }\big) \Big]&\ge \E_x \left[ \int_0^{\tau^{\kappa/2}_{B(0,\eps_0R)\setminus \overline{B(0,\eps)}}  } \int_{B(0,2^{1/d}R)\setminus B(0,R)} \frac{\sA(d, -\alpha)
		\phi_R(y)}{|X^{\kappa/2}_s-y|^{d+\alpha}} \,dy \,ds \right]\nn\\
	&\ge \frac{\sA(d, -\alpha)
		\delta_0 R^\alpha}{ (2^{1/\alpha}+\eps_0)^{d+\alpha}R^{d+\alpha}} \E_x \big[\tau^{\kappa/2}_{B(0,\eps_0R)\setminus \overline{B(0,\eps)}}  \big]\int_{B(0,2^{1/d}R)\setminus B(0,R)}dy \nn\\
	&= \frac{ \sA(d, -\alpha)
		\delta_0 A_{d-1} }{d  (2^{1/\alpha}+\eps_0)^{d+\alpha}} \E_x \big[\tau^{\kappa/2}_{B(0,\eps_0R)\setminus \overline{B(0,\eps)}} \big].
\end{align}
Since $\eps\in (0,\eps_0R)$ is arbitrary, combining \eqref{e:Dynkin-sub} with \eqref{e:Dynkin-sub-2} and  applying the monotone convergence theorem, we conclude that
\begin{align*}
	\E_x \big[\tau^{\kappa/2}_{B(0,\eps_0R)}  \big]   = \lim_{\eps \to 0}  \E_x \big[\tau^{\kappa/2}_{B(0,\eps_0R)\setminus \overline{B(0,\eps)}} \big] \le  \frac{d  (2^{1/\alpha}+\eps_0)^{d+\alpha}} { \sA(d, -\alpha)
		\delta_0 A_{d-1} } \psi(|x|).
\end{align*} 
\qed

Denote by $\zeta^\kappa$ the lifetime of $X^\kappa$.

\begin{lem}\label{l:survival}
	There exists $C=C(d,\alpha,\beta_1,\beta_2,\Lambda)>0$  such that for all $R\in (0,1]$ and $x\in \R^d_0$ with $|x|<\eps_0R$,
	\begin{align}\label{e:survival}
		\P_x \big(\tau^\kappa_{B(0,\eps_0R)}<\zeta^\kappa \big) \le 	\P_x \big(\tau^{\kappa/2}_{B(0,\eps_0R)}<\zeta^{\kappa/2} \big) =  \P_x \Big(X^{\kappa/2}_{\tau^{\kappa/2}_{B(0,\eps_0R)}} \in \R^d_0 \Big) \le \frac{C\psi(|x|)}{\psi(R)}.
	\end{align}
\end{lem}
\pf The first inequality in \eqref{e:survival} is evident. We now present the proof of the second inequality.  For any $z \in B(0,R)\setminus B(0,\eps_0R)$, by  \eqref{e:H-psi-compare}, 
$
\phi_R(z) \ge  C_0^{-1}\psi(\eps_0R).
$
Hence, by  Markov's inequality, \eqref{e:psi-scale} and \eqref{e:Dynkin-sub}, it holds  that
\begin{align*}
 \P_x \Big(X^{\kappa/2}_{\tau^{\kappa/2}_{B(0,\eps_0R)}} \in B(0,R) \Big) \le  \frac{C_0}{\psi(\eps_0R)}	\E_x\Big[ \phi_R\big(X^{\kappa/2}_{\tau^{\kappa/2}_{B(0,\eps_0R)} }\big)\Big] \le \frac{C_0\Lambda\psi(|x|)}{\eps_0^{\beta_2}\psi(R)}.
\end{align*}
On the other hand, by \eqref{e:Levysystem}, we have
\begin{align}\label{e:survival-1}
	 \P_x \Big(X^{\kappa/2}_{\tau^{\kappa/2}_{B(0,\eps_0R)}} \in  B(0,R)^c \Big)   = \E_x\left[ \int_0^{\tau^{\kappa/2}_{B(0,\eps_0R)} } \int_{ B(0,R)^c} \frac{
	 \sA(d, -\alpha)
	 }{|X^{\kappa/2}_s-w|^{d+\alpha}}\,dw ds \right]. 
\end{align}
For any $z \in B(0,\eps_0R)$, by \eqref{e:psi-scale}, since $\psi(1)=1$, $R\le 1$ and $\beta_1>\alpha$, it holds that
\begin{align*}
	 \int_{ B(0,R)^c}  \frac{	 \sA(d, -\alpha)
	 }{|z-w|^{d+\alpha}}dw \le  \int_{B(z,R/2)^c} \frac{\sA(d, -\alpha)
	 }{|z-w|^{d+\alpha}}dw = \frac{2^{\alpha}
	 \sA(d, -\alpha)
	 A_{d-1}}{R^{\alpha}\psi(1)} \le \frac{c_1}{R^{\alpha-\beta_1}\psi(R)} \le\frac{c_1}{\psi(R)}.
\end{align*}
Combining this with \eqref{e:survival-1} and using Proposition \ref{p:exit-time-upper}, we obtain
\begin{align*}
 \P_x \Big(X^{\kappa/2}_{\tau^{\kappa/2}_{B(0,\eps_0R)}} \in B(0,R)^c \Big)  \le \frac{c_1}{\psi(R)}	\E_x \big[\tau^{\kappa/2}_{B(0,\eps_0R)}  \big] \le \frac{c_2\psi(|x|)}{\psi(R)}.
\end{align*}
The proof is complete. \qed

\begin{lem}\label{l:survival-t}
	There exists $C=C(d,\alpha,\beta_1,\beta_2,\Lambda)>0$  such that for all $t\in (0,1]$ and $x \in \R^d_0$,
	\begin{align}\label{e:survival-t}
		\P_x ( \zeta^{\kappa/2}>t) \le C\bigg(1 \wedge \frac{\psi(|x|)}{t}\bigg).
	\end{align}
\end{lem}
\pf If $|x| \geq \eps_0\psi^{-1}(t)$, then $\psi(|x|)/t \ge c_1$ by \eqref{e:psi-scale}, and hence \eqref{e:survival-t} holds.
If $|x|<\eps_0\psi^{-1}(t)$, then using Lemma \ref{l:survival} and Markov's inequality in the second line below, and Proposition \ref{p:exit-time-upper} in the third, we obtain
\begin{align*}
	\P_x ( \zeta^{\kappa/2}>t)&\le 	\P_x \big(\tau^{\kappa/2}_{B(0,\eps_0\psi^{-1}(t))} <\zeta^{\kappa/2} \big) + \P_x \big(\tau^{\kappa/2}_{B(0,\eps_0\psi^{-1}(t))} >t \big)\\
	&\le c_2 t^{-1}\psi(|x|) + t^{-1}\E_x[\tau^{\kappa/2}_{B(0,\eps_0\psi^{-1}(t))}]\le c_3 t^{-1}\psi(|x|).
\end{align*}
 \qed

\begin{lem}\label{l:EP}
	There exists  $C=C(d,\alpha,\beta_1,\beta_2,\Lambda)>0$   such that for all $r \in (0,\eps_0]$,  $ t\in (0,1]$ and $x \in \R^d_0$ with $|x|<r$,
	\begin{align*}
		\P_x\big(\tau^{\kappa/2}_{B(0,3r)} <t \wedge \zeta^{\kappa/2} \big) \le \frac{C\psi(|x|)}{r^\alpha} \bigg( 1\vee \frac{t}{\psi(r)}\bigg).
	\end{align*}
\end{lem}
\pf Using the strong Markov property, we see that
\begin{align*}
	&\P_x\big(\tau^{\kappa/2}_{B(0,3r)} <t \wedge \zeta^{\kappa/2} \big) \le 	\P_x\big( X^{\kappa/2}_{\tau^{\kappa/2}_{B(0,r)}}  \in B(0,2r)^c \big)  + \P_x\big( X^{\kappa/2}_{\tau^{\kappa/2}_{B(0,r)}} \in  B(0,2r), \, \tau^{\kappa/2}_{B(0,3r)} <t \wedge \zeta^{\kappa/2} \big)\\
	&\le 	\P_x\big( X^{\kappa/2}_{\tau^{\kappa/2}_{B(0,r)}} \in B(0,2r)^c \big)  + \P_x\big( X^{\kappa/2}_{\tau^{\kappa/2}_{B(0,r)}} \in  B(0,2r) \big) \sup_{z \in B(0,2r)} \P_z \big( \tau^{\kappa/2}_{B(z,r)}<t \wedge \zeta^{\kappa/2} \big)\\
	&=:I_1+I_2.
\end{align*}
For any $z \in B(0,r)$, it holds that
\begin{align}\label{e:EP-1}
	\int_{ B(0,2r)^c} \frac{\sA(d, -\alpha)
	}{|z-y|^{d+\alpha}}dy \le 	\int_{B(z,r)^c} \frac{
	\sA(d, -\alpha)
	}{|z-y|^{d+\alpha}}dy  = \frac{
\sA(d, -\alpha)
	A_{d-1}}{dr^\alpha}.
\end{align}
Using  \eqref{e:Levysystem} in the equality below,  \eqref{e:EP-1} in the first inequality and  Proposition \ref{p:exit-time-upper} in the second inequality, we obtain
\begin{align*}
	I_1&= \E_x\bigg[\int_{0}^{\tau^{\kappa/2}_{B(0,r)}} \int_{ B(0,2r)^c} \frac{
	\sA(d, -\alpha)
	}{|X^{\kappa/2}_s- y|^{d+\alpha}}dy ds \bigg] \le  \frac{c_1}{r^\alpha} \E_x \big[ \tau^{\kappa/2}_{B(0,r)}\big] \le \frac{c_2 \psi(|x|)}{r^\alpha} .
\end{align*}
Further, by Lemma \ref{l:survival} and Proposition \ref{p:EP-alpha-stable}, we have
\begin{align*}
	I_2\le  \P_x\big(\tau^{\kappa/2}_{B(0,r)} <\zeta^{\kappa/2} \big) \sup_{z \in B(0,2r)} \P_z \big( \tau^{Y}_{B(z,r)}<t  \big) \le \frac{c_3\psi(|x|)t}{\psi(r)r^\alpha}.
\end{align*}
The proof is complete.\qed

	\section{Small time estimates}\label{s:4}
	In this section, we continue to  assume that $\kappa \in \sK^0_\alpha(\psi,\Lambda)$ and  that  $\eps_0=\eps_0(d,\alpha,\beta_1,\beta_2,\Lambda)\in (0,1/8]$ is the constant in \eqref{e:def-eps-0}.	The goal
	of this section is to establish the following  two-sided small time heat kernel estimates.
	
	\begin{thm}\label{t:smalltime}
	Suppose that $\kappa \in \sK^0_\alpha(\psi,\Lambda)$.	Let $T\ge1$. There exist constants $\lambda_1=\lambda_1(\beta_2,\Lambda), \lambda_2=\lambda_2(\beta_2,\Lambda)>0$ and $C= C(d,\alpha,\beta_1,\beta_2,\Lambda,T)\ge1$ such that  for all $t \in (0,T]$ and $x,y \in \R^d_0$,
	\begin{align}\label{e:smalltime-upper}
 p^{\kappa}(t,x,y) &\le C\bigg(1 \wedge \frac{\psi(|x|)}{t}\bigg) \bigg(1 \wedge \frac{\psi(|y|)}{t}\bigg)\nn\\
	&\quad \times  \bigg[\, e^{-\lambda_1t/\psi(|x|\vee |y|)}t^{-d/\alpha}\bigg( 1 \wedge \frac{t^{1/\alpha}}{|x-y|}\bigg)^{d+\alpha} +  \frac{t^2}{\psi^{-1}(t)^{d+2\alpha}}\bigg(1 \wedge \frac{\psi^{-1}(t)}{|x-y|}\bigg)^{d+2\alpha}\, \bigg]
\end{align}
and
	\begin{align}\label{e:smalltime-lower}
	p^\kappa(t,x,y)&\ge C^{-1}\bigg(1 \wedge \frac{\psi(|x|)}{t}\bigg) \bigg(1 \wedge \frac{\psi(|y|)}{t}\bigg) \nn\\
	&\quad \times \bigg[\, e^{-\lambda_2t/\psi(|x|\vee |y|)}t^{-d/\alpha}\bigg( 1 \wedge \frac{t^{1/\alpha}}{|x-y|}\bigg)^{d+\alpha} +  \frac{t^2}{\psi^{-1}(t)^{d+2\alpha}}\bigg(1 \wedge \frac{\psi^{-1}(t)}{|x-y|}\bigg)^{d+2\alpha}\, \bigg].
\end{align}
	\end{thm} 
	
	The proofs for \eqref{e:smalltime-upper} and \eqref{e:smalltime-lower}  will be  	given at the end of Subsections  \ref{ss:smalltime-upper} and \ref{ss:smalltime-lower}, respectively. We first
	record a consequence of \eqref{e:smalltime-upper}, see Corollary \ref{c:UHK-general} below.

	\begin{lem}\label{l:dominant}
		For any $T\ge 1$,	there exists $C=C(d,\alpha,\beta_1,\Lambda,T)>0$ such that for all $t \in (0,T]$ and $r>0$,
		\begin{align*}
			\frac{t^2}{\psi^{-1}(t)^{d+2\alpha}}\bigg(1 \wedge \frac{\psi^{-1}(t)}{r}\bigg)^{d+2\alpha}\le Ct^{-d/\alpha}\bigg( 1 \wedge \frac{t^{1/\alpha}}{r}\bigg)^{d+\alpha} .
		\end{align*}
	\end{lem}
	\pf Let $t\in (0,T]$ and $r>0$.  If $r> (t/T)^{1/\beta_1}$, then   $r^\alpha > (t/T)^{\alpha/\beta_1} \ge t/T$. Hence, 
	\begin{align*}
		\frac{t^2}{\psi^{-1}(t)^{d+2\alpha}}\bigg(1 \wedge \frac{\psi^{-1}(t)}{r}\bigg)^{d+2\alpha} \le 	\frac{t^2}{r^{d+2\alpha}} \le \frac{Tt}{r^{d+\alpha}}\le c_1t^{-d/\alpha}\bigg( 1 \wedge \frac{t^{1/\alpha}}{r}\bigg)^{d+\alpha}.
	\end{align*}
	If $(t/T)^{1/\alpha}\le r \le (t/T)^{1/\beta_1}$, then by \eqref{e:psi-inverse}, we get
	\begin{align*}
		&\frac{t^2}{\psi^{-1}(t)^{d+2\alpha}}\bigg(1 \wedge \frac{\psi^{-1}(t)}{r}\bigg)^{d+2\alpha} \le \frac{t^2}{\psi^{-1}(t/T)^{d+2\alpha}}\le  \frac{t^2}{(t/(\Lambda T))^{(d+2\alpha)/\beta_1}} \\
		&\le \frac{\Lambda^{(d+2\alpha)/\beta_1}t^2}{r^{d+2\alpha}} \le \frac{\Lambda^{(d+2\alpha)/\beta_1}Tt}{r^{d+\alpha}}\le c_2t^{-d/\alpha}\bigg( 1 \wedge \frac{t^{1/\alpha}}{r}\bigg)^{d+\alpha}.
	\end{align*}
	If $r<(t/T)^{1/\alpha}$, then by \eqref{e:psi-inverse}, since $t/(\Lambda T) \le 1$ and $\beta_1>\alpha$, we see that
	\begin{align*}
		&	\frac{t^2}{\psi^{-1}(t)^{d+2\alpha}}\bigg(1 \wedge \frac{\psi^{-1}(t)}{r}\bigg)^{d+2\alpha} \le \frac{t^2}{(t/(\Lambda T))^{(d+2\alpha)/\beta_1}}\\
		& \le \frac{t^2}{(t/(\Lambda T))^{(d+2\alpha)/\alpha}}\le c_3t^{-d/\alpha}\bigg( 1 \wedge \frac{t^{1/\alpha}}{r}\bigg)^{d+\alpha}.
	\end{align*}
	The proof is complete. \qed

	\begin{cor}\label{c:UHK-general}
		There exists $C=C(d,\alpha,\beta_1,\beta_2,\Lambda)>0$   such that for all $t >0$ and $x,y \in \R^d_0$,
		\begin{equation}\label{e:UHK-general}
			p^{\kappa}(t,x,y) \le C\bigg(1 \wedge \frac{\psi(|x|)}{t \wedge 1}\bigg) \bigg(1 \wedge \frac{\psi(|y|)}{t\wedge 1}\bigg)\,\wt q(t,x,y).
		\end{equation}
	\end{cor} 
	\pf Let $t>0$ and $x,y \in \R^d_0$. If $t\le 2$, then  using  	\eqref{e:smalltime-upper}
	and Lemma \ref{l:dominant}, we get \eqref{e:UHK-general}. Suppose that $t>2$. Using the semigroup property of $p^\kappa$ in the first line below, the inequality $p^\kappa(s,v,w) \le q(s,v,w)$ for all $s>0$ and $v,w\in \R^d_0$,  \eqref{e:HKE-alpha-stable} and  \eqref{e:UHK-general} with $t=1$  in the second,  and \eqref{e:wtq-semigroup} in the third, we obtain
	\begin{align*}
		p^\kappa(t,x,y)&= \int_{\R^d_0}\int_{\R^d_0} p^\kappa(1,x,v) p^\kappa(t-2,v,w) p^\kappa(1, w,y) dvdw\\
		&\le c_1(1 \wedge \psi(|x|))(1 \wedge \psi(|y|))  \int_{\R^d_0}\int_{\R^d_0}  \wt q(1,x,v) \wt q(t-2,v,w) \wt q(1, w,y)  dvdw\\
		&\le c_2(1\wedge \psi(|x|)) ( 1\wedge \psi(|y|)) \wt q(t,x,y).
	\end{align*}
	The proof is complete. \qed

	\subsection{Small time upper heat kernel estimates}\label{ss:smalltime-upper}
		\begin{lem}\label{l:upper-HKE-0}
		For any $T\ge 1$, there exists  $C=C(d,\alpha,\beta_1,\beta_2,\Lambda,T)>0$  such that for all $t \in (0,T]$ and $x,y \in \R^d_0$,
		\begin{align}\label{e:upper-HKE-0}
		p^{\kappa}(t,x,y)\le 	p^{\kappa/2}(t,x,y) \le \frac{C\psi(|x|)}{t^{d/\alpha+1}}.
		\end{align}
	\end{lem} 
	\pf Clearly, $p^\kappa(t,x,y)\le p^{\kappa/2}(t,x,y)$.  We now prove the second inequality in \eqref{e:upper-HKE-0}. If $|x| \ge \eps_0 \psi^{-1}(t/T)$, then by \eqref{e:psi-scale}, $\psi(|x|) \ge c_1t$. Hence, since $p^{\kappa/2}(t,x,y)\le q(t,x,y)$,  \eqref{e:upper-HKE-0} follows from \eqref{e:HKE-alpha-stable}. Suppose that  $|x| < \eps_0 \psi^{-1}(t/T)$. 	
	Using the semigroup property in the equality below, \eqref{e:HKE-alpha-stable} and Lemma \ref{l:survival-t} in the second inequality, we get
	\begin{align*}
			p^{\kappa/2}(t,x,y) &= \int_{\R^d_0} 	p^{\kappa/2}(t/(2T),x,z)\,	p^{\kappa/2}((2T-1)t/(2T),z,y) \,dz\\
			&\le \sup_{z \in \R^d_0} 	q((2T-1)t/(2T),z,y)\, \P_x\big(X^{\kappa/2}_{t/(2T)} \in \R^d_0\big)\\
			&\le 2Tc_2  (t/2)^{-d/\alpha} \psi(|x|)/t,
	\end{align*}
	proving that the second inequality in \eqref{e:upper-HKE-0} holds.
	The proof is complete.	\qed
	
	For any Borel set $E \subset \R^d$ with positive Lebesgue measure and 	$f\in L^1(E)$, we use the usual notation $\fint_E f := \int_E f/ \int_E 1$.
	
	\begin{lem}\label{l:upper-HKE-1}
	For any $T\ge 1$, there exists  $C=C(d,\alpha,\beta_1,\beta_2,\Lambda,T)>0$ such that for all $t \in (0,T]$ and $x,y \in \R^d_0$ with $|x|< \eps_0\psi^{-1}(t/T)/8$ and $|y|> 8|x|$,
		\begin{align*}
		 p^{\kappa}(t,x,y)\le p^{\kappa/2}(t,x,y) \le \frac{C\psi(|x|)}{|x-y|^{d+\alpha}}\bigg( 1\vee \frac{t}{\psi(|y|)}\bigg).
		\end{align*}
	\end{lem} 
	\pf Let $t \in (0,T]$ and $x,y \in \R^d_0$ be such that $|x|< \eps_0\psi^{-1}(t/T)/8$ and $|y|> 8|x|$. Set
	$$
	r_0:=\frac{|y| \wedge (\eps_0\psi^{-1}(t/T))}{8} \quad \text{and} \quad U:=B(0,r_0).
	$$
	 By the lower semi-continuity of $p^{\kappa/2}$ and the strong Markov property of $X^{\kappa/2}$, we have
\begin{align*}
&p^{\kappa/2}(t,x,y) \le \liminf_{\delta\to 0} \fint_{B(y,\delta)}	p^{\kappa/2}(t,x,v)dv\\
 &\le\limsup_{\delta\to 0} \E_x\left[ \fint_{B(y,\delta)} p^{\kappa/2}(t-\tau^{\kappa/2}_{U} , X^{\kappa/2}_{\tau_{U}} , v)dv : \tau^{\kappa/2}_{U} <t \wedge \zeta^\kappa,\, X^{\kappa/2}_{\tau^{\kappa/2}_U} \in B(y, |x-y|/2)\right]\nn\\
	& \quad +\limsup_{\delta\to 0} \E_x\left[ \fint_{B(y,\delta)} p^{\kappa/2}(t-\tau^{\kappa/2}_{U} , X^{\kappa/2}_{\tau^{\kappa/2}_U}, v)dv : \tau^{\kappa/2}_{U} <t \wedge \zeta^\kappa,\, X^{\kappa/2}_{\tau^{\kappa/2}_{U}} \in B(y, |x-y|/2)^c\right] \nn\\
	&=:I_1+I_2.
\end{align*}
For all $z \in U$ and $w \in B(y,|x-y|/2)$, since $|x|<r_0 \le |y|/8$ so that $|x|<r_0 \le |x-y|/7$,
\begin{equation}\label{e:upper-HKE-1}
	|z-w| \ge |w| -|z|  \ge |x-y|-|x|-\frac{|x-y|}{2}  - r_0 \ge \frac{3|x-y|}{14}.
	\end{equation} 
Using  \eqref{e:Levysystem} in the equality below, \eqref{e:upper-HKE-1} and the symmetry of $p$ in the first inequality  and Proposition \ref{p:exit-time-upper} in the last, we obtain
\begin{align*}
	I_1&=\limsup_{\delta\to 0} \int_0^t \int_{U} \int_{B(y,|x-y|/2)} \fint_{B(y,\delta)}  p^{\kappa/2,U}(s, x, z) \frac{c_{d,\alpha}}{|z-w|^{d+\alpha}} p^{\kappa/2}(t-s, w, v) \, dv dw dz ds \\
	&\le\limsup_{\delta\to 0} \frac{(14/3)^{d+\alpha}c_{d,\alpha}}{|x-y|^{d+\alpha}} \int_0^t \int_{U}p^{\kappa/2,U}(s, x, z) dz  \fint_{B(y,\delta)} \int_{B(y,|x-y|/2)}  p^{\kappa/2}(t-s, v, w) dw    dv \, ds\\
	&\le \frac{(14/3)^{d+\alpha}c_{d,\alpha}}{|x-y|^{d+\alpha}}  \int_0^t \P_x (\tau^{\kappa/2}_{U} > s)  ds \,\liminf_{\delta\to 0} \fint_{B(y,\delta)}   dv \\
	&\le \frac{(14/3)^{d+\alpha}c_{d,\alpha}\E_x[\tau^{\kappa/2}_{U}]}{|x-y|^{d+\alpha}}  \le \frac{c_1  \psi(|x|)}{|x-y|^{d+\alpha}} .
\end{align*}
For $I_2$, we observe that for any $\delta\in (0,|x-y|/4)$, by \eqref{e:HKE-alpha-stable}, 
\begin{align*}
	&\sup_{s \in (0,t], \, z \in B(y,|x-y|/2)^c, \, v \in B(y,\delta)} p^{\kappa/2}(s,z,v) \\
	& \le 	\sup_{s \in (0,t], \, z \in  B(y,|x-y|/2)^c, \, v \in B(y,\delta)} q(s,z,v) \le 	c_2\sup_{z \in  B(y,|x-y|/2)^c, \, v \in B(y,\delta)} \frac{t}{|z-v|^{d+\alpha}} 
	\le \frac{c_2t}{(|x-y|/4)^{d+\alpha}}.
\end{align*}
	Using this, Lemma \ref{l:survival} and \eqref{e:psi-scale}, we get that
	\begin{align*}
		I_2 &\le \frac{4^{d+\alpha}c_2t}{|x-y|^{d+\alpha}} \P_x \Big(X^{\kappa/2}_{\tau^{\kappa/2}_U} \in \R^d_0\Big) \limsup_{\delta\to 0}\fint_{B(y,\delta)} dv \\ 
		&\le  \frac{c_3 t \psi(|x|)}{|x-y|^{d+\alpha} \psi(r_0/\eps_0)} \le  \frac{c_4  \psi(|x|)}{|x-y|^{d+\alpha}} \bigg( 1\vee \frac{t}{\psi(|y|)}\bigg).
	\end{align*}
 The proof is complete.  \qed

\begin{lem}\label{l:upper-HKE-2}
For any $T\ge 1$, there exists  $C=C(d,\alpha,\beta_1,\beta_2,\Lambda,T)>0$ such that for all $t \in (0,T]$ and $x,y \in \R^d_0$ with $|x|< \eps_0\psi^{-1}(t/T)$ and  $|x|\le |y|$,
	\begin{equation}\label{e:upper-HKE-2}
	p^{\kappa}(t,x,y)	\le 	p^{\kappa/2}(t,x,y) \le \frac{C\psi(|x|)}{t}\bigg( 1\vee \frac{t}{\psi(|y|)}\bigg)t^{-d/\alpha}\bigg( 1 \wedge \frac{t^{1/\alpha}}{|x-y|}\bigg)^{d+\alpha}.
	\end{equation}
\end{lem}
\pf 
It suffices to prove the second inequality in \eqref{e:upper-HKE-2}.
 When $|x-y|\le 9t^{1/\alpha}$,  \eqref{e:upper-HKE-2}  follows from Lemma \ref{l:upper-HKE-0}. 
 Suppose that $|x-y|>9t^{1/\alpha}$. If $|y|>8|x|$, then \eqref{e:upper-HKE-2} follows from Lemma \ref{l:upper-HKE-1}. If $|y|\le 8|x|$, then we get $|x|\le |y| \le 8|x| <\psi^{-1}(t/T)$ so that by \eqref{e:psi-scale},
\begin{align*}
	\frac{\psi(|x|)}{t} \bigg( 1\vee \frac{t}{\psi(|y|)}\bigg)  = \frac{\psi(|x|)}{\psi(|y|)} \ge c_1.
\end{align*}
Hence, using \eqref{e:HKE-alpha-stable}, we get that
\begin{align*}
	p^{\kappa/2}(t,x,y)\le q(t,x,y)\le \frac{c_2t}{|x-y|^{d+\alpha}} \le \frac{c_1^{-1}c_2\psi(|x|)}{t}  \bigg( 1\vee \frac{t}{\psi(|y|)}\bigg) \frac{t}{|x-y|^{d+\alpha}}.
\end{align*} 
The proof is complete.\qed
 
 \begin{lem}\label{l:upper-HKE-3}
 	Let $T\ge1$. There exist  $\lambda_1=\lambda_1(\beta_2,\Lambda)>0$ and $C= C(d,\alpha,\beta_1,\beta_2,\Lambda,T)>0$ such that for all $t \in (0,T]$ and $x,y \in \R^d_0$ with $|x|< \eps_0\psi^{-1}(t/T)/8$ and  $|x|\le |y|< \eps_0\psi^{-1}(t/T)$,
 	\begin{align*}
 		p^{\kappa}(t,x,y) \le \frac{C\psi(|x|)\psi(|y|)}{t^2} \bigg[ e^{-\lambda_1t/\psi(|y|)}t^{-d/\alpha}\bigg( 1 \wedge \frac{t^{1/\alpha}}{|x-y|}\bigg)^{d+\alpha} +  \frac{t^2}{\psi^{-1}(t)^{d+2\alpha}}  \bigg].
 	\end{align*}
 \end{lem}
 \pf Define  $M_s:=\sup_{0\le u\le s} |X^\kappa_u|$ for $s>0$.  Let $n_0\in \N$ be such that $2^{n_0-1}|y|<\eps_0\psi^{-1}(t/T)\le 2^{n_0}|y|$.  By the lower-semicontinuity of $p^{\kappa}$, we have
 \begin{align*}
& p^{\kappa}(t,x,y) \le \liminf_{\delta\to 0}\fint_{B(y,\delta)}	p^\kappa(t,x,v) dv\\ &\le\limsup_{\delta\to 0} \frac{c_1}{\delta^d}\P_x\left( X^\kappa_t \in B(y,\delta) :  M_t \le 8|y| \right)  + \limsup_{\delta\to 0}\frac{c_1}{\delta^d}\P_x\left(X^\kappa_t \in B(y,\delta) :  M_t > 2^{n_0+3}|y| \right)\\
 &\quad\;\;   +  \sum_{m=1}^{n_0}\limsup_{\delta\to 0} \frac{c_1}{\delta^d}\P_x\left(X^\kappa_t \in B(y,\delta) : 2^{m+2}|y| < M_t \le 2^{m+3}|y| \right)\\
 &=:I_1+I_2+I_3.
 \end{align*}
From  \eqref{e:Feynman-Kac-0}, we get that for any $\delta>0$,
 \begin{align}\label{e:upper-HKE3-1}
 	\P_x\left(X^\kappa_t \in B(y,\delta) :  M_t \le 8|y| \right) &\le 	\P_x(X^{\kappa/2}_t \in B(y,\delta)  )  \exp \left( - \frac{t}{2} \inf_{|z|\le 8|y|} \kappa(z) \right) \nn\\
 	&= \int_{B(y,\delta)}p^{\kappa/2}(t,x,v)dv \exp \left( - \frac{t}{2} \inf_{|z|\le 8|y|} \kappa(z) \right).
 	\end{align}
By \eqref{e:kappa-cond-1} and \eqref{e:psi-scale}, we have 
$$\frac{t}{2}\inf_{|z|\le 8|y|}\kappa(z) \ge \frac{t}{2\Lambda \psi(8|y|)} \ge \frac{c_2t}{\psi(|y|)},$$
for some $c_2=c_2(\beta_2, \Lambda)>0$.
Using this and  Lemma \ref{l:upper-HKE-2} in the first  inequality
below and $\sup_{u \ge 1} u^2e^{-2^{-1}c_2u} <\infty$ in the third,  we get from \eqref{e:upper-HKE3-1} that  
\begin{align*}
I_1&\le \frac{c_3\psi(|x|)e^{-c_2t/\psi(|y|)}}{t}\limsup_{\delta\to 0} \fint_{B(y,\delta)}\bigg( 1\vee \frac{t}{\psi(|v|)}\bigg) \bigg( t^{-d/\alpha} \wedge \frac{t}{|x-v|^{d+\alpha}}\bigg)dv\\
&= \frac{c_3\psi(|x|)e^{-c_2t/\psi(|y|)}}{\psi(|y|)} \bigg( t^{-d/\alpha} \wedge \frac{t}{|x-y|^{d+\alpha}}\bigg)\\
&\le \frac{c_4\psi(|x|)\psi(|y|)e^{-2^{-1}c_2t/\psi(|y|)}}{t^2} \bigg( t^{-d/\alpha} \wedge \frac{t}{|x-y|^{d+\alpha}}\bigg).
\end{align*}
For $I_2$,  since $2^{n_0}|y|\ge \eps_0\psi^{-1}(t/T)$, using the strong Markov property, we see that for any $\delta>0$,
\begin{align*}
	&\P_x\left(X^\kappa_t \in B(y,\delta) :  M_t >2^{n_0+3}|y| \right)   \le  \P_x\left(X^\kappa_t \in B(y,\delta) :  \tau^{\kappa}_{B(0,8\eps_0\psi^{-1}(t/T))} <t \wedge \zeta^\kappa \right) \\
	&=\E_x\left[ \P_{X^\kappa_{\tau^{\kappa}_{B(0,8\eps_0\psi^{-1}(t/T))}}}\Big(X^\kappa_{t-\tau^{\kappa}_{B(0,8\eps_0\psi^{-1}(t/T))}} \in B(y,\delta) \Big) :  \tau^{\kappa}_{B(0,8\eps_0\psi^{-1}(t/T))} <t \wedge \zeta^\kappa \right]\\
	&\le \P_x\left( \tau^{\kappa}_{B(0,8\eps_0\psi^{-1}(t/T))} <t \wedge \zeta^\kappa \right) \sup_{0<s\le t, \, z \in B(0,8\eps_0\psi^{-1}(t/T))^c} \P_z \left(X^\kappa_{s} \in B(y,\delta) \right).
\end{align*}
Hence, using  \eqref{e:psi-scale} and Lemmas   \ref{l:EP} and \ref{l:upper-HKE-2}, since $|y| <\eps_0\psi^{-1}(t/T)$,  we obtain
\begin{align*}
	I_2&\le \frac{c_5\psi(|x|)}{\psi^{-1}(t/T)^\alpha}\limsup_{\delta\to0} \sup_{0<s\le t, \, z \in B(0,8\eps_0\psi^{-1}(t/T))^c} \fint_{B(y,\delta)} p^\kappa(s,z,v)dv\\
	&\le \frac{c_6\psi(|x|)}{\psi^{-1}(t/T)^\alpha} \sup_{0<s\le t, \, z \in B(0,8\eps_0\psi^{-1}(t/T))^c}  \frac{\psi(|y|)}{s}\bigg( 1\vee \frac{s}{\psi(|z|)}\bigg)  s^{-d/\alpha} \bigg(1\wedge \frac{s^{1/\alpha}}{|y-z|}\bigg)^{d+\alpha}\\
		&\le \frac{c_7\psi(|x|)}{\psi^{-1}(t/T)^\alpha} \sup_{0<s\le t, \, z \in B(0,8\eps_0\psi^{-1}(t/T))^c}  \frac{\psi(|y|)}{s} \frac{s}{|y-z|^{d+\alpha}}\\
			&\le \frac{c_7\psi(|x|)\psi(|y|)}{\psi^{-1}(t/T)^\alpha (7\eps_0\psi^{-1}(t/T))^{d+\alpha}} \le \frac{c_8\psi(|x|)\psi(|y|)}{\psi^{-1}(t)^{d+2\alpha}}.
\end{align*}
It remains to estimate $I_3$. Let $1\le  m\le n_0$. Using the formula  \eqref{e:Feynman-Kac-0} in the first inequality below, the strong Markov property, \eqref{e:kappa-cond-1} and \eqref{e:psi-scale} in the second, and Lemma \ref{l:EP} and $\psi(2^{m+2}|y|)\le \psi(8\eps_0\psi^{-1}(t/T))\le t/T$  in the fourth, we obtain
\begin{align}\label{e:UHK-I3-1}
	&\P_x\left(X^\kappa_t \in B(y,\delta) : 2^{m+2}|y| < M_t \le 2^{m+3}|y| \right)\nn\\
	&\le \P_x\left( X^{\kappa/2}_t \in B(y,\delta) : 2^{m+2}|y| < M_t \le 2^{m+3}|y| \right)  \exp \left( - \frac{t}{2} \inf_{|z|\le 2^{m+3}|y|} \kappa(z) \right)\nn\\
	&\le \E_x\bigg[ \P_{X^{\kappa/2}_{\tau^{\kappa/2}_{B(0,2^{m+2}|y|)}}}\Big(X^{\kappa/2}_{t-\tau^{\kappa/2}_{B(0,2^{m+2}|y|)}} \in B(y,\delta) \Big)  :  \tau^{\kappa/2}_{B(0,2^{m+2}|y|)} <t \wedge \zeta^{\kappa/2} \bigg] \exp \left( - \frac{c_9t}{ \psi(2^{m+2}|y|)}\right)\nn\\
		&\le\P_x\left(\tau^{\kappa/2}_{B(0,2^{m+2}|y|)} <t \wedge \zeta^{\kappa/2} \right) \sup_{0<s\le t, \, z \in B(0,2^{m+2}|y|)^c} \P_z \left(X^{\kappa/2}_{s} \in B(y,\delta) \right)   \exp \left( - \frac{c_9t}{ \psi(2^{m+2}|y|)}\right)\nn\\
		&\le  \frac{c_{10}\psi(|x|)t}{(2^{m+2}|y|)^{\alpha}\psi(2^{m+2}|y|)} \sup_{0<s\le t, \, z \in B(0,2^{m+2}|y|)^c} \int_{B(y,\delta)} p^{\kappa/2}(s,z,v)dv  \exp \left( - \frac{c_9t}{ \psi(2^{m+2}|y|)}\right).
\end{align}
By Lemma \ref{l:upper-HKE-2} and the inequality $\psi(2^{m+2}|y|)\le t/T$, we have
\begin{align}\label{e:UHK-I3-2}
&\limsup_{\delta\to0} \frac{1}{\delta^d}\sup_{0<s\le t, \, z \in B(0,2^{m+2}|y|)^c} \int_{B(y,\delta)} p^{\kappa/2}(s,z,v)dv\nn\\
&\le \sup_{0<s\le t, \, z \in B(0,2^{m+2}|y|)^c}  \frac{c_{11}\psi(|y|)}{s}\bigg( 1\vee \frac{s}{\psi(|z|)}\bigg)  \frac{s}{|y-z|^{d+\alpha}}  \le \frac{ c_{11}\psi(|y|)t}{(2^{m+1}|y|)^{d+\alpha} \psi(2^{m+2}|y|)}.
\end{align}
Combining \eqref{e:UHK-I3-1} 
with \eqref{e:UHK-I3-2},  and using $\sup_{a \ge 1} a^{3+(d+2\alpha)/\beta_1}e^{-c_9a}<\infty$ and \eqref{e:psi-scale}, we deduce that
\begin{align*}
	&\limsup_{\delta\to0} \frac{1}{\delta^d}\P_x\left(X^\kappa_t \in B(y,\delta) : 2^{m+2}|y| < M_t \le 2^{m+3}|y| \right) \\
	&\le  \frac{2^{d+\alpha}c_{12}\psi(|x|)\psi(|y|)t^2}{(2^{m+2}|y|)^{d+2\alpha}\psi(2^{m+2}|y|)^2}  \exp \left( - \frac{c_9t}{ \psi(2^{m+2}|y|)}\right)\\
	&\le  \frac{c_{13}\psi(|x|)\psi(|y|)t^2}{(2^{m+2}|y|)^{d+2\alpha}\psi(2^{m+2}|y|)^2}  \bigg( \frac{\psi(2^{m+2}|y|)}{t}\bigg)^{3+(d+2\alpha)/\beta_1} \\
		&\le  \frac{c_{14}\psi(|x|)\psi(|y|) \psi(2^{m+2}|y|)}{t\psi^{-1}(t)^{d+2\alpha}}.
\end{align*}
Using this, \eqref{e:psi-scale} and $2^{n_0+2}|y| < 8\eps_0\psi^{-1}(t/T) \le \psi^{-1}(t/T)$, we conclude that
\begin{align*}
	I_2&\le \frac{c_{14}\Lambda\psi(|x|)\psi(|y|) \psi(2^{n_0+2}|y|)}{t\psi^{-1}(t)^{d+2\alpha}} \sum_{m=1}^{n_0} 2^{-(n_0-m)\beta_1} \le \frac{c_{15}\psi(|x|)\psi(|y|)}{\psi^{-1}(t)^{d+2\alpha}}.
\end{align*}
The proof is complete. \qed

Note that for any $t \in (0,T]$, by \eqref{e:psi-scale}, \begin{align}\label{e:psi-inverse}
	\psi^{-1}(t/T)\ge \psi^{-1}(1) (t/(\Lambda T))^{1/\beta_1} = (t/(\Lambda T))^{1/\beta_1} \ge (t/(\Lambda T))^{1/\alpha}.
\end{align}

\noindent \textbf{Proof of  Theorem \ref{t:smalltime} (Upper estimates).}  Let $t\in (0,T]$ and $x,y \in \R^d_0$. Without loss of generality,  we assume that $|x|\le |y|$.  We deal with three cases separately.

\smallskip

Case 1: $|x| \ge \eps_0\psi^{-1}(t/T)/8$. By \eqref{e:HKE-alpha-stable}, we get
\begin{align*}
	p^{\kappa}(t,x,y)\le q(t,x,y)\le c_1 t^{-d/\alpha}\bigg( 1 \wedge \frac{t^{1/\alpha}}{|x-y|}\bigg)^{d+\alpha} .
\end{align*}
Further, by \eqref{e:psi-scale}, we have $\psi(|y|)/t \ge \psi(|x|)/t \ge \psi(\eps_0\psi^{-1}(t/T)/8)/t \ge c_2$. Thus, \eqref{e:smalltime-upper} holds.

\smallskip

Case 2: $|x|<\eps_0\psi^{-1}(t/T)/8$ and $|y| \ge \eps_0\psi^{-1}(t/T)$. In this case, by \eqref{e:psi-scale}, we have $\psi(|x|)/t\le c_3$ and $\psi(|y|)/t\ge c_4$. Moreover, by \eqref{e:psi-inverse},
\begin{align*}
	|x-y| \ge 7\eps_0\psi^{-1}(t/T)/8  \ge 7\eps_0 (t/(\Lambda T))^{1/\alpha}/8.
\end{align*}
Hence, the right-hand side of \eqref{e:smalltime-upper} is bounded below by
\begin{align*}
	\frac{c_5\psi(|x|)}{t^{1+d/\alpha}} \bigg(  \frac{t^{1/\alpha}}{|x-y|}\bigg)^{d+\alpha} = \frac{c_5\psi(|x|)}{|x-y|^{d+\alpha}}.
\end{align*}
Therefore, from  Lemma \ref{l:upper-HKE-1}, since $\psi(|y|)/t\ge c_4$, we conclude that \eqref{e:smalltime-upper} holds in this case.
 
 \smallskip
 
 Case 3: $|x|<\eps_0\psi^{-1}(t/T)/8$ and $|y| < \eps_0\psi^{-1}(t/T)$. In this case, since $\psi(|x|)\le \psi(|y|)\le t$ and $|x-y| \le |x|+|y|<\psi^{-1}(t)$, the right-hand side of \eqref{e:smalltime-upper} is equal to
 \begin{align*}
 	\frac{c_6\psi(|x|)\psi(|y|)}{t^2} \bigg[ \,e^{-\lambda_1t/\psi( |y|)}t^{-d/\alpha}\bigg( 1 \wedge \frac{t^{1/\alpha}}{|x-y|}\bigg)^{d+\alpha} +  \frac{t^2}{\psi^{-1}(t)^{d+2\alpha}}\,\bigg].
 \end{align*}
 Thus, \eqref{e:smalltime-upper} follows from Lemma \ref{l:upper-HKE-3}. The proof is complete. \qed

 \subsection{Small time lower heat kernel estimates}\label{ss:smalltime-lower}
 
 \begin{lem}\label{l:lower-HKE-1}
 For any $T\ge 1$,	there exists   $C=C(d,\alpha,\beta_1,\beta_2,\Lambda,T)>0$  such that for all $t \in (0,T]$ and $x,y \in \R^d_0$ with $|x|\le \psi^{-1}(t/T)$ and $|y|\ge 2\psi^{-1}(t/T)$,
 	\begin{align*}
 		p^\kappa(t,x,y) \ge \frac{C\psi(|x|)}{|x-y|^{d+\alpha}}.
 	\end{align*}
 \end{lem} 
 \pf Fix $t\in (0,T]$ and $x,y \in \R^d_0$ such that  $|x|\le  \psi^{-1}(t/T)$ and $|y|\ge 2\psi^{-1}(t/T)$. Let $U:=B(0,3|y|) \setminus B(0, |x|/2)$ and $r_t:=(t/(\Lambda T))^{1/\alpha}$. By the strong Markov property and the joint continuity of $p^{\kappa,U}$, we have
 \begin{align}\label{e:lower-HKE-1}
 	p^{\kappa,U}(t,x,y) &= \E_x\left[ p^{\kappa,U}(t-\tau^\kappa_{B(x,|x|/2)} , X^\kappa_{ \tau^\kappa_{B(x,|x|/2)}} , y) : \tau^\kappa_{B(x,|x|/2)} <t \right]\nn\\
 	&\ge \E_x\left[ p^{\kappa,U}(t-\tau^\kappa_{B(x,|x|/2)} , X^\kappa_{ \tau^\kappa_{B(x,|x|/2)}} , y) : \tau^\kappa_{B(x,|x|/2)} \le t/2, \,  X^\kappa_{ \tau^\kappa_{B(x,|x|/2)}} \in B(y,r_t)\right] \nn\\
 	& \ge  \P_x\left( \tau^\kappa_{B(x,|x|/2)} \le t/2, \; X^\kappa_{ \tau^\kappa_{B(x,|x|/2)}} \in B(y, r_t)\right) \inf_{s \in[t/2,t], \, z \in B(y,r_t)} p^{\kappa,U}(s,z,y). 
 \end{align}
 For all $s \in[t/2,t]$ and $z \in B(y,r_t)$, 
 by \eqref{e:psi-inverse},  we have   $|y| \wedge |z| > 2\psi^{-1}(t/T)-r_t \ge \psi^{-1}(t/T)$ and $|y| \vee |z| < |y|  + r_t \le 3|y|/2$. Hence, using Proposition \ref{p:interior-HKE-lower} (with  $R=3|y|$,  $r= \psi^{-1}(t/T)/2$ and   $a=2^{\beta_2}\Lambda T$) and \eqref{e:HKE-alpha-stable}, we get
 \begin{align}\label{e:lower-HKE-2}
& \inf_{s \in[t/2,t], \, z \in B(y,r_t)} p^{\kappa,U}(s,z,y) \ge   \inf_{s \in[t/2,t], \, z \in B(y,r_t)} p^{\kappa, B(0,3|y|) \setminus B(0, \psi^{-1}(t/T)/2)}(s,z,y) \nn\\
 	&\qquad\ge c_1 \inf_{s\in [t/2,t]} \inf_{z \in B(y,r_t)} q(s,y,z) \ge c_2 \inf_{s\in [t/2,t]} s^{-d/\alpha}= c_2 t^{-d/\alpha}.
 \end{align}
Besides, we  note that for all $w \in B(x,|x|/2)$ and $z \in B(y, r_t)$, since $|y|\ge 2|x|$ and $|y|\ge 2r_t$,
 \begin{align}\label{e:lower-HKE-off-1}
 	|w-z| \le 2|x| + |y| + r_t \le  5|y|/2 \le   5(|y|-|x|) \le 5|y-x|.
 \end{align}
 Further, by Lemma \ref{l:exit-time-lower} and \eqref{e:psi-scale}, 
 \begin{align}\label{e:lower-HKE-off-2}
 	\E_x\big[ \tau^\kappa_{B(x,|x|/2)}  \wedge (t/2)\big] \ge  \big(\psi(|x|/2) \wedge (t/2)\big) \,\P_x \big( \tau^\kappa_{B(x,|x|/2)} \ge \psi(|x|/2)\big) \ge c_3 \psi(|x|).
 \end{align}
 Using \eqref{e:Levysystem} in the equality below, \eqref{e:lower-HKE-off-1} in the first inequality and \eqref{e:lower-HKE-off-2} in the second inequality, we obtain
 \begin{align*}
 	 &\P_x\left( \tau^\kappa_{B(x,|x|/2)} \le t/2, \; X^\kappa_{ \tau^\kappa_{B(x,|x|/2)}} \in B(y, r_t)\right)	= \E_x\left[ \int_0^{\tau^\kappa_{B(x,|x|/2)}  \wedge (t/2)} \int_{B(y,r_t)} 
	 \frac{\sA(d,-\alpha)}{|X^\kappa_s - z|^{d+\alpha}}dz ds\right] \\
	&\ge \frac{\sA(d,-\alpha)
	\E_x[ \tau^\kappa_{B(x,|x|/2)}  \wedge (t/2)]}{5^{d+\alpha}  |x-y|^{d+\alpha}}  \int_{B(y,r_t)}dz 
	\ge \frac{c_4t^{d/\alpha} \psi(|x|) }{|x-y|^{d+\alpha}}.
 \end{align*}
 Combining this with \eqref{e:lower-HKE-1} and \eqref{e:lower-HKE-2}, and using the inequality $p^\kappa(t,x,y) \ge p^{\kappa,U}(t,x,y)$, we arrive at the result.
 \qed

 \begin{lem}\label{l:lower-HKE-2}
 For any $T\ge 1$,	there exists   $C=C(d,\alpha,\beta_1,\beta_2,\Lambda,T)>0$  such that for all $t \in (0,T]$ and $x,y \in \R^d_0$ with $|x|\le \psi^{-1}(t/(2T))$ and $|x|\le |y|\le  2\psi^{-1}(t/T)$,
 	\begin{align}\label{e:lower-HKE2}
 		p^{\kappa}(t,x,y) &\ge \frac{C\psi(|x|)\psi(|y|)}{\psi^{-1}(t/T)^{d+2\alpha}}.
 	\end{align}
 \end{lem} 
 \pf  Let $z_0 \in \R^d$ be such that $|z_0|=3\psi^{-1}(t/T)$. By Lemma \ref{l:lower-HKE-1},  for any $z \in B(z_0,\psi^{-1}(t/T))$,
 \begin{align}\label{e:lower-HKE-2-1}
 	p^\kappa(t/2,x,z)  \ge \frac{c_1\psi(|x|)}{|x-z|^{d+\alpha}}\ge  \frac{c_1\psi(|x|)}{(5\psi^{-1}(t/T))^{d+\alpha}}.
 \end{align}
Similarly, by  Lemma \ref{l:lower-HKE-1} and the symmetry of $p^\kappa(t,\cdot,\cdot)$, if $|y|\le \psi^{-1}(t/(2T))$, then  $p^{\kappa}(t/2,z,y) \ge c_2\psi(|y|)/\psi^{-1}(t/T)^{d+\alpha}$  for any $z \in B(z_0,\psi^{-1}(t/T))$. If $ \psi^{-1}(t/(2T)) < |y|\le 2\psi^{-1}(t/T)$, then  by using \eqref{e:psi-scale},  Proposition \ref{p:interior-HKE-lower} (with $R=8\psi^{-1}(t/T)$, 
$r=\psi^{-1}(t/(2T))/2$ and $a=2^{\beta_2}\Lambda T$), 
\eqref{e:HKE-alpha-stable} and \eqref{e:psi-inverse}, we get that for any $z \in B(z_0,\psi^{-1}(t/T))$,
\begin{align*}
	&p^\kappa(t/2,z,y) \ge p^{\kappa, B(0,8\psi^{-1}(t/T)) \setminus B(0, \psi^{-1}(t/(2T))/2)}(t/2,z,y) \\
	&\ge c_3q(t/2,z,y) \ge \frac{c_4t/2}{(6\psi^{-1}(t/T))^{d+\alpha}} \ge \frac{c_5\psi(|y|)}{\psi^{-1}(t/T)^{d+\alpha}}.
\end{align*}
Hence, in both cases, we have
  \begin{align}\label{e:lower-HKE-2-2}
 	p^\kappa(t/2,z,y)  \ge \frac{c_6\psi(|y|)}{\psi^{-1}(t/T)^{d+\alpha}} \quad \text{for all} \;\, z \in B(z_0,\psi^{-1}(t/T)).
 \end{align}
Using the semigroup property, \eqref{e:lower-HKE-2-1} and  \eqref{e:lower-HKE-2-2}, we conclude that
\begin{align*}
	p^\kappa(t,x,y) \ge \int_{B(z_0,\psi^{-1}(t/T))} p^\kappa(t/2,x,z) p^{\kappa}(t/2,z,y)dz \ge \frac{c_5\psi(|x|)\psi(|y|)}{\psi^{-1}(t/T)^{d+2\alpha}}.
\end{align*} 
\qed

\begin{lem}\label{l:lower-HKE-3-case2}
	Let $T\ge1$. There exist  $\lambda_2=\lambda_2(\beta_2,\Lambda)>0$ and $C= C(d,\alpha,\beta_1,\beta_2,\Lambda,T)>0$ such that for all $t \in (0,T]$ and $x,y \in \R^d_0$  with $|y| \ge |x| \vee (t/(2\Lambda T))^{1/\alpha}$ and $|x-y|\le (t/(2\Lambda T))^{1/\alpha}/4$,
	\begin{align*}
		p^{\kappa}(t,x,y) &\ge C t^{-d/\alpha} e^{-\lambda_2t/\psi( |y|)}.
	\end{align*}
\end{lem} 
\pf  From  \eqref{e:Feynman-Kac}, using \eqref{e:kappa-cond-1} and \eqref{e:psi-scale},  we get that
\begin{align}\label{e:lower-HKE3-2}
	&p^{\kappa}(t,x,y) \ge p^{\kappa, B(y,|y|/2)}(t,x,y) \ge  q^{B(y,|y|/2)}(t,x,y) \exp\Big( - t \sup_{z \in B(y,|y|/2)} \kappa(z) \Big) \nn\\
	& \ge  q^{B(y,|y|/2)}(t,x,y) \exp\Big( - t \sup_{|z| \ge |y|/2} \kappa(z) \Big) \ge   q^{B(y,|y|/2)}(t,x,y)  e^{-\lambda_2 t/\psi(|y|)},
\end{align}
for some constant $\lambda_2=\lambda_2(\beta_2, \Lambda)>0$.
On the other hand, since $|y| \ge 4|x-y|$ and $|y| \ge (t/(2\Lambda T))^{1/\alpha}$, applying Proposition \ref{p:DHKE-1} with $k=2(2\Lambda T)^{1/\alpha}$, we obtain
\begin{align*}
	q^{B(y,|y|/2)}(t,x,y) &\ge c_1 \bigg( 1 \wedge \frac{\delta_{B(y,|y|/2)}(x)^\alpha}{t}\bigg)^{1/2}\bigg( 1 \wedge \frac{\delta_{B(y,|y|/2)}(y)^\alpha}{t}\bigg)^{1/2} t^{-d/\alpha}  \bigg( 1 \wedge \frac{t^{1/\alpha}}{|x-y|}\bigg)^{d+\alpha}\\
	&=c_1 \bigg( 1 \wedge \frac{(|y|/2-|x-y|)^\alpha}{t}\bigg)^{1/2}\bigg( 1 \wedge \frac{(|y|/2)^\alpha}{t}\bigg)^{1/2} t^{-d/\alpha}
	\ge c_2t^{-d/\alpha}.
\end{align*}
Combining this with \eqref{e:lower-HKE3-2},  we arrive at  the result.\qed

 \begin{lem}\label{l:lower-HKE-3}
 For any $T \ge 1$, there exists   $C= C(d,\alpha,\beta_1,\beta_2,\Lambda,T)>0$ such that for all $t \in (0,T]$ and $x,y \in \R^d_0$ with  $|x| \le  |y|\le 2\psi^{-1}(t/T)$,
 	\begin{align}\label{e:lower-HKE3}
 	p^{\kappa}(t,x,y) &\ge \frac{C\psi(|x|)\psi(|y|)}{t^{2+d/\alpha}} e^{-\lambda_2t/\psi( |y|)}\bigg( 1 \wedge \frac{t^{1/\alpha}}{|x-y|}\bigg)^{d+\alpha},
 \end{align}
 where $\lambda_2=\lambda_2(\beta_2,\Lambda)>0$ is the constant in Lemma \ref{l:lower-HKE-3-case2}.  \end{lem} 
\pf  Let $b_t:=(t/(2\Lambda T))^{1/\alpha}$   and  let $\lambda_2>0$ be the constant in Lemma \ref{l:lower-HKE-3-case2}. 
Note that $b_t\le \psi^{-1}(t/(2T))$ by \eqref{e:psi-inverse}. We consider the following three cases separately.

\smallskip

Case 1: $|y|\le  b_t$.  By \eqref{e:psi-scale}, it holds that 
\begin{align*}
	\frac{t}{\psi(|y|)} \ge \frac{t}{\Lambda |y|^{\beta_1}} \ge c_1t^{-(\beta_1-\alpha)/\alpha}.
\end{align*}
Using this and  the inequality $\sup_{u\in (0,1]} u^{(d+2\alpha)/\beta_2-2-d/\alpha} e^{-c_1\lambda_2u^{-(\beta_1-\alpha)/\alpha}}=c_2<\infty$, we see that 
\begin{align}\label{e:lower-HKE3-1}
\frac{\psi(|x|)\psi(|y|)}{t^{2+d/\alpha}}	e^{-\lambda_2 t/\psi(|y|)} \le \frac{\psi(|x|)\psi(|y|)}{t^{2+d/\alpha}} e^{-c_1\lambda_2t^{-(\beta_1-\alpha)/\alpha}} \le \frac{c_2\psi(|x|)\psi(|y|)}{t^{(d+2\alpha)/\beta_2}}.
\end{align}
On the other hand, since $|x|\le |y|\le b_t \le \psi^{-1}(t/(2T))$,  from Lemma \ref{l:lower-HKE-2} and \eqref{e:psi-scale}, we obtain
\begin{align*}
	p^\kappa (t,x,y) \ge  \frac{c_3\psi(|x|)\psi(|y|)}{\psi^{-1}(t/T)^{d+2\alpha}} \ge   \frac{c_4\psi(|x|)\psi(|y|)}{(t/T)^{(d+2\alpha)/\beta_2}}.
\end{align*}
Combining this with \eqref{e:lower-HKE3-1}, we conclude that \eqref{e:lower-HKE3} holds in this case.

 \smallskip
 
 Case 2: $|y|>    b_t$ and  $|x-y|\le  b_t/4$. Since $\psi(|x|)\psi(|y|)/t^2 \le (2^{\beta_2}\Lambda)^2$ by \eqref{e:psi-scale}, the result follows from Lemma \ref{l:lower-HKE-3-case2}.

 \smallskip
 
 Case 3:  $|y|> b_t$ and  $|x-y|> b_t/4$.   Let 
 $$
 z_0:= \bigg(\frac{|y|+2^{-3-1/\alpha}b_t}{|y|}\bigg) y.
 $$
 For any $z \in B(z_0,2^{-3-1/\alpha}b_t)$, we have $|z| \ge |y|\ge|x|$,
 \begin{align}\label{e:lower-HKE3-dist}
 	|y-z| \le |y-z_0| + |z_0-z| < 2^{-2-1/\alpha}b_t =  (t/(4\Lambda T))^{1/\alpha}/4
 \end{align}
   and  $|x-z| \le |x-y| + |y-z|< 2|x-y|.$ Hence, applying Proposition \ref{p:interior-HKE-lower}     (with $a= 2^{\beta_2+1}\Lambda$) 
   and using \eqref{e:HKE-alpha-stable}, we get that for all $z \in B(z_0,2^{-3-1/\alpha}b_t)$,
 \begin{align}\label{e:lower-HKE3-4}
 	p^\kappa(\psi(|x|)/2,x,z)&\ge 	p^{\kappa, B(0,2(|z_0|+b_t))\setminus B(0,|x|/2)}(\psi(|x|)/2,x,z)\nn\\
 	&\ge c_5q(\psi(|x|)/2,x,z) \ge \frac{c_6\psi(|x|)/2}{(2|x-y|)^{d+\alpha}}.
 \end{align}
 Besides, for any  $z \in B(z_0,2^{-3-1/\alpha}b_t)$, since $t-\psi(|x|)/2 \in [t/2,t]$, $|z| \ge |y|>b_t$ and \eqref{e:lower-HKE3-dist} holds,   by Lemma \ref{l:lower-HKE-3-case2}, we obtain
 \begin{equation}\label{e:lower-HKE3-5}
 	p^{\kappa} ( t- \psi(|x|)/2 , z,y ) \ge  c_7 (t-\psi(|x|)/2)^{-d/\alpha} e^{-\lambda_2 (t-\psi(|x|/2))/\psi(|z|)}  \ge  c_7 t^{-d/\alpha} e^{-\lambda_2 t/\psi(|y|)}.
 \end{equation}
Using the semigroup property, \eqref{e:lower-HKE3-4} and \eqref{e:lower-HKE3-5}, we arrive at
\begin{align*}
	p^\kappa(t,x,y) &\ge \int_{ B(z_0,2^{-3-1/\alpha}b_t)} p^\kappa(\psi(|x|)/2,x,z) p^{\kappa}(t-\psi(|x|)/2,z,y)dz\\
		&\ge \frac{c_{8}\psi(|x|)t^{-d/\alpha} e^{-\lambda_2 t/\psi(|y|) } }{|x-y|^{d+\alpha}} \int_{ B(z_0,2^{-3-1/\alpha}b_t)}  dz\\
		& =  \frac{c_{9}\psi(|x|)e^{-\lambda_2 t/\psi(|y|) } }{|x-y|^{d+\alpha}}  \ge  \frac{c_{9}\psi(|x|)\psi(|y|)e^{-\lambda_2 t/\psi(|y|) } }{t|x-y|^{d+\alpha}} .
\end{align*}
 The proof is complete. \qed

 \medskip
 
 \noindent \textbf{Proof of  Theorem \ref{t:smalltime} (Lower  estimates).} Let $t\in (0,T]$ and $x,y \in \R^d_0$. Without loss of generality, by symmetry, we assume that $|x|\le |y|$.  We consider three cases separately.
 
 \smallskip
 
 Case 1: $|x| \ge \psi^{-1}(t/(2T))$. In this case, we have  $t \le 2T\psi(|x|) \le 2^{\beta_2+1}\Lambda T\psi(|x|/2)$. Using Proposition \ref{p:interior-HKE-lower}   (with $a=2^{\beta_2+1}\Lambda T$) 
 and \eqref{e:HKE-alpha-stable}, we obtain
 \begin{align*}
 	p^\kappa(t,x,y) \ge p^{\kappa, B(0,3|y|)\setminus B(0,|x|/2)}(t,x,y) \ge c_1q(t,x,y)  \ge c_2 t^{-d/\alpha}\bigg( 1 \wedge \frac{t^{1/\alpha}}{|x-y|}\bigg)^{d+\alpha}.
 \end{align*} 
 Combining this with Lemma \ref{l:dominant}, we conclude that \eqref{e:smalltime-lower} holds.
 
 \smallskip
 
 Case 2:  $|x| < \psi^{-1}(t/(2T))$ and $|y|\ge 2\psi^{-1}(t/T)$. By \eqref{e:psi-inverse}, $|x-y| \ge \psi^{-1}(t/T) \ge c_3t^{1/\alpha}$ in this case. Hence, by Lemma \ref{l:lower-HKE-1}, we have
 \begin{align*}
 	p^{\kappa}(t,x,y) \ge \frac{c_3\psi(|x|)}{|x-y|^{d+\alpha}} \ge  \frac{c_4\psi(|x|)}{t} t^{-d/\alpha}\bigg( 1 \wedge \frac{t^{1/\alpha}}{|x-y|}\bigg)^{d+\alpha}.
 \end{align*}
Using  Lemma \ref{l:dominant} again, we obtain \eqref{e:smalltime-lower}.

 \smallskip

Case 3:  $|x| < \psi^{-1}(t/(2T))$ and $|y|< 2\psi^{-1}(t/T)$.  Note that $|x-y|\le |x|+|y|<3\psi^{-1}(t)$. Thus, from Lemmas \ref{l:lower-HKE-2} and \ref{l:lower-HKE-3}, \eqref{e:smalltime-lower} follows. 

\smallskip

The proof is complete. \qed

\section{Key proposition for large time estimates}\label{s:5}

Starting from this section, we assume that $\kappa \in \sK_{\alpha}(\psi,\Lambda)$.
Recall that  $\wt q(t,x,y)$ is defined by \eqref{e:def-wtq}. By \cite[(9)]{BJ07}, the following  3P inequality holds for $\wt q(t,x,y)$.

\begin{prop}\label{p:3P}
	There exists $C=C(d,\alpha)>0$ such that for all $t>s>0$ and $x,y,z \in \R^d$,
	\begin{align*}
		\frac{\wt q(t-s,x,z) \wt q(s,z,y)}{\wt q(t,x,y)} \le C( \wt q(t-s,x,z) + \wt q(s,y,z) ).
	\end{align*}
\end{prop}

The next proposition is the main result of this section. 
The proof will be provided at the end of this section.
\begin{prop}\label{p:one-step}
 There exists $C=C(d,\alpha,\beta_1,\Lambda)>0$ such that for all $R\ge 1$,  $t\ge R^\alpha$ and  $x,y \in B(0,2R)^c$,
 \begin{align}\label{e:one-step}
& \int_0^t \int_{B(0,R)^c} \frac{ q^{B(0,R)^c}(s,x,z) q^{B(0,R)^c}(t-s,z,y)}{q^{B(0,R)^c}(t,x,y) }\kappa(z)dzds \le \frac{C}{R^{\beta_1-\alpha}}.
 \end{align}
\end{prop}

 In the remainder of this section, we  fix  $R\ge 1$. Note that  by \eqref{e:kappa-cond-3} and  \eqref{e:psi-scale},
 \begin{align}\label{e:psi-bound}\kappa(z)\le \frac{\Lambda}{\psi(|z|)} \le \frac{\Lambda^2}{|z|^{\beta_1}}	 \quad \text{for all} \;\, z \in B(0,R)^c.
 \end{align}
 Further, we see that
\begin{align}\label{e:one-step-range}
 \delta_{B(0,R)^c}(z)  \ge |z|/2 \ge R \quad \text{for all} \;\, z \in B(0,2R)^c.
\end{align}

We  consider the cases $d>\alpha$, $d=1<\alpha$ and $d=\alpha=1$ separately.

\subsection{The case of $d>\alpha$}

When $d>\alpha$, we have  for all $x,y \in \R^d$ with $x \neq y$,
\begin{align}\label{e:Green-alpha-stable}
	\int_0^\infty  \wt q(s,x,y) ds = \int_0^{|x-y|^\alpha} \frac{s}{|x-y|^{d+\alpha}} ds + \int_{|x-y|^\alpha}^\infty s^{-d/\alpha}ds =  \frac{d+\alpha}{2(d-\alpha)|x-y|^{d-\alpha}}.
\end{align}

\begin{lem}\label{l:one-step-transient-key}
	Suppose that $d>\alpha$. There exists $C=C(d,\alpha,\beta_1)>0$ independent of $R$ such that
	\begin{align*}
		\int_{B(0,R)^c} \frac{dz}{|w-z|^{d-\alpha}|z|^{\beta_1}} \le  \frac{C}{R^{\beta_1-\alpha}} \quad \text{for all} \;\, w \in B(0,2R)^c.
	\end{align*}
\end{lem}
\pf Let $w \in B(0,2R)^c$.  We decompose the integral as follows:
\begin{align*}
	&\int_{B(0,R)^c} \frac{dz}{|w-z|^{d-\alpha}|z|^{\beta_1}}\\
	& = \bigg(	\int_{B(w,|w|/2) \setminus B(0,R)} +	\int_{B(w,3|w|)\setminus (B(w,|w|/2) \cup B(0,R))} + 	\int_{B(w,3|w|)^c} \bigg) \frac{dz}{|w-z|^{d-\alpha}|z|^{\beta_1}}\\
	&=:I_1+I_2+I_3.
\end{align*}
For any $z \in B(w,|w|/2)$, we have $|z|\ge |w|/2$. Using this, we get that
\begin{align*}
	I_1& \le \frac{1}{(|w|/2)^{\beta_1}}\int_{B(w,|w|/2)} \frac{dz}{|w-z|^{d-\alpha}} = \frac{A_{d-1} (|w|/2)^\alpha}{\alpha(|w|/2)^{\beta_1}}  = \frac{c_1}{|w|^{\beta_1-\alpha}} \le  \frac{c_1}{(2R)^{\beta_1-\alpha}}.
\end{align*}
Set $\wt \beta_1:= \beta_1 \wedge ((d+\alpha)/2)$. Then $\alpha<\wt \beta_1<d$. For $I_2$, we have
\begin{align*}
	I_2& \le \frac{1}{(|w|/2)^{d-\alpha}}	\int_{B(w,3|w|)\setminus (B(w,|w|/2) \cup B(0,R))}\frac{dz}{|z|^{\beta_1}}\nn\\
	& \le \frac{1}{R^{\beta_1- \wt \beta_1}(|w|/2)^{d-\alpha}}	\int_{B(0,4|w|)\setminus B(0,R)} \frac{dz}{|z|^{\wt\beta_1}}\nn\\
	&\le  \frac{c_2 (4|w|)^{d-\wt \beta_1}}{R^{\beta_1- \wt \beta_1}|w|^{d-\alpha}}  \le  \frac{c_3}{R^{\beta_1-\alpha}} .
\end{align*}
For any $z \in B(w,3|w|)^c$, we have $|z| \ge 2|w|$ and $|w-z| \ge |z|-|w| \ge |z|$. Hence,
\begin{align*}
	I_3&\le \int_{B(0,2|w|)^c}  \frac{dz}{|z|^{d-\alpha+\beta_1}} = \frac{c_4}{(2|w|)^{\beta_1-\alpha}} \le  \frac{c_4}{(4R)^{\beta_1-\alpha}}.
\end{align*}
The proof is complete. \qed

\begin{lem}\label{l:one-step-transient}
When $d>\alpha$, \eqref{e:one-step} holds true.
\end{lem}
\pf  Let $t \ge R^\alpha$ and $x,y \in B(0,2R)^c$. By Proposition \ref{p:DHKE-2}(i), we have
\begin{align}\label{e:one-step-1}
	q^{B(0,r)^c}(t,x,y) \ge c_1\wt q(t,x,y).
\end{align}
On the other hand, using \eqref{e:HKE-alpha-stable} and Proposition \ref{p:3P} in the first inequality below, \eqref{e:psi-bound} and \eqref{e:Green-alpha-stable} in the third, and Lemma \ref{l:one-step-transient-key} in the last, we see that
\begin{align}\label{e:one-step-2}
	&\frac{1}{\wt q(t,x,y)}	\int_0^t \int_{B(0,R)^c} q(s,x,z) q(t-s,z,y)\kappa(z)dzds \nn\\
	&\le c_2 	\int_0^t \int_{B(0,R)^c} (\wt q(s,x,z) +  \wt q(t-s,y,z))\kappa(z)dzds\nn\\
		&\le 2c_2 \sup_{w \in B(0,2R)^c}	\int_0^\infty \int_{B(0,R)^c} \wt q(s,w,z) \kappa(z)dzds\nn\\
			&\le c_3 \sup_{w \in B(0,2R)^c} \int_{B(0,R)^c} \frac{dz}{|w-z|^{d-\alpha}|z|^{\beta_1}}\nn\\
			&\le c_4 R^{-(\beta_1-\alpha)}.
\end{align}
Combining \eqref{e:one-step-1} with \eqref{e:one-step-2}, we arrive at \eqref{e:one-step}.  \qed

\subsection{The case of $d=1<\alpha$}

Throughout this subsection, we assume that $d=1<\alpha$. Recall that $R\ge1$ is fixed. Define
\begin{align*}
	h_R(s,z):=  1 \wedge \frac{|z|^{\alpha-1}  R^{(2-\alpha)/2}}{  s^{(\alpha-1)/\alpha} (s\wedge R^\alpha)^{(2-\alpha)/(2\alpha)} } \quad \text{for $s>0$ and $z \in B(0,R)^c$.} 
\end{align*}
By Proposition \ref{p:DHKE-2}(ii) and \eqref{e:one-step-range}, there exists $C=C(d,\alpha)>0$ such that 
\begin{align}\label{e:one-step2-1}
	q^{B(0,R)^c}(t,x,y) \ge C h_R(t,x)h_R(t,y) \wt q(t,x,y) \quad \text{for all} \;\, t \ge R^\alpha \text{ and } x,y \in B(0,2R)^c.
\end{align}
For any $z\in B(0,R)^c$ and $s\le R^\alpha$, it holds that  $h_R(s,z)=1$. Hence, 
\begin{align*}
	h_R(s,z)=  1 \wedge \frac{|z|^{\alpha-1} }{  s^{(\alpha-1)/\alpha}  }, \quad s >0, \; z \in B(0,R)^c.
\end{align*}
Note that, for each fixed $z \in B(0,R)^c$, the map $s\mapsto h_R(s,z)$ is non-increasing and
\begin{align}\label{e:h-scaling}
	s^{(\alpha-1)/\alpha}h_R(s,z) \le u^{(\alpha-1)/\alpha}h_R(u,z) \quad \text{for all} \;\, u\ge s>0.
\end{align}

\begin{lem}\label{l:one-step-strong-recurrent-key1}
	Suppose that $d=1<\alpha$.	There exists $C=C(\alpha,\beta_1,\Lambda)>0$ independent of $R$  such that  for all $t \ge R^\alpha$ and  $v\in B(0,2R)^c$,
	\begin{align*}
		\int_0^{t} \int_{B(0,R)^c} h_R(s,v)h_R(s,z)h_R(t,z)  \kappa(z)dzds\le \frac{C  t^{1/\alpha}h_R(t,v)}{R^{\beta_1-\alpha}}.
	\end{align*} 
\end{lem}
\pf Let $t\ge R^\alpha$ and $v\in R(0,2R)^c$.  Using   \eqref{e:psi-bound} and the fact that  $h_R\le 1$ in the first inequality below, and  \eqref{e:h-scaling} in the second, we obtain
\begin{align*}
		&\int_0^{t} \int_{B(0,R)^c} h_R(s,v)h_R(s,z)h_R(t,z)  \kappa(z)dzds\le \frac{\Lambda^2}{t^{(\alpha-1)/\alpha}} \int_{0}^t \int_{B(0,R)^c} \frac{h_R(s,v)|z|^{\alpha-1}}{|z|^{\beta_1}}dzds\\
	&\le  \Lambda^2  h_R(t,v)\int_{0}^t \int_{B(0,R)^c} \frac{1}{s^{(\alpha-1)/\alpha}  |z|^{\beta_1-\alpha+1}}dzds 
	= c_1 R^{-\beta_1+\alpha} t^{1/\alpha} h_R(t,v).
\end{align*} \qed

\begin{lem}\label{l:one-step-strong-recurrent-key2}

	Suppose that $d=1<\alpha$. There exists $C=C(\alpha,\beta_1,\Lambda)>0$ independent of $R$ such that  for all $t \ge R^\alpha$ and  $v\in B(0,2R)^c$,
	\begin{align*}
			\int_0^{t} \int_{B(0,R)^c} h_R(s,v)h_R(s,z)h_R(t,z) \wt q(s,v,z) \kappa(z)dzds\le   \frac{C  h_R(t,v)}{R^{\beta_1-\alpha}}.
	\end{align*}
\end{lem}
\pf Let $t\ge R^\alpha$ and $v\in R(0,2R)^c$.  Using \eqref{e:psi-bound}, \eqref{e:h-scaling}, and the inequalities $h_R\le1$ and $\wt q(s,\cdot,\cdot)\le s^{-1/\alpha}$, we get
\begin{align*}
	&	\int_0^{t} \int_{B(0,R)^c} h_R(s,v)h_R(s,z)h_R(t,z) \wt q(s,v,z) \kappa(z)dzds\\
	&\le \frac{\Lambda^2}{t^{(\alpha-1)/\alpha}}	\int_0^{t} \int_{B(0,R)^c} \frac{h_R(s,v)
		 \wt q(s,v,z)}{|z|^{\beta_1-\alpha+1}}dzds\\
	& \le 	\Lambda^2	h_R(t,v)\int_0^{t} \int_{B(0,R)^c} \frac{h_R(s,z)\wt q(s,v,z)}{s^{(\alpha-1)/\alpha}|z|^{\beta_1-\alpha+1}} dzds \\
	&\le\Lambda^2	h_R(t,v)\int_0^{R^\alpha} \int_{ B(0,R)^c} \frac{\wt q(s,v,z)}{s^{(\alpha-1)/\alpha}|z|^{\beta_1-\alpha+1}} dzds  + \Lambda^2	h_R(t,v)\int_{R^\alpha}^{t} \int_{B(0,R)^c} \frac{h_R(s,z)}{s|z|^{\beta_1-\alpha+1}} dzds \\
	&=:\Lambda^2 h_R(t,v)(I_1+I_2).
\end{align*}
By using \eqref{e:wtq-integral}, we have
\begin{align*}	I_1&\le \frac{1}{R^{\beta_1-\alpha+1}}\int_0^{R^\alpha} \int_{ B(0,R)^c} \frac{\wt q(s,v,z)}{s^{(\alpha-1)/\alpha}} dzds \le \frac{c_1}{R^{\beta_1-\alpha+1}}\int_0^{R^\alpha}  \frac{ds}{s^{(\alpha-1)/\alpha}} = \frac{\alpha c_1}{R^{\beta_1-\alpha}}.\end{align*}
On the other hand, we observe that
\begin{align*}
		I_2 	&\le \int_{R^\alpha}^t \int_{B(0,s^{1/\alpha})\setminus B(0,R)} \frac{|z|^{\alpha-1}}{s^{1+(\alpha-1)/\alpha} |z|^{\beta_1-\alpha+1}} dzds +\int_{R^\alpha}^t \int_{B(0,s^{1/\alpha})^c} \frac{1}{s|z|^{\beta_1-\alpha+1}} dzds\\	&=:I_{2,1}+I_{2,2}.
\end{align*}
Set $\eta_1:=((\alpha-1) \wedge (\beta_1-\alpha))/2$. For $I_{2,1}$, we have
\begin{align*}
	I_{2,1}
	&\le 
	\int_{R^\alpha}^t \int_{B(0,s^{1/\alpha})\setminus B(0,R)} \frac{(s^{1/\alpha})^{\alpha-1-\eta_1}}{s^{1+(\alpha-1)/\alpha} |z|^{\beta_1-\alpha+1-\eta_1 }} dzds\le \frac{c_2
	}{R^{\beta_1-\alpha-\eta_1}}\int_{R^\alpha}^t \frac{1}{s^{1+\eta_1/\alpha}}ds\le \frac{c_3
	}{R^{\beta_1-\alpha}}.
\end{align*}
For $I_{2,2}$, we see that
\begin{align*}
	I_{2,2} =c_4 
	\int_{R^\alpha}^t  \frac{ds}{s^{1 + (\beta_1-\alpha)/\alpha}} \le \frac{c_5
	}{R^{\beta_1-\alpha}}.
\end{align*}
The proof is complete. \qed 

\begin{lem}\label{l:one-step-strong-recurrent}
	When $d=1<\alpha$, \eqref{e:one-step} holds true.
\end{lem}
\pf 
Let $t\ge R^\alpha$ and $x,y \in B(0,2R)^c$.  
 Using Proposition \ref{p:DHKE-2}(ii)  and Proposition \ref{p:3P} in the first inequality below,  \eqref{e:h-scaling} and the inequality $\wt q(t-s,w,z)\le (t/2)^{-1/\alpha}$ for all $s \in (0,t/2]$ and $w,z \in \R^d$ in the second, and Lemmas \ref{l:one-step-strong-recurrent-key1} and \ref{l:one-step-strong-recurrent-key2} in the last, we obtain
\begin{align*}
	&\frac{1}{\wt q(t,x,y)}	\int_0^t \int_{B(0,R)^c} 	q^{B(0,R)^c}(s,x,z) 	q^{B(0,R)^c}(t-s,z,y)\kappa(z)dzds \nn\\
		&\le c_2 	\int_0^{t} \int_{B(0,R)^c} h_R(s,x)h_R(s,z)h_R(t-s,y)h_R(t-s,z)(\wt q(s,x,z) +  \wt q(t-s,y,z))\kappa(z)dzds  \nn\\
	&\le 	 c_3h_R(t,y)	\int_0^{t/2} \int_{B(0,R)^c} h_R(s,x)h_R(s,z)h_R(t,z) \wt q(s,x,z) \kappa(z)dzds \nn\\
	&\quad  + c_3t^{-1/\alpha}h_R(t,y)	\int_0^{t/2} \int_{B(0,R)^c} h_R(s,x)h_R(s,z)h_R(t,z)  \kappa(z)dzds \nn\\
		&\quad  + c_3h_R(t,x)	\int_0^{t/2} \int_{B(0,R)^c} h_R(s,y)h_R(s,z)h_R(t,z) \wt q(s,y,z) \kappa(z)dzds \nn\\
			&\quad  + c_3t^{-1/\alpha}h_R(t,x)	\int_0^{t/2} \int_{B(0,R)^c} h_R(s,y)h_R(s,z)h_R(t,z)  \kappa(z)dzds \nn\\
			&\le c_4R^{-\beta_1+\alpha} h_R(t,x) h_R(t,y).
 \end{align*}
  Combining this with \eqref{e:one-step2-1}, we conclude that  \eqref{e:one-step} holds when $d=1<\alpha$.
  \qed

\subsection{The case of $d=1=\alpha$}

In this subsection, we assume that $d=1=\alpha$.  Recall that $\Lg r$ is  defined in \eqref{e:Log}.
Note that
\begin{align}\label{e:Log-prop-1}
	\Lg sr \le \Lg s + \Lg r \quad \text{for all} \;\, r,s>0,
\end{align}
\begin{align}\label{e:Log-prop-2}
	\Lg r \ge 1 \quad \text{for all} \;\, r \ge 1
\end{align}
and
\begin{align}\label{e:Log-prop-3}
 \Lg r \le \frac{r\Lg a}{a} \le r \quad \text{for all} \;\, r \ge a \ge 1.
\end{align}
Moreover, for any $\eps>0$, there exists $c(\eps)\ge 1$ such that
\begin{align}\label{e:log-scaling}
	\frac{\Lg u}{\Lg s} \le c(\eps) \bigg(\frac{u}{s}\bigg)^{\eps} \quad \text{for all} \;\, u\ge s>0.
\end{align}
Define for $s>0$ and $z \in B(0,R)^c$,
\begin{align*}
	f_R(s,z):=  1 \wedge \frac{ R^{1/2}\,\Lg(|z|/R)  }{   (s\wedge R)^{1/2} \,	\Lg(s/R) }.
\end{align*} 
By Proposition \ref{p:DHKE-2}(iii), using \eqref{e:one-step-range}, we see that there exists $C=C(d,\alpha)>0$ such that 
\begin{align}\label{e:one-step3-1}
	q^{B(0,R)^c}(t,x,y) \ge C f_R(t,x)f_R(t,y) \wt q(t,x,y) \quad \text{for all} \;\, t \ge R \text{ and } x,y \in B(0,2R)^c.
\end{align}
For any $s\le R$ and  $z \in B(0,R)^c$, $f_R(s,z)=1$. Hence
\begin{align*}
		f_R(s,z)=  1 \wedge \frac{ \Lg(|z|/R)  }{    \Lg(s/R) }.
\end{align*}
Consequently, for each fixed $z \in B(0,R)^c$, we have
\begin{align}\label{e:f-scaling}
f_R(s,z)\,\Lg(s/R)\le    f_R(u,z)\,\Lg(u/R)   \quad \text{for all} \;\, u\ge s >0.
\end{align}

\begin{lem}\label{l:one-step-critical-key1}
Suppose that $d=1=\alpha$. 
	There exists $C=C(\beta_1,\Lambda)>0$ independent of $R$  such that  for all $t \ge R$ and  $v\in B(0,2R)^c$,
	\begin{align*}
		\int_0^{t} \int_{B(0,R)^c} f_R(s,v)f_R(s,z)f_R(t,z)  \kappa(z)dzds\le \frac{C tf_R(t,v)}{R^{\beta_1-1}}.
	\end{align*} 
\end{lem}
\pf Let $t\ge R$ and $v\in R(0,2R)^c$.  Using \eqref{e:psi-bound} and  $f_R(s, z)\le 1$
in the first inequality below, \eqref{e:f-scaling} in the second, and \eqref{e:log-scaling} twice (with $\eps=1/2$ and $\eps=(\beta_1-1)/2$) in the third, we obtain
\begin{align*}
		&\int_0^{t} \int_{B(0,R)^c} f_R(s,v)f_R(s,z)f_R(t,z)  \kappa(z)dzds\\
		 & \le \frac{\Lambda^2}{\Lg(t/R)}	\int_0^{t} \int_{B(0,R)^c} \frac{f_R(s,v)\,\Lg(|z|/R)}{|z|^{\beta_1}}dzds\\
		  & \le \Lambda^2f_R(t,v)	\int_0^{t} \int_{B(0,R)^c} \frac{\Lg(|z|/R)}{|z|^{\beta_1} \,\Lg(s/R)}dzds\\
		    & \le \frac{c_1t^{1/2}f_R(t,v)\,\Lg 1}{R^{(\beta_1-1)/2}\,\Lg (t/R)}	\int_0^{t} \frac{ds}{s^{1/2}} \int_{B(0,R)^c} \frac{dz}{|z|^{(\beta_1+1)/2}}\\
		        & = \frac{c_2tf_R(t,v) }{R^{\beta_1-1}\,\Lg (t/R)} \le  \frac{c_2tf_R(t,v)}{R^{\beta_1-1}}.
\end{align*}
\qed

\begin{lem}\label{l:one-step-critical-key2}
Suppose that $d=1=\alpha$. 
	There exists $C=C(\beta_1,\Lambda)>0$ independent of $R$ such that  for all $t \ge R$ and  $v\in B(0,2R)^c$,
	\begin{align*}
		\int_0^{t} \int_{B(0,R)^c} f_R(s,v)f_R(s,z)f_R(t,z) \wt q(s,v,z) \kappa(z)dzds\le \frac{C f_R(t,v)}{R^{\beta_1-1} }.
	\end{align*}
\end{lem}
\pf Let $t\ge R$ and $v\in R(0,2R)^c$.  Using \eqref{e:psi-bound}, \eqref{e:f-scaling}, and the inequalities $f_R\le1$ and $\wt q(s,\cdot,\cdot)\le s^{-1}$, we get
\begin{align*}
	&	\int_0^{t} \int_{B(0,R)^c} f_R(s,v)f_R(s,z)f_R(t,z) \wt q(s,v,z) \kappa(z)dzds\\
	&\le \frac{\Lambda^2}{\Lg(t/R)}	\int_0^{t} \int_{B(0,R)^c} \frac{f_R(s,v)f_R(s,z) \wt q(s,v,z) \,\Lg(|z|/R)}{|z|^{\beta_1}}dzds\\
	& \le 	\Lambda^2	f_R(t,v)\int_0^{t} \int_{B(0,R)^c} \frac{f_R(s,z)\wt q(s,v,z) \,\Lg(|z|/R)}{|z|^{\beta_1} \Lg(s/R)} dzds \\
	&\le\Lambda^2	f_R(t,v)\int_0^{R} \int_{ B(0,R)^c} \frac{\wt q(s,v,z)\,\Lg(|z|/R)}{|z|^{\beta_1}\,\Lg(s/R)} dzds \\
	&\quad  + \Lambda^2	f_R(t,v)\int_{R}^{t} \int_{B(0,R)^c} \frac{ (\Lg(|z|/R))^2}{s|z|^{\beta_1}(\Lg(s/R))^2} dzds \\
	&=:\Lambda^2 f_R(t,v)(I_1+I_2).
\end{align*}
Using \eqref{e:log-scaling} twice (with $\eps=1/2$ and $\eps=\beta_1$) in the first inequality below and \eqref{e:wtq-integral} in the second, we obtain
\begin{align*}
	I_1&\le \frac{c_1R^{1/2}
		\Lg 1}{R^{\beta_1} \, \Lg 1}	\int_0^{R} \int_{ B(0,R)^c} \frac{\wt q(s,v,z)}{s^{1/2}} dzds\le \frac{c_2 
	 }{R^{\beta_1-1/2}}	\int_0^{R}  \frac{ds}{s^{1/2}} =\frac{2c_2 
	 }{R^{\beta_1-1}}.
\end{align*}
For $I_{2}$, using \eqref{e:log-scaling}  (with $\eps=(\beta_1-1)/4$), we see that
\begin{align*}
	I_{2} 	&\le \frac{c_3 
		(\Lg 1)^2}{R^{(\beta_1-1)/2}}\int_{R}^t  \int_{ B(0,R)^c}\frac{ 1}{ s (\Lg (s/R))^2 \,|z|^{(\beta_1+1)/2}} dzds\nn\\
	&= \frac{c_4 
	}{R^{\beta_1-1}}\int_{R}^t \frac{1}{ s (\Lg (s/R))^2} ds \nn\\
	&\le  \frac{c_4
	}{R^{\beta_1-1}}\int_{R}^\infty \frac{e}{ (s + (e-1)R) (\Lg (s/R))^2} ds=   \frac{ec_4 
	}{R^{\beta_1-1} \,\Lg 1}.
\end{align*}
The proof is complete. \qed 

\begin{lem}\label{l:one-step-critical}
	When $d=1=\alpha$, \eqref{e:one-step} holds true.
\end{lem}
\pf Let $t\ge R$ and $x,y \in B(0,2R)^c$.  
Using Proposition \ref{p:DHKE-2}(iii) and Proposition \ref{p:3P} in the first inequality below,  \eqref{e:f-scaling} and the inequality $\wt q(t-s,w,z)\le (t/2)^{-1}$ for all $s \in (0,t/2]$ and $w,z \in \R^d$ in the second, and Lemmas \ref{l:one-step-critical-key1} and \ref{l:one-step-critical-key2} in the last, we obtain
\begin{align*}
	&\frac{1}{\wt q(t,x,y)}	\int_0^t \int_{B(0,R)^c} 	q^{B(0,R)^c}(s,x,z) 	q^{B(0,R)^c}(t-s,z,y)\kappa(z)dzds \nn\\
	&\le c_2 	\int_0^{t} \int_{B(0,R)^c} f_R(s,x)f_R(s,z)f_R(t-s,y)f_R(t-s,z)(\wt q(s,x,z) +  \wt q(t-s,y,z))\kappa(z)dzds  \nn\\
	&\le 	 c_3f_R(t,y)	\int_0^{t/2} \int_{B(0,R)^c} f_R(s,x)f_R(s,z)f_R(t,z) \wt q(s,x,z) \kappa(z)dzds \nn\\
	&\quad  + c_3t^{-1}f_R(t,y)	\int_0^{t/2} \int_{B(0,R)^c} f_R(s,x)f_R(s,z)f_R(t,z)  \kappa(z)dzds \nn\\
	&\quad  + c_3f_R(t,x)	\int_0^{t/2} \int_{B(0,R)^c} f_R(s,y)f_R(s,z)f_R(t,z) \wt q(s,y,z) \kappa(z)dzds \nn\\
	&\quad  + c_3t^{-1}f_R(t,x)	\int_0^{t/2} \int_{B(0,R)^c} f_R(s,y)f_R(s,z)f_R(t,z)  \kappa(z)dzds \nn\\
	&\le \frac{c_4 f_R(t,x)f_R(t,y)}{R^{\beta_1-1}}.
\end{align*}
Combining this with \eqref{e:one-step3-1}, we conclude that  \eqref{e:one-step} holds.  \qed

Now, the proof of Proposition \ref{p:one-step} is straightforward.  

\medskip

\noindent \textbf{Proof of Proposition \ref{p:one-step}.} The result  follows from Lemmas \ref{l:one-step-transient}, \ref{l:one-step-strong-recurrent} and \ref{l:one-step-critical}. \qed

\section{Large time estimates}\label{s:6}

Throughout this section, we continue to assume that $\kappa \in \sK_\alpha(\psi, \Lambda)$. The goal of this section is to establish the following large time estimates for $p^\kappa(t,x,y)$. 

\begin{thm}\label{t:largetime}
	Suppose that $\kappa \in \sK_\alpha(\psi, \Lambda)$.  Then	there exist comparison constants depending only on  $d,\alpha,\beta_1,\beta_2$ and $\Lambda$  such that the following estimates hold for all $t\ge 2$ and  $x,y \in \R^d_0$:
	\begin{align*}
		\frac{p^\kappa(t,x,y)}{\wt q(t,x,y)} \asymp \begin{cases}
			(1\wedge \psi(|x|)) ( 1\wedge \psi(|y|)) &\mbox{ if } d>\alpha,\\[6pt]
			\left( 1 \wedge \frac{\psi(|x|) \wedge |x|^{\alpha-1} }{t^{(\alpha-1)/\alpha}} \right) \left( 1 \wedge \frac{\psi(|y|) \wedge |y|^{\alpha-1} }{t^{(\alpha-1)/\alpha}} \right) &\mbox{ if } d=1<\alpha,\\[9pt]
			\left( 1 \wedge \frac{\psi(|x|) \wedge \Lg |x| }{\Lg t} \right) \left( 1 \wedge \frac{\psi(|y|) \wedge \Lg |y| }{\Lg t} \right) &\mbox{ if } d=1=\alpha.
		\end{cases}
	\end{align*}
\end{thm}

\subsection{Large time lower heat kernel estimates}

 Recall that for any open subset $U$ of $\R^d_0$ with $\overline U \subset \R^d_0$, 
$$
p^{\kappa,U}(t,x,y)=q^{U}(t,x,y) +  \sum_{k=1}^\infty p^{\kappa,U}_{k}(t,x,y)$$ where  $p^{\kappa,U}_{k}(t,x,y)$, $k\ge 1$, are defined by \eqref{e:p-construction},
and $p^\kappa(t,x,y)$ is defined as the increasing limit of $p^{\kappa,B(0,1/n)^c}(t,x,y)$ as $n \to \infty$.

\begin{lem}\label{l:largetime-lower-bound}
	There exists 
	$R_1=R_1(d,\alpha,\beta_1,\Lambda)\ge 2$  such that for all $t\ge R_1^\alpha$ and  $x,y \in B(0,2R_1)^c$,
	\begin{align*}
		p^\kappa(t,x,y) \ge \frac12q^{B(0,R_1)^c}(t,x,y).
	\end{align*}
\end{lem}
\pf By Proposition \ref{p:one-step}, there exists $R_1=R_1(d,\alpha,\beta_1,\Lambda)\ge 2$ such that for all $t \ge R_1^\alpha$ and $x,y \in B(0,2R_1)^c$,
\begin{align}\label{e:largetime-lower-bound}
	& \int_0^t \int_{B(0,R_1)^c} \frac{ q^{B(0,R_1)^c}(s,x,z) q^{B(0,R_1)^c}(t-s,z,y)}{q^{B(0,R_1)^c}(t,x,y) }\kappa(z)dzds \le \frac13.
\end{align}
Let $t\ge R_1^\alpha$ and  $x,y \in B(0,2R_1)^c$.  By using induction and \eqref{e:largetime-lower-bound},  we see that 
$$|p^{\kappa,B(0,R_1)^c}_{k}(t,x,y)| \le 3^{-k} q^{B(0,R_1)^c}(t,x,y) \quad \text{for all} \;\, k \ge 1.$$ 
 Therefore, we conclude that
\begin{align*}
p^{\kappa}(t,x,y) 
\ge 	p^{\kappa,B(0,R_1)^c}(t,x,y)\ge q^{B(0,R_1)^c}(t,x,y) - q^{B(0,R_1)^c}(t,x,y) \sum_{k=1}^\infty 3^{-k} = \frac12q^{B(0,R_1)^c}(t,x,y).
\end{align*} \qed

\noindent \textbf{Proof of  Theorem \ref{t:largetime} (Lower estimates).}   Let $t\ge 2$ and $x,y \in \R^d_0$. Without loss of generality, we assume that $|x|\le |y|$.  Let  $R_1\ge 2$ be the constant in Lemma \ref{l:largetime-lower-bound} and $z_0 \in \R^d$ be such that $|z_0|=3R_1$. We deal with four cases separately.

\smallskip

Case 1:  $1\le t\le 4R_1^\alpha$  and $|y|< 8R_1$. Note that $|x-y|\le |x|+|y|<16R_1$ and $t \asymp \psi^{-1}(t) \asymp 1$. In particular, $\wt q(t,x,y) \asymp 1$. Applying Theorem \ref{t:smalltime} with $T=4R_1^\alpha$, we obtain
\begin{align}
	p^{\kappa}(t,x,y) &\ge c_1\bigg(1 \wedge \frac{\psi(|x|)}{t}\bigg) \bigg(1 \wedge \frac{\psi(|y|)}{t}\bigg) \frac{t^2}{\psi^{-1}(t)^{d+2\alpha}}\bigg(1 \wedge \frac{\psi^{-1}(t)}{|x-y|}\bigg)^{d+2\alpha}\\
	&\ge c_2(1\wedge \psi(|x|)) ( 1\wedge \psi(|y|)\\
	&\ge c_3(1\wedge \psi(|x|)) ( 1\wedge \psi(|y|) \wt q(t,x,y).
\end{align}
Since $t^{(\alpha-1)/\alpha} \asymp \Lg t \asymp 1$ in this case, this yields the desired lower bound.

\smallskip

Case 2:  $1\le t\le 4R_1^\alpha$ and $|y|\ge 8R_1$. Applying Theorem \ref{t:smalltime} with $T=4R_1^\alpha$, we obtain
\begin{align}
	p^{\kappa}(t,x,y) &\ge c_4\bigg(1 \wedge \frac{\psi(|x|)}{t}\bigg) \bigg(1 \wedge \frac{\psi(|y|)}{t}\bigg)e^{-\lambda_2t/\psi( |y|)}\,\wt q(t,x,y)\\
	&\ge c_4\bigg(1 \wedge \frac{\psi(|x|)}{4R_1^\alpha}\bigg) \bigg(1 \wedge \frac{\psi(|y|)}{4R_1^\alpha}\bigg)e^{-4\lambda_2R_1^\alpha/\psi(8R_1)}\,\wt q(t,x,y) \\
	&\ge c_5(1\wedge \psi(|x|)) ( 1\wedge \psi(|y|))\,\wt q(t,x,y).
\end{align}
Using  $t^{(\alpha-1)/\alpha} \asymp \Lg t \asymp 1$ for $t \in [1,4R_1^\alpha]$, we get the desired lower bound.

\smallskip

Case 3:  $t > 4R_1^\alpha$ and $|y| <8R_1$.  For all $z \in B(z_0,R_1)$, we have $|z| >2R_1$ and 
$$|x-z| \vee |y-z| \le |y| + |z| <12R_1.$$ Hence, by Theorem \ref{t:smalltime}, we get that  for any $z \in B(z_0,R_1)$, 
\begin{align}\label{e:largetime-LHK-1}
	p^\kappa(1,x,z)&\ge c_6(1 \wedge \psi(|x|)) (1 \wedge \psi(2R_1))e^{-\lambda_2/\psi(2R_1)}\left( 1 \wedge 
	(12R_1)^{-1}\right)^{d+\alpha} \nn\\
	&\ge c_7(1 \wedge \psi(|x|))  \ge c_8 \psi(|x|)
\end{align}
and $p^\kappa(1,z,y)\ge c_8\psi(|y|).$ 
Besides,  for all $z,w \in B(z_0,R_1)$, since $t-2 \ge t/2 \ge R_1^\alpha$, $|z| \wedge |w| >2R_1$ and $|z-w|<2R_1 \le 2(t-2)^{1/\alpha}$, from Lemma \ref{l:largetime-lower-bound} and Proposition \ref{p:DHKE-2}, we deduce that
\begin{align*}
	p^\kappa(t-2,z,w) &\ge \frac12 q^{B(0,R_1)^c}(t-2,z,w) \\
	&\ge   c_9 t^{-d/\alpha}\times  \begin{cases}
		1 &\mbox{ if } d>\alpha,\\[2pt]
		1\wedge \left(  (2R_1)^{\alpha-1} / t^{(\alpha-1)/\alpha}\right)^2 &\mbox{ if } d=1<\alpha,\\[2pt]
		1\wedge \left( \Lg (2R_1) / \Lg t\right)^2 &\mbox{ if } d=1=\alpha
	\end{cases}\\
	&\ge   c_{10} t^{-d/\alpha}\times  \begin{cases}
		1 &\mbox{ if } d>\alpha,\\[2pt]
		t^{-2(\alpha-1)/\alpha} &\mbox{ if } d=1<\alpha,\\[2pt]
		(\Lg t)^{-2} &\mbox{ if } d=1=\alpha.
	\end{cases}
\end{align*}
Therefore, by using the semigroup property, we arrive at
\begin{align*}
	p^\kappa(t,x,y) &\ge \int_{B(z_0,R_1)} \int_{B(z_0,R_1)} p^\kappa(1,x,z) p^\kappa(t-2,z,w) p^\kappa(1,w,y) dzdw\\
	&\ge  c_8^2 c_{10} \psi(|x|) \psi(|y|) t^{-d/\alpha} \int_{B(z_0,R_1)} dzdw \times  \begin{cases}
		1 &\mbox{ if } d>\alpha,\\[2pt]
		t^{-2(\alpha-1)/\alpha} &\mbox{ if } d=1<\alpha,\\[2pt]
		(\Lg t)^{-2} &\mbox{ if } d=1=\alpha.
	\end{cases}\\
	&\ge   c_{11} 	\psi(|x|)\psi(|y|) \wt q(t,x,y) \times  \begin{cases}
		1&\mbox{ if } d>\alpha,\\[2pt]
		t^{-2(\alpha-1)/\alpha} &\mbox{ if } d=1<\alpha,\\[2pt]
		(\Lg t)^{-2} &\mbox{ if } d=1=\alpha.
	\end{cases}
\end{align*}

Case 4:  $t > 4R_1^\alpha$ and $|y| \ge 8R_1$. If $|x| \ge 2R_1$, then the lower bound follows from Lemma \ref{l:largetime-lower-bound} and Proposition \ref{p:DHKE-2}. Suppose that $|x|<2R_1$. Note that \eqref{e:largetime-LHK-1} is still valid. Further, for any $z \in B(z_0,R_1)$, since $|z|<4R_1< 2t^{1/\alpha}$, by  \eqref{e:wtq-compare-2} and \eqref{e:wtq-compare},
\begin{align*}
	\wt q(t-1,z,y)\ge c_{12} \wt q(t,z,y) \ge  c_{13} \wt q(t,0,y) \ge  c_{14}\wt q(t,x,y).
\end{align*}
Thus, by Lemma \ref{l:largetime-lower-bound} and Proposition \ref{p:DHKE-2}, we get that for any $z \in B(z_0,R_1)$,
\begin{align*}
	p^\kappa(t-1, z,y) &\ge \frac12 q^{B(0,R_1)^c}(t-1,z,y) \\
	&\ge  c_{15}\wt q(t-1,z,y)\times  \begin{cases}
		1 &\mbox{ if } d>\alpha,\\[4pt]
		\left( 1\wedge 	\frac{(2R_1)^{\alpha-1}}{t^{(\alpha-1)/\alpha}} \right) \left( 1 \wedge \frac{|y|^{\alpha-1}}{t^{(\alpha-1)/\alpha}} \right) &\mbox{ if } d=1<\alpha,\\[6pt]
		\left( 1 \wedge \frac{ \Lg (2R_1) }{\Lg t} \right) \left( 1 \wedge \frac{\Lg |y| }{\Lg t} \right) &\mbox{ if } d=1=\alpha
	\end{cases}\\[3pt]
	&\ge  c_{16}\wt q(t,x,y)\times  \begin{cases}
		1 &\mbox{ if } d>\alpha,\\[4pt]
		t^{-(\alpha-1)/\alpha}\left( 1 \wedge \frac{|y|^{\alpha-1}}{t^{(\alpha-1)/\alpha}} \right) &\mbox{ if } d=1<\alpha,\\[6pt]
		(\Lg t)^{-1} \left( 1 \wedge \frac{\Lg |y| }{\Lg t} \right) &\mbox{ if } d=1=\alpha.
	\end{cases}
\end{align*}
Combining this with \eqref{e:largetime-LHK-1} and using  the semigroup property, we conclude that
\begin{align*}
	p^\kappa (t,x,y) &\ge \int_{B(z_0,R_1)} p^\kappa(1,x,z) p^\kappa(t-1,z,y) dz\\
	& \ge c_{17}\int_{B(z_0,R_1)}  dz\,\psi(|x|) \wt q(t,x,y) \times  \begin{cases}
		1 &\mbox{ if } d>\alpha,\\[4pt]
		t^{-(\alpha-1)/\alpha}\left( 1 \wedge \frac{|y|^{\alpha-1}}{t^{(\alpha-1)/\alpha}} \right) &\mbox{ if } d=1<\alpha,\\[6pt]
		(\Lg t)^{-1} \left( 1 \wedge \frac{\Lg |y| }{\Lg t} \right) &\mbox{ if } d=1=\alpha.
	\end{cases}
\end{align*}
The proof is complete. \qed

\subsection{Large time upper heat kernel estimates}

As preparation for the proof of the upper estimates in Theorem \ref{t:largetime}, we first prove some lemmas.

From \eqref{e:psi-scale}, we see that $\psi(r) \le \Lambda r^{\beta_1}$ for all $r\in (0,1]$ and $\psi(r) \ge \Lambda^{-1} r^{\beta_1}$ for all $r\in [1,\infty)$. Thus, for each fixed $a \in (0,\beta_1)$ and $R>0$, 
 there exist comparison constants depending on $a$ and  $R$ such that 
\begin{align}\label{e:decay-psi-a}
	\psi(r) \wedge r^{a} \asymp \begin{cases}
		\psi(r) &\mbox{ if } r \le R,\\
		r^{a} &\mbox{ if } r>R.
	\end{cases}
\end{align}

\begin{lem}\label{l:largetime-2-upper}
Suppose that $d=1<\alpha$.  Then there exists $C=C(\alpha,\beta_1,\beta_2,\Lambda)\ge 1$ such that for all $t\ge 1$ and  $x,y \in \R^1_0$,
	\begin{align}\label{e:largetime-2-upper}
		\frac{p^\kappa(t,x,y)}{\wt q(t,x,y)} \le C \bigg( 1 \wedge \frac{\psi(|x|) \wedge |x|^{\alpha-1} }{t^{(\alpha-1)/\alpha}} \bigg).
	\end{align}
\end{lem}
\pf  Let $t \ge 1$ and $x,y \in \R^1_0$.  By Proposition \ref{p:DHKE-2}(ii), we get that for all $s>0$ and $v,w \in \R^1_0$,
\begin{align}\label{e:largetime-2-upper-1}
	&\frac{p^\kappa(s,v,w)}{\wt q(s,v,w)}= \lim_{n \to \infty} 	\frac{p^{\kappa,B(0,1/n)^c} (s,v,w) }{\wt q(s,v,w)}\le \lim_{n \to \infty}    \frac{q^{B(0,1/n)^c} (s,v,w) }{\wt q(s,v,w)}  \nn\\
	& \le c_1 \lim_{n \to \infty} \bigg( 1 \wedge \frac{ (|v|-n^{-1})^{\alpha-1} ( |v| \wedge n^{-1})^{(2-\alpha)/2}}{  s^{(\alpha-1)/\alpha} (s\wedge n^{-\alpha})^{(2-\alpha)/(2\alpha)} }\bigg) = c_1  \bigg( 1 \wedge \frac{ |v|^{\alpha-1}}{  s^{(\alpha-1)/\alpha}}\bigg).
\end{align}
Hence, by \eqref{e:decay-psi-a}, we deduce that  \eqref{e:largetime-2-upper} holds if $|x| > 1$.

 Suppose that $|x|\le 1$. Using the semigroup property in the first line below and Corollary \ref{c:UHK-general} in the second, we obtain
\begin{align*}
	p^\kappa(t,x,y)& \le \bigg(  \int_{B(0,2)} + \int_{B(0,t^{1/\alpha})\setminus B(0,2)}  + \int_{B(0,t^{1/\alpha})^c} \bigg) \, p^\kappa(1/2,x,z) p^\kappa(t- 1/2,z,y) dz \\
	&\le c_2 \psi(|x|)  \bigg(  \int_{B(0,2)} + \int_{B(0,t^{1/\alpha})\setminus B(0,2)}  + \int_{B(0,t^{1/\alpha})^c} \bigg)\,  \wt q(1/2,x,z) p^\kappa(t-1/2,z,y) dz\\
	&=:c_2\psi(|x|)(I_1+I_2+I_3).
\end{align*}
By  \eqref{e:largetime-2-upper-1}, \eqref{e:wtq-compare-2} and  \eqref{e:wtq-compare}, for any $z \in B(0,t^{1/\alpha})$,
\begin{align}\label{e:largetime-2-upper-2}
	p^\kappa(t-1/2,z,y) \le \frac{c_1 |z|^{\alpha-1}\wt q(t-1/2,z,y)}{(t-1/2)^{(\alpha-1)/\alpha}}\le \frac{c_3 |z|^{\alpha-1}\wt q(t,0,y)}{t^{(\alpha-1)/\alpha}} \le \frac{c_4 |z|^{\alpha-1}\wt q(t,x,y)}{t^{(\alpha-1)/\alpha}}.
\end{align}
Using this,  we get
\begin{align*}
	I_1 \le  \frac{2c_4  \wt q(t,x,y)}{t^{(\alpha-1)/\alpha}} \int_{B(0,2)} |z|^{\alpha-1}dz  = \frac{c_5\wt q(t,x,y)}{t^{(\alpha-1)/\alpha}}.
\end{align*} 
For $I_2$, we note that for any $z \in B(0,2)^c$,
\begin{align}\label{e:largetime-2-upper-3}
	\wt q(1/2,x,z) \le \frac{1/2}{|x-z|^{\alpha+1}} \le \frac{2^\alpha}{|z|^{\alpha+1}}.
\end{align} 
Using \eqref{e:largetime-2-upper-2} and  \eqref{e:largetime-2-upper-3}, we obtain
\begin{align*}
	I_2\le \frac{2^\alpha c_4 \wt q(t,x,y)}{t^{(\alpha-1)/\alpha}} \int_{B(0,t^{1/\alpha})\setminus B(0,2)}  \frac{dz}{|z|^2} \le \frac{c_6\wt q(t,x,y)}{t^{(\alpha-1)/\alpha}}.
\end{align*}
For $I_3$,  we consider the cases $|y|\le 2t^{1/\alpha}$ and $|y|\le 2t^{1/\alpha}$ separately. If $|y|\le 2t^{1/\alpha}$, then
using \eqref{e:largetime-2-upper-3}, we see that
\begin{align*}
	I_3\le \frac{2^\alpha}{(t^{1/\alpha})^{\alpha+1}} \int_{B(0,t^{1/\alpha})^c} p^\kappa(t-1/2,z,y)dz \le \frac{2^\alpha }{t^{1+1/\alpha}} \le  \frac{2^{1+2\alpha} \wt q(t,0,y)}{t} \le  \frac{2^{1+2\alpha}  \wt q(t,0,y)}{t^{(\alpha-1)/\alpha}}.
\end{align*}
Suppose that $|y|>2t^{1/\alpha}$. Then  for any $z \in B(0,|y|/2)\setminus B(0,t^{1/\alpha})$, by \eqref{e:HKE-alpha-stable}  and \eqref{e:wtq-compare-2}, we get
\begin{align*}
	p^\kappa(t-1/2,z,y) \le c_7 \wt q (t,z,y) = \frac{c_7t}{|y-z|^{\alpha+1}} \le \frac{2^{\alpha+1}c_7t}{|y|^{\alpha+1}}.
\end{align*}
Using this and  \eqref{e:largetime-2-upper-3}, we get that
\begin{align*}
	I_3&\le  2^\alpha \bigg( \int_{B(0,|y|/2) \setminus B(0,t^{1/\alpha})}   + \int_{B(0,|y|/2)^c} \bigg) \frac{p^\kappa(t-1/2,z,y)}{|z|^{\alpha+1}}dz  \\
	&\le \frac{c_8t}{|y|^{\alpha+1}} \int_{B(0,|y|/2) \setminus B(0,t^{1/\alpha})}  \frac{dz}{|z|^{\alpha+1}}  + \frac{2^\alpha}{(|y|/2)^{\alpha+1}}   \int_{B(0,|y|/2)^c}p^\kappa(t-1/2,z,y) dz\\
	&\le \frac{c_9}{|y|^{\alpha+1}} = \frac{c_9 \wt q(t,0,y)}{t} \le \frac{c_9 \wt q(t,0,y)}{t^{(\alpha-1)/\alpha}}.
\end{align*}
Thus, in both cases, using \eqref{e:wtq-compare}, we deduce that
\begin{align*}
	I_3\le c_{10}t^{-(\alpha-1)/\alpha}  \wt q(t,x,y).
\end{align*}
Now combining the estimates for $I_1$, $I_2$ and $I_3$ above and using \eqref{e:decay-psi-a}, we  conclude that \eqref{e:largetime-2-upper} is also valid in this case. The proof is complete. \qed

For $n \ge 1$, we denote by  $\Lg^n  r:=\Lg \circ \cdots \circ \Lg r$ the $n$-th iterated function of $\Lg$.

\begin{lem}\label{l:criticial-upper-1}
Suppose that $d=1=\alpha$.  Let $n \in \N$. If there exists $a_0\ge 1$ such that
	\begin{align}\label{e:criticial-upper-1-ass}
		\frac{p^\kappa(t,z,y)}{\wt q(t,0,y)} \le 
		\frac{a_0 \psi(|z|)\, \Log^n \,t }{\Lg t}, \quad \mbox{for all } t \ge 1 \mbox{ and  } z,y \in \R^1_0 \mbox{ with } |z|\le 1,
	\end{align}
then there exists $C=C(\Lambda)>0$ such that
	\begin{align}\label{e:criticial-upper-1}
		\frac{p^\kappa(t,x,y)}{\wt q(t,0,y)} \le 
		\frac{C \Lg a_0(\Lg |x| + \Log^{n+1} \,t ) }{\Lg t}, \quad \mbox{for all } t \ge 1 \mbox{ and  } x,y \in \R^1_0 \mbox{ with } |x|\le t.
	\end{align}
\end{lem}
\pf Let $t\ge 1$ and $x,y\in \R^1_0$ with $|x|\le t$. Set
$$
\eps:= \bigg( \frac{(\Lg |x|) \wedge (\Log^n\,t) }{a_0\Log^n\,t}\bigg)^{1/\beta_1} \in (0,1].
$$
 For all $w \in \R^1_0$, using Proposition \ref{p:DHKE-2}(iii) in the first inequality below, \eqref{e:Log-prop-1}  and $\beta_1>\alpha=1$ in the second,   $\Lg r> \log (e-1)$ for all $r>0$ in the third and \eqref{e:log-scaling} in the last,  we obtain
\begin{align}\label{e:criticial-upper-1-Dirichlet}
\frac{ q^{B(0,\eps)^c}(t,x,w)}{\wt q(t,x,w)} &\le \frac{c_1  ( (|x|-\eps)_+ \wedge \eps)^{1/2}\, \Lg (|x|/\eps)}{\eps^{1/2}\Lg (t/\eps)}\nn\\
 &  \le \frac{c_1( \Lg |x| + \Lg (  \Log^n \, t/((\Lg |x|) \wedge (\Log^n\,t) ) + \Lg a_0)  }{\Lg t} \nn\\
 &\le  \frac{ c_1 (\Lg a_0) (\Lg |x|  + \Lg(\Log^{n}\,t / (\log (e-1))))}{(\log(e-1))\,\Lg t}\nn\\
 &\le  \frac{c_2 (\Lg a_0)(\Lg |x|  + \Log^{n+1} \, t)}{\Lg t}.
\end{align}
We deal with the cases $|y|\le 2t$ and $|y|>2t$ separately.

\smallskip

Case 1: $|y|\le 2t$. By \eqref{e:criticial-upper-1-ass}, \eqref{e:psi-scale} (with $\psi(1)=1$) and \eqref{e:log-scaling}, we get that for all $w \in \R^1_0$,
\begin{align*}
	\sup_{s\in (t/2,t],\, z\in B(0,\eps)\setminus \{0\}} p^\kappa (s,z,w)  \le \sup_{s\in (t/2,t]} \frac{a_0 \psi(\eps)\Log^n \,s }{s\,\Lg s}   \le  \frac{c_3a_0 \eps^{\beta_1} \Log^n \,t }{t\,\Lg t} \le  \frac{c_3\Lg |x| }{t\,\Lg t}.
\end{align*}
Combining this with \eqref{e:criticial-upper-1-Dirichlet}
 and  applying Lemma \ref{l:Dirichlet-upper}, since $\Lg a_0 \ge \Lg 1 = 1$ and $t \ge 1$, we deduce that
 \begin{align*}
 	p^\kappa(t,x,y) &\le  \frac{c_2 (\Lg a_0)(\Lg |x|  + \Log^{n+1} \, t)}{\Lg t} + \frac{2c_3\Lg |x| }{t\,\Lg t}\\
 	&\le  \frac{(c_2+2c_3)(\Lg a_0)(\Lg |x|  + \Log^{n+1} \, t)}{\Lg t}.
 \end{align*}

Case 2: $|y|>2t$. By the strong Markov property, we see that for a.e. $w \in B(0,2t)^c$,
\begin{align*}
	p^\kappa(t,x,w) \le p^{\kappa, B(0,\eps)^c}(t,x,w) + \E_x[ p^{\kappa}(t- \tau^\kappa_{B(0,\eps)^c} , X^\kappa_{ \tau^\kappa_{B(0,\eps)^c} },w):\tau^\kappa_{B(0,\eps)^c}<t]=:I_1+I_2.
\end{align*}
For all $w \in B(0,2t)^c$,  we have $|w-x| \ge |w|/2$ so that  $\wt q(t,x,w) \le t|x-w|^{-2} \le 4t |w|^{-2}$. Hence,  by \eqref{e:criticial-upper-1-Dirichlet}, we obtain
\begin{align}\label{e:criticial-upper-1-I1}
	I_1&\le  q^{B(0,\eps)^c}(t,x,w) \le  \frac{4c_2 (\Lg a_0)t(\Lg |x|  + \Log^{n+1} \, t)}{|w|^2\Lg t}.
\end{align}
For $I_2$, we note that
\begin{align}\label{e:criticial-upper-1-I2}
	I_2\le \sup_{s\in (0,t],\, z \in B(0,\eps)\setminus \{0\} }p^\kappa(s,z,w).
\end{align}
For all $s\in (0,1]$ and $z\in B(0,\eps)\setminus\{0\}$,  since $|w| >2t \ge 2$, we get from   \eqref{e:HKE-alpha-stable} that
\begin{align}\label{e:criticial-upper-1-I2-1}
	p^\kappa(s,z,w)\le q(s,z,w) \le c_4 \wt q(s,z,w) \le \frac{c_4s}{|w-z|^2} \le \frac{4c_4}{|w|^2} \le \frac{4c_4t}{|w|^2 \Lg t},
\end{align}
where we used \eqref{e:Log-prop-3} in the last inequality.  Besides, for all $s\in (1,t]$ and $z\in B(0,\eps)\setminus\{0\}$, using \eqref{e:criticial-upper-1-ass} and \eqref{e:psi-scale}, we see that
\begin{align}\label{e:criticial-upper-1-I2-2}
	p^\kappa(s,z,w) \le \frac{a_0\wt q(s,0,w)\psi(|z|) \,\Log^n \,t }{\Lg t}  \le \frac{a_0\Lambda \eps^{\beta_1}t \,\Log^n \,t }{|w|^2\Lg t}  \le  \frac{\Lambda t \, \Lg |x|}{|w|^2\Lg t}.
\end{align}
By \eqref{e:criticial-upper-1-I2}, \eqref{e:criticial-upper-1-I2-1} and \eqref{e:criticial-upper-1-I2-2}, we deduce that
\begin{align*}
I_2 \le \frac{4\Lambda c_4 t\,  \Lg |x|}{|w|^2\Lg t}.
\end{align*}
Combining this with \eqref{e:criticial-upper-1-I1}, we obtain
\begin{align*}
	p^\kappa(t,x,w) \le \frac{c_5(\Lg a_0) t (\Lg |x|  + \Log^{n+1} \, t)}{|w|^2\Lg t} \quad \text{for a.e.} \;\, w \in B(0,2t)^c,
\end{align*}
where $c_5:=4c_2  + 4\Lambda c_4$. By the lower semi-continuity of $p^\kappa$, it follows that
\begin{align*}
&	p^\kappa(t,x,y) \le \liminf_{\delta \to 0} \fint_{B(y,\delta)} p^\kappa(t,x,w)dw\\
& \le \frac{c_5(\Lg a_0) t (\Lg |x|  + \Log^{n+1} \, t)}{|y|^2\Lg t} =  \frac{c_5(\Lg a_0)  \wt q(t,0,y)(\Lg |x|  + \Log^{n+1} \, t)}{\Lg t}.
\end{align*}

Combining Case 1 and Case 2, we conclude that 
 \eqref{e:criticial-upper-1} holds with $C_1:= c_2+2c_3 + c_5$. \qed

\begin{lem}\label{l:criticial-upper-2}
Suppose that $d=1=\alpha$. Let $n \in \N$. If there exists $b_0\ge 1$ such that 
	\begin{align}\label{e:criticial-upper-2-ass}
		\frac{p^\kappa(t,x,y)}{ \wt q(t,0,y)} \le \frac{b_0 (\Lg |x| + \Log^n \,t )}
		{\Lg t}, \quad \mbox{ for all } t \ge 1 \mbox{ and } x,y \in \R^1_0 \mbox{ with } |x|\le t, 
	\end{align}
then there exists $C=C(\beta_1,\beta_2,\Lambda)\ge1$ such that
	\begin{align}\label{e:criticial-upper-2}
		\frac{p^\kappa(t,z,y)}{\wt q(t,0,y)} \le \frac{C_2b_0 \psi(|z|) \,\Log^n \,t }
		{\Lg t}, \quad \mbox{ for all } t\ge 1 \mbox{ and } z,y \in \R^1_0 \mbox{ with } |z|\le 1.
	\end{align} 
\end{lem}
\pf  Let $t\ge 1$ and $z,y \in \R^1_0$ with $|z|\le 1$. If $t\le 10$, then by Corollary \ref{c:UHK-general} and  \eqref{e:wtq-compare}, we get
\begin{align*}
	p^\kappa(t,z,y) \le  c_1 \psi(|z|) \wt q(t,z,y) \le c_2\psi(|z|)\wt q(t,0,y)\le \frac{(c_2\Lg 10)b_0 \psi(|z|) \wt q(t,0,y) \Log^n \, t}{\Lg t}.
\end{align*}
Hence, by taking $C_2$ larger than $c_2\Lg 10$, \eqref{e:criticial-upper-2} holds in this case.

Suppose that $t>10$. By \eqref{e:Log-prop-2} and \eqref{e:Log-prop-3},  we see that  $1\le \Log^n  t \le \Log^{n-1}(t/4) \le t/4$. By the semigroup property,
\begin{align*}
	p^\kappa(t,z,y)&= \bigg( \int_{B(0,\Log^n \, t)} + \int_{B(0,t/2)\setminus B(0,\Log^n \, t)} + \int_{B(0,t/2)^c} \bigg)\, p^\kappa(\Log^n  \,t,z,w)\,p^\kappa(t- \Log^n\, t,w,y)dw\nn\\
	&=:I_1+I_2+I_3. 
\end{align*}
Since $t/2 \le t- \Log^n\, t < t$,  we have $\wt q(t-\Log^n \, t,0,y) \le c_3 \wt q(t,0,y)$ and by \eqref{e:log-scaling}, $\Log (t- \Log^n \, t) \ge c_4 \Lg t$. Thus, we get from   \eqref{e:criticial-upper-2-ass} that for all $w \in B(0,t/2)$,
\begin{equation}\label{e:criticall-upper-2-1}
	p^\kappa(t- \Log^n\, t,w,y) \le \frac{c_3b_0 \wt q(t,0,y)(\Lg |w| + \Log^n \,t )}{c_4\Lg t}.
\end{equation}
By Corollary \ref{c:UHK-general}, we have for all $w \in B(0,\Log^n \, t)$,
\begin{align*}
	p^\kappa(\Log^n \,t,z,w) \le \frac{c_5 \psi(|z|)}{\Log^n \,t}.
\end{align*}
Using this and \eqref{e:criticall-upper-2-1}, we obtain
\begin{align*}
	I_1& \le  \frac{c_6b_0  \psi(|z|)\wt q(t,0,y)}{\Lg t \cdot \Log^n \, t} \int_{B(0,\Log^n \, t)} (\Lg |w| + \Log^n \,t ) dw \\
&\le  \frac{c_6b_0  \psi(|z|) \wt q(t,0,y) (\Lg^{n+1}t + \Log^n \,t )}{\Lg t}  \le  \frac{2c_6b_0 \psi(|z|) \wt q(t,0,y) \Log^n \,t }{\Lg t} ,
\end{align*}
where we used \eqref{e:Log-prop-3} in the last inequality. For  $I_2$, using  Corollary \ref{c:UHK-general} and \eqref{e:wtq-compare}, since $|z|\le 1 \le \Log^n \, t$, we see that for all $w \in B(0,t/2) \setminus B(0, \Log^n t)$, 
\begin{align*}
		p^\kappa(\Log^n \,t,z,w) &\le c_7 \psi(|z|) \wt q(\Log^n \,t,0,w)   \le \frac{c_7  \psi(|z|)\Log^n \,t}{|w|^2}.
\end{align*}
Thus, using \eqref{e:criticall-upper-2-1} and \eqref{e:log-scaling} (with $\eps=1/2$), we  obtain
\begin{align*}
		I_2 &\le  \frac{c_8b_0  \psi(|z|) \wt q(t,0,y)}{\Lg t} \int_{B(0,t/2)\setminus B(0,\Log^n \, t)}  \frac{(\Lg |w| + \Log^n \,t )\Log^n \,t}{|w|^2} dw \\
	&\le  \frac{c_9b_0 \psi(|z|) \wt q(t,0,y) (\Lg^{n+1}t + \Log^n \,t )}{\Lg t \cdot (\Log^{n} \, t)^{1/2}}  \int_{B(0,t/2)\setminus B(0,\Log^n \, t)}  \frac{\Log^n \,t}{|w|^{3/2}} dw \\
	& \le  \frac{c_{10}b_0  \psi(|z|)\wt q(t,0,y) \Log^n \,t }{\Lg t}.
\end{align*}
For $I_3$,  we note that for all $w \in B(0,t/2)^c$,  by  Corollary \ref{c:UHK-general} and \eqref{e:wtq-compare},
\begin{align*}
	p^\kappa(\Log^n \,t,z,w) &\le c_{11} \psi(|z|)\wt q(\Log^n\,t,0,w) \le \frac{c_{11} \psi(|z|)\Log^n \,t }{|w|^2}
\end{align*}
and by \eqref{e:HKE-alpha-stable}, since $t- \Log^n\, t\in [t/2,t]$,
\begin{align*}
		p^\kappa(t- \Log^n\, t,w,y) \le c_{12}\wt q(t- \Log^n\, t,w,y) \le c_{13} \wt q(t,w,y) = c_{13} \wt q(t,y,w).
\end{align*}
Hence,   by \eqref{e:wtq-integral}, if $|y|\le t$, then
\begin{align*}
	I_3 &\le c_{14} \psi(|z|)\Log^n \,t \int_{B(0,t/2)^c} \frac{\wt q(t,y,w)}{|w|^2}dw\\
	&\le  \frac{c_{14} \psi(|z|)\Log^n \,t}{(t/2)^2}  \int_{B(0,t/2)^c}  \wt q(t,y,w)dw\le \frac{c_{15} \psi(|z|)\Log^n \,t}{t^2} 
\end{align*}
and if $|y|>t$, then 
\begin{align*}
	I_3 &\le c_{14} \psi(|z|)\Log^n \,t \int_{B(0,t/2)^c} \frac{\wt q(t,y,w)}{|w|^2}dw\\
	&\le \frac{c_{14} \psi(|z|)t\,\Log^n \,t}{(|y|/2)^2} \int_{B(0,|y|/2)\setminus B(0,t/2)} \frac{dw}{|w|^2} +  \frac{c_{14} \psi(|z|)\Log^n \,t}{(|y|/2)^2} \int_{B(0,|y|/2)^c} \wt q(t,y,w) dw \\
	&\le \frac{c_{16}\psi(|z|) \Log^n \,t}{|y|^2}.
\end{align*}
Therefore, in both cases, using \eqref{e:Log-prop-3}, we get
\begin{align*}
	I_3 \le\frac{ (c_{15}\vee c_{16})\psi(|z|) \wt q(t,0,y)\Log^n \,t}{t} \le\frac{ (c_{15}\vee c_{16})b_0\psi(|z|) \wt q(t,0,y)\Log^n \,t}{\Lg t}. 
\end{align*}
The proof is complete. \qed

\begin{lem}\label{l:largetime-3-upper}
Suppose that $d=1=\alpha$.  Then there exists $C=C(\beta_1,\beta_2,\Lambda)\ge 1$ such that for all $t\ge 1$ and  $x,y \in \R^1_0$,
	\begin{align}\label{e:largetime-3-upper}
		\frac{p^\kappa(t,x,y)}{\wt q(t,x,y)} \le C\bigg( 1 \wedge \frac{\psi(|x|) \wedge \Lg |x| }{\Lg t} \bigg).
	\end{align}
\end{lem}
\pf By Corollary \ref{c:UHK-general} and \eqref{e:wtq-compare}, there exists $c_1=c_1(\beta_1,\beta_2,\Lambda)\ge 1$ such that  for all $s \ge 1$ and $z,w\in \R^1_0$ with $|z|\le 1$,
\begin{align*}
	p^\kappa(s,z,w) \le c_1\psi(|z|) \wt q(s,0,w).
\end{align*}
Hence, \eqref{e:criticial-upper-1-ass} holds with $a_0=c_1$ and $n=1$. Applying Lemma \ref{l:criticial-upper-1}, we get that	 for all $s\ge 1$ and $v,w \in \R^1_0$ with $|v|\le s$,
\begin{align*}
	\frac{p^\kappa(s,v,w)}{ \wt q(s,0,w)} &\le \frac{C_1 (\Lg c_1)(\Lg |v| + \Log^{2} s) }{\Lg s}.
\end{align*}
 Then using Lemma \ref{l:criticial-upper-2}, we see that for all $s \ge 1$ and $z,w\in \R^1_0$ with $|z|\le 1$,
\begin{align*}
	\frac{p^\kappa(s,z,w)}{\wt q(s,0,w)} &\le \frac{C_1 C_2 (\Lg c_1) \psi(|z|) \Log^{2} s }{\Lg s}
\end{align*}
Iterating this procedure, we deduce that the following inequalities hold for all $n \ge 1$ and $s\ge 1$:
\begin{align}\label{e:largetime-3-upper-iteration-1}
	\frac{p^\kappa(s,z,w)}{\wt q(s,0,w)}\le \frac{a_n \psi(|z|) \Log^{n+1} s   }{\Lg s} \quad \text{for all $z,w \in \R^1_0$ with $|z|\le 1$}
\end{align}
and
\begin{align}\label{e:largetime-3-upper-iteration-2}
		\frac{p^\kappa(s,v,w)}{\wt q(s,0,w)}\le \frac{b_n(\Lg |v| + \Log^{n+1} s  ) }{\Lg s} \quad \text{for all $v,w \in \R^1_0$ with $|v|\le s$},
\end{align}
where the sequences $(a_n)_{n \ge 1}$ and  $(b_n)_{n\ge 1}$ are defined by  $b_1:=C_1(\log c_1)$,
$$a_{n}:=C_2b_n \quad \text{and} \quad b_{n+1}:=C_1(\Lg a_{n}) = C_1(\Lg (C_2b_n)), \quad n \ge 1. $$
By \eqref{e:log-scaling} (with $\eps=1/2$), we have for all $r \ge 1$,
$$
 C_1(\Lg (C_2r)) \le c_2C_1(C_2  r )^{1/2} = c_3r^{1/2},
$$
where $c_3:=c_2C_1C_2^{1/2} \ge 1$. Hence, for all $n \ge 1$, if $b_n \ge 4c_3^2$, then $b_{n+1} \le c_3 b_n^{1/2} \le b_n/2$. It follows that $\limsup_{n \to \infty} b_n \le 4c_3^2$.  Taking  $\limsup_{n \to \infty}$ in \eqref{e:largetime-3-upper-iteration-1} and \eqref{e:largetime-3-upper-iteration-2}, since $\lim_{n\to \infty} \Log^{n+1} s = 1$ for all $s \ge 1$ and $\Lg r \ge \log(e-1)$ for all $r>0$, we conclude that  for all $s\ge 1$,
\begin{align}\label{e:largetime-3-upper-conclusion-1}
	\frac{p^\kappa(s,z,w)}{\wt q(s,0,w)}\le \frac{4C_2c_3^2 \psi(|z|)  }{\Lg s} \quad \text{for all $z,w \in \R^1_0$ with $|z|\le 1$}
\end{align}
and
\begin{align}\label{e:largetime-3-upper-conclusion-2}
	\frac{p^\kappa(s,v,w)}{\wt q(s,0,w)}\le \frac{4c_3^2 (1 + (\log(e-1))^{-1})(\Lg |v|) }{\Lg s}  \quad \text{for all $v,w \in \R^1_0$ with $|v|\le s$}.
\end{align}

Now, we prove \eqref{e:largetime-3-upper}. Let $t\ge 1$ and $x,y \in \R^1_0$. If $|x|\le 1$, then using \eqref{e:largetime-3-upper-conclusion-1} and \eqref{e:wtq-compare}, we get that
\begin{align*}
	p^\kappa(t,x,y) \le \frac{4C_2c_3^2 \psi(|x|) \wt q(t,0,y)  }{\Lg t} \le   \frac{c_4\psi(|x|) \wt q(t,x,y)  }{\Lg t}.
\end{align*}
By \eqref{e:psi-scale}  and \eqref{e:log-scaling}, we have $
\psi(|x|) \le \Lambda |x|^{\beta_1} \le c_5 \Lg |x| \le c_5\Lg t$. Thus, \eqref{e:largetime-3-upper} holds when $|x|\le 1$.   If $1<|x|\le t$, then
 using \eqref{e:largetime-3-upper-conclusion-2} and \eqref{e:wtq-compare}, we get that
\begin{align*}
	p^\kappa(t,x,y) \le \frac{c_6 (\Lg |x|)  \wt q(t,0,y)  }{\Lg t} \le   \frac{c_7 (\Lg |x|)  \wt q(t,x,y)  }{\Lg t}.
\end{align*}
In this case, by  \eqref{e:log-scaling} and  \eqref{e:psi-scale}, we have
\begin{align*}
	\Lg |x| \le c_8|x|^{\beta_1} \le c_8\Lambda \psi(|x|) \quad \text{and} \quad \Lg |x| \le \Lg t.
\end{align*}
Thus, \eqref{e:largetime-3-upper} holds. If $|x| \ge t$, then by \eqref{e:psi-scale} and \eqref{e:log-scaling} (with $\eps=\beta_1$),   $$\psi(|x|) \wedge \Lg |x| \ge \psi(t) \wedge \Lg t \ge \Lambda t^{\beta_1} \wedge \Lg t \ge c_9 \Lg t.$$ 
Hence, since $p^\kappa(t,x,y) \le q(t,x,y) \le c_{10}\wt q(t,x,y)$ by \eqref{e:HKE-alpha-stable}, we conclude that  \eqref{e:largetime-3-upper} holds. 
 \qed

\medskip

\noindent \textbf{Proof of Theorem \ref{t:largetime} (Upper estimates).}  When $d>\alpha$,   the upper bound follows  from  Corollary \ref{c:UHK-general}.   When $d=1<\alpha$,  using the semigroup property and symmetry of $p^\kappa$ in the first line below, Lemma \ref{l:largetime-2-upper} in the second and \eqref{e:wtq-semigroup} in the third, we get that  for all $t \ge 2$ and $x,y \in \R^1_0$,
\begin{align*}
	p^\kappa(t,x,y) &= \int_{\R^1_0} 	p^\kappa(t/2,x,z)	p^\kappa(t/2,y,z) dz\\
	&\le c_1\bigg( 1 \wedge \frac{\psi(|x|) \wedge |x|^{\alpha-1} }{(t/2)^{(\alpha-1)/\alpha}} \bigg)\bigg( 1 \wedge \frac{\psi(|y|) \wedge |y|^{\alpha-1} }{(t/2)^{(\alpha-1)/\alpha}} \bigg) \int_{\R^1_0} 	\wt q(t/2,x,z)	\wt q(t/2,y,z) dz\\	&\le c_2\bigg( 1 \wedge \frac{\psi(|x|) \wedge |x|^{\alpha-1} }{t^{(\alpha-1)/\alpha}} \bigg)\bigg( 1 \wedge \frac{\psi(|y|) \wedge |y|^{\alpha-1} }{t^{(\alpha-1)/\alpha}} \bigg) \wt q(t,x,y).
\end{align*}
When $d=\alpha=1$, using the semigroup property and symmetry of $p^\kappa$ in the first line below, Lemma \ref{l:largetime-3-upper} in the second, and \eqref{e:log-scaling} and  \eqref{e:wtq-semigroup} in the third, we deduce  that  for all $t \ge 2$ and $x,y \in \R^1_0$,
\begin{align*}
	p^\kappa(t,x,y) &= \int_{\R^1_0} 	p^\kappa(t/2,x,z)	p^\kappa(t/2,y,z) dz\\
	&\le c_3 \bigg( 1 \wedge \frac{\psi(|x|) \wedge \Lg |x| }{\Lg (t/2)} \bigg)  \bigg( 1 \wedge \frac{\psi(|y|) \wedge \Lg |y| }{\Lg (t/2)} \bigg)  \int_{\R^1_0} 	\wt q(t/2,x,z)	\wt q(t/2,y,z) dz\\	&\le c_4\bigg( 1 \wedge \frac{\psi(|x|) \wedge \Lg |x| }{\Lg t} \bigg)  \bigg( 1 \wedge \frac{\psi(|y|) \wedge \Lg |y| }{\Lg t} \bigg)  \wt q(t,x,y).
\end{align*}
The proof of the theorem is complete.  \qed

\section{Green function estimates}\label{s:7}

 Define
\begin{align*}
G^\kappa(x,y):=\int_0^\infty p^\kappa(t,x,y)dt, \quad x,y \in \R^d_0.
\end{align*} The goal of this section is to prove the following two-sided Green function estimates.

\begin{thm}\label{t:Green}
Suppose that  $\kappa \in \sK_\alpha(\psi,\Lambda)$.  Then  there   exist comparison constants depending only on $d,\alpha,\beta_1,\beta_2$ and $\Lambda$ such that the following estimates hold for all $x,y \in \R^d_0$:

\smallskip

\noindent (i) If $d>\alpha$, then \begin{align}\label{e:Green-1}
 	G^\kappa(x,y) \asymp	\left( 1\wedge \frac{\psi(|x|)\wedge  1}{(|x-y|\wedge 1)^{\alpha}}\right) \left( 1\wedge \frac{\psi(|y|)\wedge 1}{(|x-y|\wedge 1)^{\alpha}}\right) \frac{1}{|x-y|^{d-\alpha}}.
 \end{align}
 \noindent (ii) If $d=1<\alpha$, then
 	\begin{align}\label{e:Green-2}
 	G^\kappa(x,y) &\asymp 	\left( 1\wedge \frac{\psi(|x|)\wedge 1}{(|x-y|\wedge 1)^\alpha}\right) 	\left( 1\wedge \frac{\psi(|y|)\wedge 1}{(|x-y|\wedge 1)^\alpha}\right) 	\left( 1\wedge \frac{|x| \vee 1}{|x-y| \vee 1}\right)^{\alpha-1} 	\left( 1\wedge \frac{|y| \vee 1}{|x-y| \vee 1}\right)^{\alpha-1} \nn\\
 	&\;\;\quad \times \big( (\psi(|x|) \wedge  |x|^\alpha) \vee |x-y|^\alpha \big)^{(\alpha-1)/(2\alpha)}  \big( (\psi(|y|) \wedge  |y|^\alpha) \vee |x-y|^\alpha \big)^{(\alpha-1)/(2\alpha)}.
 \end{align}
 \noindent (iii) If $d=1=\alpha$, then
	\begin{align}\label{e:Green-3}
		G^\kappa(x,y) &\asymp 	\left( 1\wedge \frac{\psi(|x|)\wedge 1}{|x-y|\wedge 1}\right) 	\left( 1\wedge \frac{\psi(|y|)\wedge 1}{|x-y|\wedge 1}\right) 	\left( 1\wedge \frac{\Lg |x|}{\Lg |x-y|}\right)^{1/2} \left( 1\wedge \frac{\Lg |y|}{\Lg |x-y|}\right)^{1/2}  \nn\\
		&\;\;\quad \times \left[   \Log \left( \frac{\psi(|x|) \wedge |x| }{|x-y|\wedge 1}\right)  \Log \left( \frac{\psi(|y|) \wedge |y| }{|x-y| \wedge 1}\right) \right]^{1/2}.
	\end{align}
\end{thm}

The following lemma will be used in the proof of Theorem \ref{t:Green}.

\begin{lem}\label{l:Green}
Suppose that  $\kappa \in \sK_\alpha(\psi,\Lambda)$.  Then  there   exist comparison constants depending only on $d,\alpha,\beta_1,\beta_2$ and $\Lambda$ such that the following estimates hold for all $x,y \in \R^d_0$ with $|x-y| \vee \psi(|x| \wedge|y|)\le 1$:

\smallskip 
\noindent (i) If $|x-y|^\alpha > \psi(|x| \wedge |y|)$, then
\begin{align*}
	G^\kappa(x,y) \asymp  \frac{\psi(|x|) \psi(|y|)}{|x-y|^{d+\alpha}}.
\end{align*}
\noindent (ii) If $|x-y|^\alpha \le \psi(|x| \wedge |y|)$, then
\begin{align}\label{e:Green-small-case2}
G^\kappa(x,y) \asymp \int_{|x-y|^\alpha}^{2\psi (|x|\wedge |y|)}   t^{-d/\alpha}dt  \asymp \begin{cases}
|x-y|^{-d+\alpha} &\mbox{ if } d>\alpha,\\[3pt]
 \psi(|x| \wedge |y|)^{(\alpha-1)/\alpha} &\mbox{ if } d=1<\alpha,\\[3pt]
\Log \big(\psi(|x| \wedge |y|)/|x-y|\big) &\mbox{ if } d=1=\alpha.
\end{cases}
\end{align} 
\end{lem}
\pf Let $x,y \in \R^d_0$ be such that $|x-y| \vee \psi(|x| \wedge|y|) \le 1$. Without loss of generality, we assume that $|x|\le |y|$.

\smallskip

(i) Suppose that $|x-y|^\alpha>\psi(|x|)$.  By  \eqref{e:psi-scale}, since $\beta_1>\alpha$, we have
\begin{align}\label{e:Green-small-1}
	&\psi(|y|) \le \psi(|x|+|x-y|) \le \psi(2|x|) \vee \psi(2|x-y|) \nn\\
	&\le 2^{\beta_2}\Lambda (\psi(|x|) \vee \psi(|x-y|)) \le 2^{\beta_2}\Lambda (\psi(|x|) \vee (\Lambda|x-y|^{\beta_1})) \le 2^{\beta_2}\Lambda^2  |x-y|^\alpha .
\end{align}
Applying Theorem \ref{t:smalltime} with $T=2^{\beta_2+1}\Lambda^2$, we get
\begin{align*}
	G^\kappa(x,y)&\ge c_1\int_{ \psi(|y|)}^{2\psi(|y|)}  \bigg(1 \wedge \frac{\psi(|x|)}{t}\bigg) \bigg(1 \wedge \frac{\psi(|y|)}{t}\bigg)  e^{-\lambda_2t/\psi( |y|)}t^{-d/\alpha}\bigg( 1 \wedge \frac{t^{1/\alpha}}{|x-y|}\bigg)^{d+\alpha} dt\\	
	&\ge c_2\psi(|x|)\psi(|y|)\int_{ \psi(|y|)}^{2\psi(|y|)} t^{-2-d/\alpha}\bigg( \frac{t^{1/\alpha}}{|x-y|}\bigg)^{d+\alpha} dt\\	
	& = \frac{c_2\psi(|x|) \psi(|y|)}{|x-y|^{d+\alpha}} \int_{ \psi(|y|)}^{2\psi(|y|)} \frac{dt}{t} =  \frac{(\log 2)c_2\psi(|x|) \psi(|y|)}{|x-y|^{d+\alpha}}.
\end{align*}
On the other hand, by using Theorems \ref{t:smalltime} and \ref{t:largetime},  we have
\begin{align*}
	G^\kappa(x,y) &\le  c_3 \int_0^{\psi(|x|)}  \bigg( \frac{t}{|x-y|^{d+\alpha}}   +\frac{t^2}{|x-y|^{d+2\alpha}}\bigg)  dt  \\
	&\quad + c_3 \psi(|x|) \int_{\psi(|x|)}^{\psi(|y|)}  \bigg( \frac{1}{|x-y|^{d+\alpha}}   +\frac{t}{|x-y|^{d+2\alpha}}\bigg)  dt   \\
	&\quad + c_3\psi(|x|)\psi(|y|) \int_{\psi(|y|)}^{2^{\beta_2}\Lambda^2 |x-y|^\alpha}   \bigg( \frac{e^{-\lambda_1t/\psi( |y|)}}{t |x-y|^{d+\alpha} }   + \frac{1}{|x-y|^{d+2\alpha}} \bigg) dt \\
	&\quad  + c_3\psi(|x|)\psi(|y|) \int_{2^{\beta_2}\Lambda^2 |x-y|^\alpha}^{2^{\beta_2}\Lambda^2 }   \bigg( \frac{e^{-\lambda_1t/\psi( |y|)}}{t^{2+d/\alpha} } + \frac{1}{\psi^{-1}(t)^{d+2\alpha}}\bigg)  dt\\
	&\quad +  c_3 \psi(|x|)\psi(|y|) \int_{2^{\beta_2}\Lambda^2 }^\infty \frac{dt}{t^{d/\alpha}}\\
	&=:c_3(I_1+I_2+I_3+I_4+I_5).
\end{align*}
Since $\psi(|x|)<|x-y|^\alpha$ and $|x|\le |y|$, we have
\begin{align*}
	I_1\le  \frac{\psi(|x|)^2}{2|x-y|^{d+\alpha}}   +\frac{\psi(|x|)^3}{3|x-y|^{d+2\alpha}} \le  \frac{5\psi(|x|)^2}{6|x-y|^{d+\alpha}}  \le  \frac{5\psi(|x|)\psi(|y|)}{6|x-y|^{d+\alpha}} .
\end{align*}
For $I_2$, using  \eqref{e:Green-small-1}, we obtain
\begin{align*}
	I_2 \le \frac{\psi(|x|)\psi(|y|)}{|x-y|^{d+\alpha}}   +\frac{\psi(|x|)\psi(|y|)^2}{2|x-y|^{d+2\alpha}}   \le  \frac{(2^{\beta_2}\Lambda^2+1)\psi(|x|)\psi(|y|)}{|x-y|^{d+\alpha}} .
\end{align*}
For $I_3$, using    the inequality $e^{-r}\le 1/r$ for $r>0$, we see that
\begin{align*}
	I_3\le \frac{\psi(|x|)\psi(|y|)^2}{\lambda_1|x-y|^{d+\alpha}}   \int_{\psi(|y|)}^{\infty} \frac{dt}{t^2} +  \frac{\psi(|x|)\psi(|y|)}{|x-y|^{d+\alpha}} = \frac{(\lambda_1^{-1}+1)\psi(|x|)\psi(|y|)}{|x-y|^{d+\alpha}}.
\end{align*} 
For $I_4$, using \eqref{e:psi-inverse}, we get
\begin{align*}
	I_4 \le \psi(|x|)\psi(|y|)\int_{2^{\beta_2}\Lambda^2 |x-y|^\alpha}^{2^{\beta_2}\Lambda^2 }   \left( \frac{1}{t^{2+d/\alpha} } +  \frac{c_4}{t^{2+d/\alpha}}\right)  dt\le   \frac{c_5\psi(|x|)\psi(|y|)}{|x-y|^{d+\alpha}}.
\end{align*}
For $I_5$, since $|x-y|\le 1$, we have
\begin{align*}
	I_5\le c_6 \psi(|x|)\psi(|y|)\le   \frac{c_6\psi(|x|)\psi(|y|)}{|x-y|^{d+\alpha}}.
\end{align*}
The proof of (i) is complete.

\smallskip

(ii) The second comparison in \eqref{e:Green-small-case2} is straightforward. We now  prove  the first comparison.   Applying Theorem \ref{t:smalltime} with $T=2$, we get
\begin{align*}
	G^\kappa(x,y)&\ge c_1\int_{|x-y|^\alpha}^{2\psi(|x|)}  \bigg(1 \wedge \frac{\psi(|x|)}{t}\bigg) \bigg(1 \wedge \frac{\psi(|y|)}{t}\bigg)  e^{-\lambda_2t/\psi( |y|)}t^{-d/\alpha}\bigg( 1 \wedge \frac{t^{1/\alpha}}{|x-y|}\bigg)^{d+\alpha} dt\\	
	& \ge  2^{-2}c_1e^{-2\lambda_2}\int_{|x-y|^\alpha}^{2\psi(|x|)} t^{-d/\alpha} dt.
\end{align*}
On the other hand, using Corollary \ref{c:UHK-general} and  Theorem \ref{t:largetime},  we see that 
\begin{align*}
	G^\kappa(x,y) &\le  c_2 \int_0^{|x-y|^\alpha}  \frac{t}{|x-y|^{d+\alpha}}   dt  + c_2  \int_{|x-y|^\alpha}^{2\psi(|x|)} \frac{dt}{t^{d/\alpha}}    + c_2\psi(|x|)\int_{2\psi(|x|)}^{2}   \frac{dt}{t^{1+d/\alpha}} \\
	&\quad + c_2\1_{\{d=1<\alpha\}}\psi(|x|)\psi(|y|)\int_{2}^{\infty} \frac{dt}{t^{1+(\alpha-1)/\alpha}}  + c_2\1_{\{d=1=\alpha\}}\psi(|x|)\psi(|y|)\int_{2}^{\infty} \frac{dt}{t (\Lg t)^2}  \\
	&=: c_2(I_1 + I_2+I_3 +\1_{\{d=1<\alpha\}} I_4 +\1_{\{d=1=\alpha\}} I_5).
\end{align*}
Observe that
\begin{align*}
	I_2 \ge \int_{|x-y|^\alpha}^{2|x-y|^\alpha} \frac{dt}{t^{d/\alpha}}  \ge \frac{1}{2^{d/\alpha}|x-y|^{d-\alpha}}
\end{align*}
and
\begin{align}\label{e:Green-small-I2}
	I_2 \ge \int_{\psi(|x|)}^{2\psi(|x|)} \frac{dt}{t^{d/\alpha}}  \ge \frac{1}{2^{d/\alpha}\psi(|x|)^{d/\alpha-1}}.
\end{align}
Hence,  
\begin{align*} 
	I_1 =2^{-1}|x-y|^{\alpha-d}    \le 2^{d/\alpha-1} I_2 \qquad \text{and} \qquad  I_3 \le (\alpha/d) \psi(|x|)^{1-d/\alpha} \le 2^{d/\alpha} (\alpha/d) I_2.
\end{align*}
By  using \eqref{e:psi-scale} and the fact that $\psi(r) \le \Lambda r^{\beta_1} \le \Lambda r^\alpha$ for all $r \in (0,1]$, we have
\begin{equation}\label{e:Green-small-2}
	\psi(|x|)\le 	\psi(|y|) \le  2^{\beta_2}\Lambda (\psi(|x|) \vee \psi(|x-y|))  \le  2^{\beta_2}\Lambda (\psi(|x|) \vee (\Lambda|x-y|^{\alpha})) \le 2^{\beta_2}\Lambda^2 \psi(|x|).
\end{equation}
When $d=1<\alpha$, using  \eqref{e:Green-small-2} and \eqref{e:Green-small-I2}, since $\psi(|x|)\le 1$, we get that
\begin{align*}
	I_4 = c_3 \psi(|x|)\psi(|y|) \le 2^{\beta_2}\Lambda^2 c_3 \psi(|x|)^2 \le 2^{\beta_2}\Lambda^2 c_3 \psi(|x|)^{1-1/\alpha} \le c_4 I_2.
\end{align*}
Similarly, when $d=1=\alpha$,  using  \eqref{e:Green-small-2}, \eqref{e:Green-small-I2} and $\psi(|x|)\le 1$, we obtain
\begin{align*}
	I_5 = c_5\psi(|x|)\psi(|y|) \le 2^{\beta_2}\Lambda^2 c_5 \psi(|x|)^2  \le 2^{\beta_2}\Lambda^2 c_5  \le c_6 I_2.
\end{align*}
The proof is complete. \qed

\noindent \textbf{Proof of Theorem \ref{t:Green}.} Let $x,y \in \R^d_0$. Without loss of generality, we assume that $|x|\le |y|$.

\smallskip

(i) Suppose that $d>\alpha$.  By \eqref{e:HKE-alpha-stable} and \eqref{e:Green-alpha-stable}, we have
\begin{align}\label{e:Green-upper-trivial}
	G^\kappa(x,y) \le \int_0^\infty q(t,x,y)dt \le c_1 \int_0^\infty \wt q(t,x,y)dt = \frac{c_2}{|x-y|^{d-\alpha}}.
\end{align}
We deal with four cases separately.

\smallskip

Case 1: $|x-y| \vee \psi(|x|) \le 1$. By \eqref{e:Green-small-1}, we have  $\psi(|y|) \le 2^{\beta_2} \Lambda^2|x-y|^\alpha\le 2^{\beta_2} \Lambda^2$ in this case. Moreover, by Lemma \ref{l:Green}, it holds that 
\begin{align*}
	G^\kappa(x,y) \asymp	\left( 1\wedge \frac{\psi(|x|)}{|x-y|^{\alpha}}\right) \left( 1\wedge \frac{\psi(|y|)}{|x-y|^{\alpha}}\right) \frac{1}{|x-y|^{d-\alpha}}.
\end{align*}
Hence, \eqref{e:Green-1} holds true. 

\smallskip

Case 2: $|x-y|\le 1 <\psi(|x|)$.   Applying Theorem \ref{t:smalltime} (with $T=2$),  we get
\begin{align*}
	G^\kappa(x,y)&\ge c_1\int_{|x-y|^\alpha}^{2|x-y|^\alpha}  \bigg(1 \wedge \frac{\psi(|x|)}{t}\bigg) \bigg(1 \wedge \frac{\psi(|y|)}{t}\bigg)  e^{-\lambda_2t/\psi( |y|)}t^{-d/\alpha}\bigg( 1 \wedge \frac{t^{1/\alpha}}{|x-y|}\bigg)^{d+\alpha} dt\\	
	& \ge  2^{-2}c_1e^{-2\lambda_2}\int_{|x-y|^\alpha}^{2|x-y|^\alpha} t^{-d/\alpha} dt = \frac{c_2}{|x-y|^{d-\alpha}}.
\end{align*}
Combining this with \eqref{e:Green-upper-trivial}, we get \eqref{e:Green-1}.

\smallskip

Case 3: $|x-y|>1$ and $|x|\ge 1/2$.   By Theorem \ref{t:largetime},
\begin{align*}
	G^\kappa(x,y)  \ge c_{3} (1 \wedge \psi(1/2))^2\int_{2|x-y|^\alpha}^{3|x-y|^\alpha} \frac{t}{|x-y|^{d+\alpha}}dt = \frac{c_{4}}{|x-y|^{d-\alpha}}.
\end{align*}
Combining this with \eqref{e:Green-upper-trivial},  we get \eqref{e:Green-1}.

\smallskip

Case 4: $|x-y|>1$ and $|x|<1/2$. In this case, we have $|y| \ge |y-x|-|x|>1/2$.  Applying Theorem \ref{t:largetime}, we get that
\begin{align*}
		G^\kappa(x,y)  \ge c_{5} (1 \wedge \psi(|x|))(1 \wedge \psi(1/2))\int_{2|x-y|^\alpha}^{3|x-y|^\alpha} \frac{t}{|x-y|^{d+\alpha}}dt \ge \frac{c_{6}\psi(|x|)}{|x-y|^{d-\alpha}}.
\end{align*}
On the other hand, by Corollary \ref{c:UHK-general},  Theorem \ref{t:largetime} and \eqref{e:Green-alpha-stable}, since $|x-y|>1$, we have
\begin{align*}
	G^\kappa(x,y) &\le c_{7} \int_0^2 \frac{\psi(|x|)}{|x-y|^{d+\alpha}} dt  + c_{7} \psi(|x|) \int_2^\infty \wt q(t,x,y)dt\\
	& \le \frac{2c_{7}\psi(|x|)}{|x-y|^{d+\alpha}}  + \frac{c_{8}\psi(|x|)}{|x-y|^{d-\alpha}} \le \frac{(2c_{7}+c_{7})\psi(|x|)}{|x-y|^{d-\alpha}} .
\end{align*}
The proof of (i) is complete. 

\smallskip

(ii) Suppose that $d=1<\alpha$. We deal with  five cases separately.

\smallskip

 Case 1: $\psi(|x|)<|x-y|^\alpha\le 1$.  Using \eqref{e:decay-psi-a} and  \eqref{e:Green-small-1}, we see that  the right-hand side of \eqref{e:Green-2} is comparable to
 \begin{align*}
 &\left( \frac{\psi(|x|)}{|x-y|^\alpha}\right) 	\left(  \frac{\psi(|y|)}{|x-y|^\alpha}\right) 	 (  |x-y|^\alpha )^{(\alpha-1)/(2\alpha)}  (  |x-y|^\alpha )^{(\alpha-1)/(2\alpha)} = \frac{\psi(|x|)\psi(|y|)}{|x-y|^{1+\alpha}}.
 \end{align*}
By Lemma \ref{l:Green}(i),  \eqref{e:Green-2} holds true.
 
 \smallskip
 
 Case 2:  $ |x-y|^\alpha \le \psi(|x|)\le 1$.   By using  \eqref{e:decay-psi-a} and  \eqref{e:Green-small-2}, we see that  the right-hand side of \eqref{e:Green-2} is comparable to	$\psi(|x|)^{(\alpha-1)/(2\alpha)}  \psi(|y|)^{(\alpha-1)/(2\alpha)} \asymp \psi(|x|)^{(\alpha-1)/\alpha}$. 
 By Lemma \ref{l:Green}(ii), we get  \eqref{e:Green-2}.
 
 \smallskip

Case 3: $|x-y|\le 1\le  |x|$. Note that $|x|\le |y|\le |x|+|x-y|\le 2|x|$. Hence, using \eqref{e:decay-psi-a}, we see that  the right-hand side of \eqref{e:Green-2} is comparable to $( |x|^\alpha)^{(\alpha-1)/(2\alpha)}  ( |y|^\alpha)^{(\alpha-1)/(2\alpha)} \asymp |x|^{\alpha-1}$ 
in this case. For the lower bound, using Theorem \ref{t:largetime} and \eqref{e:decay-psi-a},  we obtain
\begin{align*}
	G^\kappa(x,y)&\ge c_1 \int_{2|x|^\alpha}^{3|x|^\alpha} \bigg( 1 \wedge \frac{|x|^{\alpha-1} }{t^{(\alpha-1)/\alpha}} \bigg)^2  \frac{dt}{t^{1/\alpha}}  \ge c_2\int_{2|x|^\alpha}^{3|x|^\alpha}  \frac{dt}{t^{1/\alpha}} = c_3 |x|^{\alpha-1}.
\end{align*}
For the upper bound,  using Corollary \ref{c:UHK-general}, Theorem \ref{t:largetime} and the inequality $|x|\le |y|\le 2|x|$, we get that 
\begin{align*}
	G^\kappa(x,y) &\le  c_{4} \int_0^{|x-y|^\alpha}  \frac{t}{|x-y|^{1+\alpha}}   dt  + c_{4}  \int_{|x-y|^\alpha}^{2|x|^\alpha} \frac{dt}{t^{1/\alpha}}    + c_{4}  |x|^{\alpha-1} \int_{2|x|^\alpha}^{2|y|^\alpha} \frac{dt}{t}  \\
	&\quad + c_{4}  |x|^{\alpha-1}|y|^{\alpha-1} \int_{2|y|^\alpha}^{\infty} \frac{dt}{t^{1+(\alpha-1)/\alpha}}  \\
	&\le c_{5}  \left(|x-y|^{\alpha-1} + |x|^{\alpha-1}+ |x|^{\alpha-1} \log(|y|^\alpha/|x|^\alpha) + |x|^{\alpha-1}  \right) \\
	&\le (3+\log 2^\alpha) c_{5} |x|^{\alpha-1}.
\end{align*}

Case 4: $|x-y|>1$ and $1/2\le |x| <  |x-y|$. Using \eqref{e:decay-psi-a}, we see that  the right-hand side of \eqref{e:Green-2} is comparable to
	\begin{align*}
	\left( \frac{|x|}{|x-y|}\right)^{\alpha-1} 	\left( \frac{|y|}{|x-y|}\right)^{\alpha-1} \big(|x-y|^\alpha \big)^{(\alpha-1)/(2\alpha)}  \big(  |x-y|^\alpha \big)^{(\alpha-1)/(2\alpha)} = \frac{|x|^{\alpha-1}|y|^{\alpha-1} }{|x-y|^{\alpha-1}}.
\end{align*}
By Theorem \ref{t:largetime}, we have
\begin{align*}
	G^\kappa(x,y) & \ge c_{6}  \int_{2|x-y|^\alpha}^{3|x-y|^\alpha}  \bigg(1 \wedge \frac{|x|^{\alpha-1}}{t^{(\alpha-1)/\alpha}}\bigg) \bigg(1 \wedge \frac{|y|^{\alpha-1}}{t^{(\alpha-1)/\alpha}}\bigg)  t^{-1/\alpha}\bigg( 1 \wedge \frac{t^{1/\alpha}}{|x-y|}\bigg)^{1+\alpha}dt\\
	&\ge   c_{7} |x|^{\alpha-1}|y|^{\alpha-1} \int_{2|x-y|^\alpha}^{3|x-y|^\alpha}  \frac{1}{t^{1-2/\alpha}|x-y|^{1+\alpha}}dt = \frac{ c_{8} |x|^{\alpha-1}|y|^{\alpha-1}}{|x-y|^{\alpha-1}}.
\end{align*}
For the upper bound,  by 
Corollary \ref{c:UHK-general} and  Theorem \ref{t:largetime}, since
 $|x|\le |y| \le |x|+|x-y|<2|x-y|$, we have
\begin{align*}
	G^\kappa(x,y) &\le  c_{9} \int_0^{8|x|^\alpha}  \frac{t}{|x-y|^{1+\alpha}}   dt  + c_{9} |x|^{\alpha-1}  \int_{8|x|^\alpha}^{8|y|^\alpha} \frac{t^{1/\alpha}}{|x-y|^{1+\alpha}} dt    \\
	&\quad  + c_{9}  |x|^{\alpha-1}  |y|^{\alpha-1} \int_{8|y|^\alpha}^{8|x-y|^\alpha} \frac{t^{-1+2/\alpha}}{|x-y|^{1+\alpha}} dt+ c_{9}  |x|^{\alpha-1}|y|^{\alpha-1} \int_{8|x-y|^\alpha}^{\infty} \frac{dt}{t^{1+(\alpha-1)/\alpha}}  \\
	&\le c_{10}  \left( \frac{|x|^{2\alpha}}{|x-y|^{\alpha+1}} + \frac{|x|^{\alpha-1} |y|^{\alpha+1} }{|x-y|^{\alpha+1}}+ \frac{|x|^{\alpha-1} |y|^{\alpha-1}}{|x-y|^{\alpha-1}}  + \frac{|x|^{\alpha-1}  |y|^{\alpha-1}}{ |x-y|^{\alpha-1}} \right) \\
	&\le\frac{c_{11}|x|^{\alpha-1}  |y|^{\alpha-1}}{ |x-y|^{\alpha-1}} .
\end{align*}

Case 5: $|x-y|>1$ and $|x|<1/2$. In this case, we have $|y| \le |x|+|x-y| \le (3/2)|x-y|$ and $|y| \ge |x-y|-|x| \ge (1/2)|x-y|$. In particular, $|y| >1/2$ so that $\psi(|y|)\wedge |y|^\alpha \asymp |y|^\alpha$ and $\psi(|y|)\wedge |y|^{\alpha-1} \asymp |y|^{\alpha-1} $ by \eqref{e:decay-psi-a}. Therefore,
 the right-hand side of \eqref{e:Green-2} is comparable to 
 	\begin{align*}
  \psi(|x|) |x-y|^{1-\alpha}  ( |x-y|^\alpha )^{(\alpha-1)/(2\alpha)}  (  |x-y|^\alpha )^{(\alpha-1)/(2\alpha)} = \psi(|x|).
 \end{align*}
 Applying Theorem \ref{t:largetime}, and using \eqref{e:decay-psi-a} and $|y| \asymp |x-y|$, we get that
\begin{align*}
	G^\kappa(x,y)  &\ge c_{12} \int_{2|x-y|^\alpha}^{3|x-y|^\alpha}  \bigg(1 \wedge \frac{\psi(|x|)}{t^{(\alpha-1)/\alpha}}\bigg) \bigg(1 \wedge \frac{|y|^{\alpha-1}}{t^{(\alpha-1)/\alpha}}\bigg)  t^{-1/\alpha}\bigg( 1 \wedge \frac{t^{1/\alpha}}{|x-y|}\bigg)^{1+\alpha}dt\\
	&\ge c_{13}\psi(|x|) \int_{2|x-y|^\alpha}^{3|x-y|^\alpha}  t^{-1}dt = (\log (3/2))c_{13} \psi(|x|).
\end{align*}
On the other hand, using Corollary \ref{c:UHK-general} and  Theorem \ref{t:largetime}, since $|y|\asymp |x-y|>1$, we obtain
\begin{align*}
	G^\kappa(x,y) &\le  c_{14} \psi(|x|)\int_0^{2|x-y|^\alpha}  \frac{1}{|x-y|^{1+\alpha}}   dt  + c_{14} \psi(|x|) |y|^{\alpha-1} \int_{2|x-y|^\alpha}^{\infty} \frac{1}{t^{1+(\alpha-1)/\alpha}} dt    \\
	&\le c_{15}  \left(  \psi(|x|) |x-y|^{-1} + \psi(|x|) |y|^{\alpha-1} |x-y|^{-\alpha+1}\right)\le c_{16}\psi(|x|).
\end{align*}
The proof of (ii) is complete. 

\smallskip

(iii) Suppose that $d=1=\alpha$. We deal with five cases separately.

\smallskip

Case 1:  $\psi(|x|)< |x-y|\le 1$.   Note that $|x| \le \psi^{-1}(1)=1$ and $|y| \le |x| + |x-y|\le 2$ in this case. Using  \eqref{e:decay-psi-a}, \eqref{e:Green-small-1}  and the fact that  $\Lg r \asymp 1$ for $r \in (0,2]$, we see that the  right-hand side of \eqref{e:Green-3} is comparable to  $\psi(|x|)\psi(|y|) |x-y|^{-2}$ in this case. Hence,  by Lemma \ref{l:Green}(i),  \eqref{e:Green-3} holds.

\smallskip

Case 2: $ |x-y| \le \psi(|x|)\le 1$. We have  $|x| \le 1$ and $|y| \le |x| + |x-y|\le 2$. By \eqref{e:decay-psi-a}, \eqref{e:Green-small-2}, \eqref{e:log-scaling}  and the fact that  $\Lg r \asymp 1$ for $r \in (0,2]$,  the right-hand side of \eqref{e:Green-3} is  comparable to
	\begin{align*}
\left[   \Log \left( \frac{\psi(|x|) }{|x-y|}\right)  \Log \left( \frac{\psi(|y|)  }{|x-y|}\right) \right]^{1/2} \asymp \Log \left( \frac{\psi(|x|) }{|x-y|}\right)
\end{align*}
 in this case. Hence, \eqref{e:Green-3} follows from Lemma \ref{l:Green}(ii).
 
 \smallskip

Case 3: $|x-y| \vee 1 \le |x|$. We have $|x|\le |y|\le 2|x|$. Using this, \eqref{e:decay-psi-a} and \eqref{e:log-scaling}, we see that the right-hand side of \eqref{e:Green-3} is comparable to
	\begin{align*}
\left[   \Log \left( \frac{|x| }{|x-y|\wedge 1}\right)  \Log \left( \frac{ |y| }{|x-y| \wedge 1}\right) \right]^{1/2} \asymp  \Log \left( \frac{|x| }{|x-y|\wedge 1}\right) 
\end{align*}
in this case. For the lower bound, by using Theorems \ref{t:smalltime} and \ref{t:largetime}, $\psi(r) \wedge \Lg r \asymp \Lg r$ for $r\ge 1$  and \eqref{e:Log-prop-1}, since $|y|\ge |x|\ge 1$,  we have
\begin{align*}
&	G^\kappa(x,y)\\
&\ge c_1 \1_{\{|x-y|\le 1\}}\int_{|x-y|}^{  2} \bigg(1 \wedge \frac{|x|}{t}\bigg) \bigg(1 \wedge \frac{|y|}{t}\bigg)  e^{-\lambda_2 t/ |y|}t^{-1}\bigg( 1 \wedge \frac{t}{|x-y|}\bigg)^{2} dt  + c_1 (\Lg |x|)^2 \int_{2|x|}^{\infty}  \frac{dt}{t(\Lg t)^2} \\
	&\ge c_2 \1_{\{|x-y|\le 1\}}\int_{|x-y|}^{ 2}  t^{-1} dt  + c_2 (\Lg |x|)^2 \int_{2|x|}^{\infty}  \frac{dt}{t(\log t)^2}\\
	& = c_2 \1_{\{|x-y|\le 1\}}\log\left(\frac{2}{|x-y|}\right) + \frac{c_2(\Lg |x|)^2}{\log(2|x|)} \\
	&\ge c_3 \1_{\{|x-y|\le 1\}} \Log\left(\frac{1}{|x-y|}\right)  + c_3 \Log |x| \ge c_3 \Log\left(\frac{|x|}{|x-y| \wedge 1}\right).
\end{align*}
For the upper bound,  using Corollary \ref{c:UHK-general}, Theorem \ref{t:largetime}, $|x| \asymp |y|$ and \eqref{e:log-scaling},   we see that 
\begin{align*}
	&G^\kappa(x,y)\\
	 &\le  c_{4} \int_0^{|x-y|}  \frac{t}{|x-y|^{2}}   dt  + c_{4}  \int_{|x-y|}^{2|x|} \frac{dt}{t}    + c_{4} \Lg |x| \int_{2|x|}^{2|y|} \frac{dt}{t}  + c_{4} (\Lg |x|)(\Lg |y|) \int_{2|y|}^{\infty} \frac{dt}{t (\Lg t)^2}  \nn\\
	&\le c_{5}  \left(1 +  \Log\left(\frac{|x|}{|x-y|}\right) + (\log 2)  \Lg|x|  + \Lg |x|  \right) \le c_{6} \Log\left(\frac{|x|}{|x-y| \wedge 1}\right).
\end{align*}

Case 4: $|x-y|>1$ and $1/2\le |x| <  |x-y|$. We see that the right-hand side of \eqref{e:Green-3} is comparable to 
\begin{align*}
		\left(  \frac{\Lg |x|}{\Lg |x-y|}\right)^{1/2} \left(\frac{\Lg |y|}{\Lg |x-y|}\right)^{1/2} \left[    (\Lg |x|)(  \Lg |y|) \right]^{1/2} = \frac{ (\Lg |x|)(  \Lg |y|) }{\Lg |x-y|} .
\end{align*}
Since $|y| \le |x|+|x-y|< 2|x-y|$, by Theorem \ref{t:largetime}, we have
\begin{align*}
	G^\kappa(x,y) & \ge c_{7}(\Lg |x|) (\Lg |y|)  \int_{2|x-y|}^{\infty}  \frac{1}{t(\Lg t)^2}dt \ge \frac{c_{8}(\Lg |x|) (\Lg |y|) }{\Lg |x-y|}. 
\end{align*}
For the upper bound, by  using Corollary \ref{c:UHK-general} and  Theorem \ref{t:largetime},  we obtain 
\begin{align*}
	&G^\kappa(x,y) \le  c_{9} \int_0^{4|x|}  \frac{t}{|x-y|^{2}}   dt  + \frac{c_{9} (\Lg |x|)}{|x-y|^2} \int_{4|x|}^{4|y|} \frac{t}{\Lg t} dt    \\
	&\;\; \qquad  \qquad  + \frac{c_{9} (\Lg |x|)(\Lg |y|)}{|x-y|^2} \int_{4|y|}^{4|x-y|} \frac{t}{(\Lg t)^2} dt  + c_{9}   (\Lg |x|)(\Lg |y|)  \int_{4|x-y|}^{\infty} \frac{dt}{t(\Lg t)^2}  \\
	&\le c_{10}  \left( \frac{|x|^{2}}{|x-y|^{2}} + \frac{(\Lg |x|)|y|^2}{|x-y|^{2} \Lg |y|}+ \frac{(\Lg |x|)(\Lg |y|)}{(\Lg|x-y|)^{2}}  +  \frac{(\Lg |x|)|y|^2}{|x-y|^{2} \Lg |y|}+ \frac{(\Lg |x|)(\Lg |y|)}{\Lg|x-y|} \right) \\
	&\le \frac{c_{11}(\Lg |x|)(\Lg |y|)}{\Lg|x-y|} .
\end{align*}

Case 5: $|x-y|>1$ and $|x|<1/2$. In this case, we have $(1/2)|x-y|\le |y| \le (3/2)|x-y|$. Hence, using \eqref{e:decay-psi-a}, \eqref{e:log-scaling} and the fact that $\Lg r \asymp 1 $ for $r \in (0,1]$,  we see that  the right-hand side of \eqref{e:Green-3} is comparable to 
\begin{align*}
\psi(|x|)	\left( \frac{\Lg |x|}{\Lg |x-y|}\right)^{1/2} \left[   (\Lg \psi(|x|))( \Lg |y|) \right]^{1/2} \asymp  \frac{\psi(|x|)(\Lg |y|)^{1/2}}{(\Lg |x-y|)^{1/2}}\asymp \psi(|x|).
\end{align*}
By Theorem \ref{t:largetime}, since $\psi(|x|) \wedge \Lg |x| \asymp \psi(|x|)$ and $|y| \asymp |x-y|>1$, we get that
\begin{align*}
	G^\kappa(x,y)  &\ge c_{12}\int_{2|x-y|}^{3|x-y|}  \bigg(1 \wedge \frac{\psi(|x|)}{\Lg t}\bigg) \bigg(1 \wedge \frac{\Lg |y|}{\Lg t}\bigg)  t^{-1}\bigg( 1 \wedge \frac{t}{|x-y|}\bigg)^{2}dt\\
	&\ge c_{13}\psi(|x|) \int_{2|x-y|}^{3|x-y|}  t^{-1}dt = (\log (3/2))c_{13} \psi(|x|).
\end{align*}
Besides, using Corollary \ref{c:UHK-general},  Theorem \ref{t:largetime} and \eqref{e:log-scaling}, since $|x-y| >1$,  we obtain
\begin{align*}
	G^\kappa(x,y) &\le  c_{14} \psi(|x|)\int_0^{1}  \frac{1}{|x-y|^{2}}   dt  +c_{14} \psi(|x|)\int_1^{2|x-y|}  \frac{t}{|x-y|^{2}(\Lg t)}   dt \\
	&\quad +  c_{14} \psi(|x|) (\Lg |y|) \int_{2|x-y|}^{\infty} \frac{1}{t(\Lg t)^2} dt    \\
	&\le c_{15}\psi(|x|)  \left( |x-y|^{-2} + (\Lg |x-y|)^{-1}+ (\Lg |y|)/ (\Lg |x-y|)\right) \\
	&\le c_{16}\psi(|x|).
\end{align*}

The proof is complete.  \qed

 \section{Appendix}\label{s:8}
 
 In this appendix, we prove the following general result used in getting large time heat kernel upper bound in the case
 $d=1=\alpha$. We believe it is of independent interest.
 
 \begin{lem}\label{l:Dirichlet-upper}
 	Let $t>0$ and $\eps>0$. 	If there exists a non-negative continuous function 
	$F_{t,\eps}$ on $\R^d_0$ such that
 	\begin{align}\label{e:Dirichlet-upper-cond}
 		\sup_{s\in (t/2,t],\, z \in B(0,\eps)\setminus \{0\}}p^{\kappa}(s, z,y) \le F_{t,\eps}(y) \quad \text{for all} \;\, y \in \R^d_0,
 	\end{align}
 then we have
 	\begin{align*}
 		p^\kappa(t,x,y) \le p^{\kappa,B(0,\eps)^c}(t,x,y) + F_{t,\eps}(x) + F_{t,\eps}(y)  \quad \text{for all} \;\, x,y \in \R^d_0.
 	\end{align*}
 \end{lem} 
 \pf 
 Let $x,y \in \R^d_0$ and $\delta>0$. Using the strong Markov property,  we see that for all $w,  u \in \R^d_0$,
 \begin{align*}
 	&	\int_{B(u,\delta)} p^\kappa(t/2,w,v)dv = \P_w(X^\kappa_{t/2} \in B(u,\delta)) \\
 	&= \P_w(X^{\kappa,B(0,\eps)^c}_{t/2} \in B(u,\delta)) +  \P_w(X^{\kappa}_{t/2} \in B(u,\delta), \; \tau^\kappa_{B(0,\eps)^c}<t/2)\\
 	&= \P_w(X^{\kappa,B(0,\eps)^c}_{t/2} \in B(u,\delta)) +  \E_w\left[ \P_{X^\kappa_{\tau^\kappa_{B(0,\eps)^c} }} ( X^{\kappa}_{t/2-\tau^\kappa_{B(0,\eps)^c}}  \in B(u,\delta) ) \, : \tau^\kappa_{B(0,\eps)^c}<t/2\right] \\
 	&\le  \int_{B(u,\delta)}p^{\kappa,B(0,\eps)^c}(t/2,w,v)dv + \sup_{s\in (0,t/2], \, z \in B(0,\eps)\setminus \{0\}}  \int_{B(u,\delta)} p^{\kappa}(s,z,v)dv.
 \end{align*}
 By the Lebesgue differentiation theorem,  for a.e. $(w,u) \in \R^d_0 \times \R^d_0$,
 \begin{align}\label{e:Dirichlet-upper-1}
 	p^\kappa(t/2,w,u) \le p^{\kappa,B(0,\eps)^c}(t/2,w,u) + \sup_{s\in (0,t/2], \, z \in B(0,\eps)\setminus \{0\}}  p^{\kappa}(s,z,u).
 \end{align}
 Using the semigroup property of $X^\kappa$ in the  equality below, \eqref{e:Dirichlet-upper-1} in the first  inequality, the symmetry of $p^\kappa$ and    \eqref{e:Dirichlet-upper-1} in the second,  the symmetry and the semigroup properties of $p^\kappa$ and $p^{\kappa,B(0,\eps)^c}$ in the third, and \eqref{e:Dirichlet-upper-cond} in the fourth, we obtain
 \begin{align*}
 	&	\fint_{B(y,\delta)} p^\kappa(t,x,v)dv =  \int_{\R^d_0}  p^\kappa(t/2,x,u) \fint_{B(y,\delta)} p^\kappa(t/2,u,v) dv du \\
 	&\le  \int_{\R^d_0}   p^{\kappa,B(0,\eps)^c}(t/2,x,u) \fint_{B(y,\delta)} p^\kappa(t/2,u,v) dv du  \\
 	&\quad + \sup_{s\in (0,t/2], \, z \in B(0,\eps)\setminus \{0\}}  \int_{\R^d_0}  p^{\kappa}(s,z,u) \fint_{B(y,\delta)} p^{\kappa}(t/2,u,v)dv du\\
 	&\le  \int_{\R^d_0}   p^{\kappa,B(0,\eps)^c}(t/2,x,u) \fint_{B(y,\delta)} p^{\kappa, B(0,\eps)^c}(t/2,v,u) dv du  \\
 	&\quad + \sup_{s\in (0,t/2], \, z \in B(0,\eps)\setminus \{0\}}  \int_{\R^d_0}  p^{\kappa,B(0,\eps)^c}(t/2,x,u) \fint_{B(y,\delta)} p^{\kappa}(s,z,u)dv du\\
 	&\quad + \sup_{s\in (0,t/2], \, z \in B(0,\eps)\setminus \{0\}}  \int_{\R^d_0}  p^{\kappa}(s,z,u) \fint_{B(y,\delta)} p^{\kappa}(t/2,u,v)dv du\\
 	&\le  \fint_{B(y,\delta)} p^{\kappa, B(0,\eps)^c}(t,x,v) dv   \\
 	&\quad + \sup_{s\in (0,t/2], \, z \in B(0,\eps)\setminus \{0\}}  \fint_{B(y,\delta)} p^{\kappa}(t/2 +s,z,x)  dv \\
 	&\quad + \sup_{s\in (0,t/2], \, z \in B(0,\eps)\setminus \{0\}}   \fint_{B(y,\delta)} p^{\kappa}(t/2 + s,z,v)dv  \\
 	&\le  \fint_{B(y,\delta)} p^{\kappa, B(0,\eps)^c}(t,x,v) dv    + F_{t,\eps}(x)  \fint_{B(y,\delta)}  dv  + \fint_{B(y,\delta)} F_{t,\eps}(v)  dv.
 \end{align*}
 Since $p^\kappa(t,x,\cdot)$ is lower semi-continuous on $\R^d_0$, we have $	p^\kappa(t,x,y) \le \liminf_{\delta\to 0 }\fint_{B(y,\delta)} p^\kappa(t,x,v)dv.$ Therefore, by the continuities of $p^{\kappa,B(0,\eps)^c}(t,x,\cdot)$ and $F_{t,\eps}$, we get the desired result. \qed 
 
 \bigskip
\noindent
{\bf Acknowledgements:} We thank Panki Kim for helpful comments. Part of the research for this paper was done while the second-named author was visiting Jiangsu Normal University, where he was partially supported by a grant from the National Natural Science Foundation of China (11931004, Yingchao Xie).

	\vskip 0.4truein

\noindent {\bf Soobin Cho:} Department of Mathematics,
University of Illinois Urbana-Champaign,
Urbana, IL 61801, U.S.A.
Email: \texttt{soobinc@illinois.edu}
	
	\medskip

\noindent {\bf Renming Song:} Department of Mathematics,
University of Illinois Urbana-Champaign,
Urbana, IL 61801, U.S.A.
Email: \texttt{rsong@illinois.edu}

\end{document}